\documentclass[preprint, 10pt]{elsarticle}

\usepackage[margin=2.5cm]{geometry}
\usepackage{setspace}
\usepackage{lmodern}
\usepackage{amsmath,amssymb}
\usepackage{graphicx}
\usepackage{subfigure}
\usepackage{epsf,epsfig}
\usepackage{bm,bbm}
\usepackage{algpseudocode}
\usepackage{algorithmicx, algorithm}
\usepackage{physics}
\usepackage{color}
\usepackage{pxfonts}
\usepackage[footnotesize,bf]{caption}
\usepackage{mathtools}
\usepackage{hhline}
\usepackage[T1]{fontenc}
\usepackage{placeins}

\usepackage{multirow,dcolumn}
\usepackage{mathrsfs}
\usepackage{verbatim}
\usepackage{comment}
\usepackage{epstopdf}
\usepackage{lineno}
\usepackage[colorlinks=true,
linkcolor=blue,
anchorcolor=blue,
citecolor=blue,
urlcolor=blue,]{hyperref}
\usepackage{ulem}

\definecolor{myyellow}{rgb}{0.75, 0.5, 0.25} 

\usepackage{amsfonts}

\usepackage{enumitem}
\setlist[enumerate]{leftmargin=.5in}
\setlist[itemize]{leftmargin=.5in}

\usepackage{amsopn}

\newcommand{\bsigma}{\boldsymbol{\sigma}}
\newcommand{\bepsilon}{\boldsymbol{\epsilon}}
\newcommand{\xx}{{\mathbf x}}
\newcommand{\bA}{{\mathbf A}}
\newcommand{\BB}{{\mathbf B}}
\newcommand{\FF}{{\mathbf F}}
\newcommand{\RR}{{\mathbf R}}
\newcommand{\MM}{{\mathbf M}}
\newcommand{\uu}{{\mathbf u}}
\newcommand{\vv}{{\mathbf v}}
\newcommand{\II}{{\mathbf I}}
\newcommand{\ff}{{\mathbf f}}
\newcommand{\nn}{{\mathbf n}}
\newcommand{\bPhi}{\boldsymbol{\Phi}}
\newcommand{\ds}{\displaystyle}
\newcommand{\ww}{{\mathbf w}}
\newcommand{\bpsi}{{\mathbf \psi}}

\usepackage{stmaryrd}
\usepackage{mathrsfs}
\usepackage{subfigure}
\usepackage{booktabs}
\usepackage{multirow}
\usepackage{amsopn}
\usepackage{verbatim}

\newtheorem{mythm}{Theorem}[section]

\newtheorem{mylem}{Lemma}[section]
\newtheorem{myassum}{Assumption}[section]
\newtheorem{myinfer}{Corollary}[section]
\newtheorem{mypf}{Proof}[section]

\DeclareMathOperator{\essinf}{essinf}
 \DeclareMathOperator{\esssup}{esssup}

\newtheorem{remark}{Remark}[section]

\graphicspath{{figures/}}
\journal{XXX}

\begin{document}
\begin{frontmatter}
    \title{Regularized coupling multiscale method for thermomechanical coupled problems}
    
    \author[TJU,tjmoe,has]{Xiaofei Guan}
    \ead{guanxf@tongji.edu.cn}
    \author[TJU]{Lijian Jiang}
    \ead{ljjiang@tongji.edu.cn}
    \author[TJU]{Yajun Wang\corref{cor}}
    \ead{1910733@tongji.edu.cn}

    \address[TJU]{School of Mathematical Sciences, Tongji University, Shanghai 200092, China}
    \address[tjmoe]{Key Laboratory of Intelligent Computing and Applications (Tongji University), Ministry of Education, Shanghai 20092, China}
    \address[has]{Institute of Mathematics, Henan Academy of Sciences, Zhengzhou 450046, China}
    \cortext[cor]{Corresponding Author}

    \begin{abstract}
    The coupling effects in multiphysics processes are often neglected in designing multiscale methods. The coupling  may be described by a non-positive definite operator, which in turn brings  significant challenges in  multiscale simulations. In the paper, we develop a regularized coupling multiscale method based on the generalized multiscale finite element method (GMsFEM) to solve coupled thermomechanical problems, and it is referred to as the coupling generalized multiscale finite element method (CGMsFEM). The method consists of defining the coupling multiscale basis functions through local regularized coupling spectral problems in each coarse-grid block, which can be implemented by a novel design of two relaxation parameters. Compared to the standard GMsFEM, the proposed method can not only accurately capture the multiscale coupling correlation effects of multiphysics problems but also greatly improve computational efficiency with fewer multiscale basis functions. In addition, the convergence analysis is also established, and the optimal error estimates are derived, where the upper bound of errors is independent of the magnitude of the relaxation coefficient. Several numerical examples for  periodic, random microstructure, and random material coefficients are presented to validate the theoretical analysis. The numerical results show that the CGMsFEM shows better robustness and efficiency than uncoupled  GMsFEM. 
    \end{abstract}


\end{frontmatter}

\section{Introduction} \label{section:intro}
Heterogeneous media have been extensively applied in practical engineering, and their physical and mechanical properties across multiple temporal and spatial scales closely interact with their underlying microstructure, properties, compositions, etc. Hence, efficient coupling multiscale modeling is of paramount importance for the performance prediction \cite{DYGC:JCP:2022}, optimization design \cite{CS:AM:2012, PP:SMO:2010}, and safety assessment \cite{DCAK:JMS:2018} of heterogeneous media. In the case of thermomechanical processes, except for microscale heterogeneity, microscale thermal expansion properties may also play a significant role at multiple spatial scales. For example, depending on the complex loading conditions (e.g., high stress or extreme temperature environments), damage or cracks may occur due to induced sharp thermal stresses \cite{NFW:JAE:2012}, which can be computationally prohibitive. For this reason, efforts should be made to develop more efficient multiscale models and numerical methods to predict the coupled multiphysics behavior of heterogeneous media for engineering practice.

Based on the first law of thermodynamics and the conservation of momentum, mathematical modeling of thermomechanical processes can lead to coupled partial differential equations (PDEs) \cite{D:ARMA:1968}, including hyperbolic mechanical equations and parabolic heat transfer equations that are mutually coupled. Considering the quasi-static case, the existence and uniqueness of the solutions for the coupled PDEs have been analyzed in \cite{PM:CPDE:1992,Z:RAN:1984}, where the inertia terms of the mechanical parts are ignored. Numerical methods have also been developed with the help of the finite element method (FEM), and the convergence analysis and error orders have been obtained in \cite{CB:CS:1989, EM:EMMNA-MMAN:2009}. Numerous studies have focused on the 2D and 3D linear thermomechanical problems \cite{ S:CMAME:1999,E:ZFAMP:2022} and nonlinear thermomechanical problems \cite{YVM:CNSNS:2022}, including coupling effects \cite{EFC:JTS:2017}, fracture analysis \cite{D:IJNME:2008}, etc. It is imperative to emphasize that although the physical mechanisms underlying the thermomechanical equations and the poroelastic equations diverge, their mathematical equivalence prevails. Numerous methodologies and conclusions demonstrate uniformity across both equation types. Consequently, this paper abstains from an extensive exploration of the poroelastic equations. Furthermore, the multiscale nature of heterogeneous materials leads to a tremendous cost when solving thermomechanical problems using direct numerical methods, which are also unreliable and notoriously ill-conditioned \cite{BDRSSS:CPC:2020}. These motivate us to develop robust and efficient multiscale methods that can not only reduce the dimensionality of coupled thermomechanical problems but also incorporate as much important microscopic physical information as possible.

Research on multiscale methods is of great significance in modeling and numerically computing heterogeneous media. This remains an active field in engineering practice because it has the ability to fully describe the intrinsic physical processes of complex systems. Great progress has been made to solve multiscale problems using various strategies, which include constrained macro simulation and the generalized  finite element method (GFEM) \cite{BO:SIAMJNA:1983}. Constrained macro simulation is a micro-macro coupling interaction method that establishes a suitable macroscopic model through solving specific microscopic problems, including the homogenization method \cite{BLP:ENH:1978,A:SJMA:1992}, the variational multiscale method \cite{HFMQ:CMAME:1998}, the heterogeneous multiscale method \cite{EELRV:CICP:2007}, the computational homogenization method \cite{OBG:CMAME:2008}, etc. The generalized finite element method was originally proposed in \cite{BO:SIAMJNA:1983} and has been applied to solve various multiscale problems by constructing different finite element spaces. In the multiscale finite element method (MsFEM) \cite{HE:S:2009}, the basis functions of the finite element space are constructed by solving the local problems on each coarse block  incorporated with the localized microscale physical information, and the corresponding mathematical theories are also developed for the MsFEM, such as convergence analysis and the error reduction method. Then the generalized multiscale finite element method (GMsFEM) \cite{EGH:JCP:2013a,CCJ:LSMGMs:2016} is designed to obtain high-accuracy multiscale basis functions by solving local spectral problems on snapshot spaces. Similarly, the multiscale spectral generalized finite element method (MS-GFEM) \cite{BLSS:CMAME:2020,MSD:SJNA:2022} is another important multiscale GFEM with local approximation spaces constructed by solving local spectral problems. Moreover, the localized orthogonal decomposition (LOD) method \cite{HM:SJSC:2014,MP:MC:2014} stems from the ideas of  the variational multiscale method \cite{HFMQ:CMAME:1998}, and the coarse space is modified with the basis functions well  approximated locally. These techniques have been extended and applied to the thermomechanical and poroelastic problems, for instance, the homogenization method \cite{FDCG:IJCM:2014,F:SJMA:1983,WCW:MMS:2015,GYT:IJNME:2016}, LOD \cite{MP:EMMNA:2017,altmann2020computational}, GMsFEM \cite{brown2016generalized,VS:NAIA:2017} and extended MsFEM \cite{ZYZZ:CS:2013,ZZLZ:IJMMD:2020}, etc.

Inspired by previous works, this paper designed an efficient regularized coupling formulation of GMsFEM named CGMsFEM. For the traditional GMsFEM \cite{brown2016generalized,VS:NAIA:2017} or the LOD method \cite{MP:EMMNA:2017,altmann2020computational}, the multiscal basis functions are separately constructed through the corresponding elliptic operators obtained by decoupling the mechanical and thermal parts of thermomachical equations. The multiscale characteristics of the thermal expansion coefficients are not considered, which will lead to some limitations of these methods. 
Treating the coupled system as a unified operator and constructing local spectral problems within the framework of the GMsFEM is a natural but impractical idea. This is because the coupling operator exhibits non-positive definite characteristics, especially in cases of strong coupling, i.e., when the coupling coefficients are relatively large, which renders traditional methods ineffective. To overcome this challenge, we constructed a regularized local spectral problem by introducing two relaxation coefficients in the CGMsFEM. Through the regularization technique, the local spectral problems will be solved by a unified operator without decoupling the thermomechanical equations. The multiscale basis functions are constructed by considering the coupling effects of the displacement and temperature fields. Thus, the coupling information can be more accurately captured. In addition, the theoretical results of the GMsFEM method \cite{EGH:JCP:2013a,ZYZZ:CS:2013,ZZLZ:IJMMD:2020} are also generalized. The convergence analysis of the CGMsFEM is obtained, and the error estimations are derived, where the rate of decay of eigenvalues for local spectral problems is obtained, which is similar to \cite{CEL:JCP:2014}. Furthermore, several numerical examples associated with periodic and random microstructure material coefficients are presented to confirm the theoretical analysis of CGMsFEM. When the same number of basis functions are chosen, the numerical results show that the CGMsFEM not only has better accuracy than the GMsFEM in the case of weak coupling but also has obvious advantages in accuracy for the strong coupling situation. The novelty of the proposed  approach is  highlighted as follows: 
\begin{itemize}
  \item[(1)] The multiscale basis functions are constructed by the coupled spectral problems incorporating local multiscale physical information, which can more accurately approximate the solutions of the original thermomechanical problems with fewer degrees of freedom.
  \item[(2)] Two relaxation coefficients are creatively designed for the local coupling spectral problems in each coarse-grid block, the proper regularity of which is ensured by adjusting the value of the relaxation coefficient. Moreover, the method can be reduced to the standard GMsFEM when the multiscale coupling correlation effect disappears, which also provides a generalized framework to design the coupling multiscale basis functions for weak or strong coupled multiphysics problems.
  \item[(3)] Through convergence analysis, it is deduced that the error of CGMsFEM is closely related to the eigenvalue decay in each local coarse block, which is consistent with GMsFEM. At the same time, it is also proven that the upper error bound is independent of the two relaxation coefficients.
\end{itemize}

The paper is organized as follows. The formulation of CGMsFEM is given in Section \ref{sec:basis}, including the problem statement, the construction of coupling multiscale basis functions, and the corresponding finite element method. In Section \ref{sec:analysis}, the convergence analysis of CGMsFEM is carried out. A few numerical examples are implemented to illustrate the efficiency and accuracy of the proposed CGMsFEM in Section \ref{sec:numexp}, and the conclusions follow in Section \ref{sec:conclusions}.

\section{Formulation of CGMsFEM}
\label{sec:basis}
\subsection{Problem statement}
\label{subsec:problem}
A \textcolor{black}{thermomechanical} problem is considered to describe the quasi-static deformation of heterogeneous media in domain $\Omega$, where $\Omega\subset\mathbb{R}^d, ~ d=2,3$ is a convex bounded polygonal or polyhedaral domain with Lipschitz continuous boundary $\partial \Omega$, and $T>0$ is a given time. It is assumed that the solid phase is a linear elastic solid and that the deformation is coupled with the temperature field \cite{DCAK:JMS:2018, PP:SMO:2010}. Then the equilibrium and heat transfer \textcolor{black}{equations} are written in standard form as follows
\begin{equation}
  \begin{cases}\vspace{2mm}
    -\nabla\cdot\big(\bsigma\left(\uu\right)-\beta\theta \II \big)=\ff,                & \text{in}~\Omega\times(0,T], \\
    \dot{\theta}-\nabla\cdot\big(\kappa\nabla\theta\big)+\beta\nabla\cdot \dot{\uu}=g, & \text{in}~\Omega\times(0,T], \\
  \end{cases}
  \label{eq:govequs}
\end{equation}
where $\bsigma$ is the stress tensor and $\II$ is the d-dimensional identity matrix.
$\uu(\xx,t):\Omega\times(0,T]\rightarrow\mathbb{R}^d$ and $\theta(\xx,t):\Omega\times(0,T]\rightarrow\mathbb{R}$ denote the displacement and temperature fields, respectively. Here, the superscript dot represents partial differentiation with respect to time $t$, and $\xx$ denotes the space coordinates. $\kappa(\xx)$ and $\beta(\xx)$ are the thermal conductivity and expansion coefficients with multiscale characteristics, $\ff$ is the body force, and $g$ is the heat source. Then the initial and boundary condition are defined by
\begin{equation}
  \left\{
  \begin{array}{llll}\vspace{2mm}
    \uu=\uu_D,            & \text{on}~\Gamma_D^u\times(0,T],      & \left(\bsigma(u)-\beta\theta \II\right)\cdot \nn=\bsigma_N, & \text{on}~\Gamma_N^u\times(0,T],      \\\vspace{2mm}
    \theta=\theta_D,      & \text{on}~\Gamma_D^\theta\times(0,T], & \left(\kappa\nabla\theta(x,t)\right)\cdot \nn=q_N,          & \text{on}~\Gamma_N^\theta\times(0,T], \\
    \theta(x,0)=\theta_0, & \text{in}~\Omega,
  \end{array}\right.
  \label{eq:boundary}
\end{equation}
where the boundary of $\Omega$ is divided into surfaces $\Gamma_D^u$ and $\Gamma_N^u$, $\Gamma_D^\theta$ and $\Gamma_N^\theta$ with $\Gamma_D^u\cap\Gamma_N^u=\emptyset$ and $\Gamma_D^\theta\cap\Gamma_N^\theta=\emptyset$. $\uu_D$ and $\theta_D$ are the prescribed displacement and temperature, and $\bsigma_N$ and $q_N$ denote surface load value and heat flux, respectively. \textcolor{black}{$\nn$ is the exterior normal to the surface $\Gamma_N^u$ and $\Gamma_N^\theta$.} The constitutive relation of stress and strain for a general thermoelastic material is given by
\begin{equation} \nonumber
  \bsigma(\uu)=2\mu\bepsilon(\uu)+\lambda\nabla\cdot \uu\II,
\end{equation}
where $\bepsilon(\uu)$ is the strain tensor defined by
\begin{equation} \nonumber
  \bepsilon(\uu)=\frac{1}{2}(\nabla \uu+\nabla \uu^T).
\end{equation}
$\mu$ and $\lambda$ are the $\text{Lam}\acute{e} $ constants defined by
\begin{equation}. \nonumber
  \mu=\frac{E}{2(1+\nu)},\lambda=\frac{E\nu}{(1+\nu)(1-2\nu)},
\end{equation}
where $E$ and $\nu$ are Young's modulus and Poisson's ratio.

Here, some definitions and  assumptions will first be stated.  Let $H^1(\Omega)$ denote the classical Sobolev space equipped with the norm $||v||_{H^1(\Omega)}^2=||v||^2_{L^2(\Omega)}+||\nabla v||^2_{L^2(\Omega)}$, and $H^{-1}(\Omega)$ represents the dual space of $H^1$. Then define $L^p([0,T];V)$ for the Bochner space with the norm
\begin{equation*}
  \begin{array}{lll}
    ||v||_{L^p([0,T];V)}  =\ds\left(\int_0^T ||v||_V^p dt\right)^{\frac{1}{p}},\quad 1\leq p<\infty, \\
    ||v||_{L^{\infty}([0,T];V)}  =\mathop{\essinf}\limits_{0\leq t \leq T}||v||_V,
  \end{array}
\end{equation*}
where $V$ ia a Banach space with the norm $||\cdot||_V$, such as $L^2(\Omega)$ and $H^{1}(\Omega)$. For the sake of simplicity, the time interval $[0,T]$ and the domain $\Omega$ will be omitted, for instance, $L^p(V)$ for $L^p([0,T];V)$. The following two \textcolor{black}{spaces} $V_u(\Omega)$ and $V_\theta(\Omega)$ are also defined as
$$V_u(\Omega):=\Big\{v|v\in \big(H^1(\Omega)\big)^d:v=0\ \mathrm{on}\ \Gamma_D^u\Big\},\quad V_\theta(\Omega):=\Big\{v|v\in H^1(\Omega):v=0\ \mathrm{on}\ \Gamma_D^\theta\Big\}.$$
The following assumption is also presented.
\begin{myassum}
  The material parameters $\lambda$, $\mu$, $\kappa$, $\beta\in L^{\infty}(\Omega)$ satisfied
  \begin{equation*}
    \begin{array}{lll}
      \ds 0 < \lambda_{\text{min}}:=\mathop{\essinf}\limits_{\xx\in \Omega} \lambda(\xx) \leq\mathop{\esssup}\limits_{\xx\in \Omega}\lambda(\xx):=\lambda_{\text{max}} < \infty, \\
      \ds  0 < \mu_{\text{min}}:=\mathop{\essinf}\limits_{\xx\in \Omega} \mu(\xx) \leq\mathop{\esssup}\limits_{\xx\in \Omega} \mu(\xx) =\mu_{\text{max}} < \infty,                \\
      \ds  0 < \kappa_{\text{min}}:=\mathop{\essinf}\limits_{\xx\in \Omega} \kappa(\xx) \leq\mathop{\esssup}\limits_{\xx\in \Omega} \kappa(\xx):=\kappa_{\text{max}} < \infty,    \\
      \ds  0 < \beta_{\text{min}}:=\mathop{\essinf}\limits_{\xx\in \Omega}\beta(\xx) \leq\mathop{\esssup}\limits_{\xx\in \Omega} \beta(\xx):=\beta_{\text{max}} < \infty.
    \end{array}
  \end{equation*}
\end{myassum}

Furthermore, the corresponding variational formulation of problem (\ref{eq:govequs}) is to find the weak solutions $\uu\in V_u(\Omega)$ and $\theta\in V_{\theta}(\Omega)$ such that
\begin{equation}
  \begin{cases}\vspace{2mm}
    a(\uu,\vv_u)-b(\vv_u,\theta) = \left\langle \ff,\vv_u \right\rangle ,                                             & \forall \vv_u \in V_u(\Omega),         \\
    c(\dot{\theta}, v_\theta) + d(\theta,v_\theta) + b(\dot{\uu},v_\theta) =  \left\langle g, v_\theta \right\rangle, & \forall v_\theta \in V_\theta(\Omega),
  \end{cases}\label{eq:thermovariation}
\end{equation}
where
\begin{equation*}
  \begin{array}{lll}\vspace{2mm}
    \ds a(\uu,\vv_u) =\int_{\Omega} \bsigma(\uu):\bepsilon(\vv_u), \quad b(\vv_u,\theta) = \int_{\Omega} \beta \theta \nabla \cdot \vv_u, \\
    \ds c(\theta,v_\theta) = \int_{\Omega} \theta v_\theta , \quad d(\theta,v_\theta) = \int_{\Omega} \kappa \nabla \theta \cdot \nabla v_\theta.
  \end{array}
\end{equation*}
 : represents the Frobenius inner product, and $\left\langle \cdot,\cdot \right\rangle$ denotes the inner product in $L^2(\Omega)$.

\subsection{Construction of coupling multiscale basis functions}
\label{subsec:basis}
The procedure of the construction of the CGMsFEM basis functions can be divided into two steps. Firstly, \textcolor{black}{coupling multiscale eigenfunctions} can be constructed by solving regularized coupling local spetral problems in coarse block, and the local computations associated with construction of local approximation spaces are independent and can be performed in parallel. \textcolor{black}{Then the partition of unity functions are obtained to construct multiscale basis functions.} Finally, the resulting global stiffness matrix can be several orders of magnitude smaller than the stiffness matrix obtained by applying FEM directly.

\begin{figure}[H]
  \centering
  \includegraphics[width=0.85\linewidth, height=0.28\textheight]{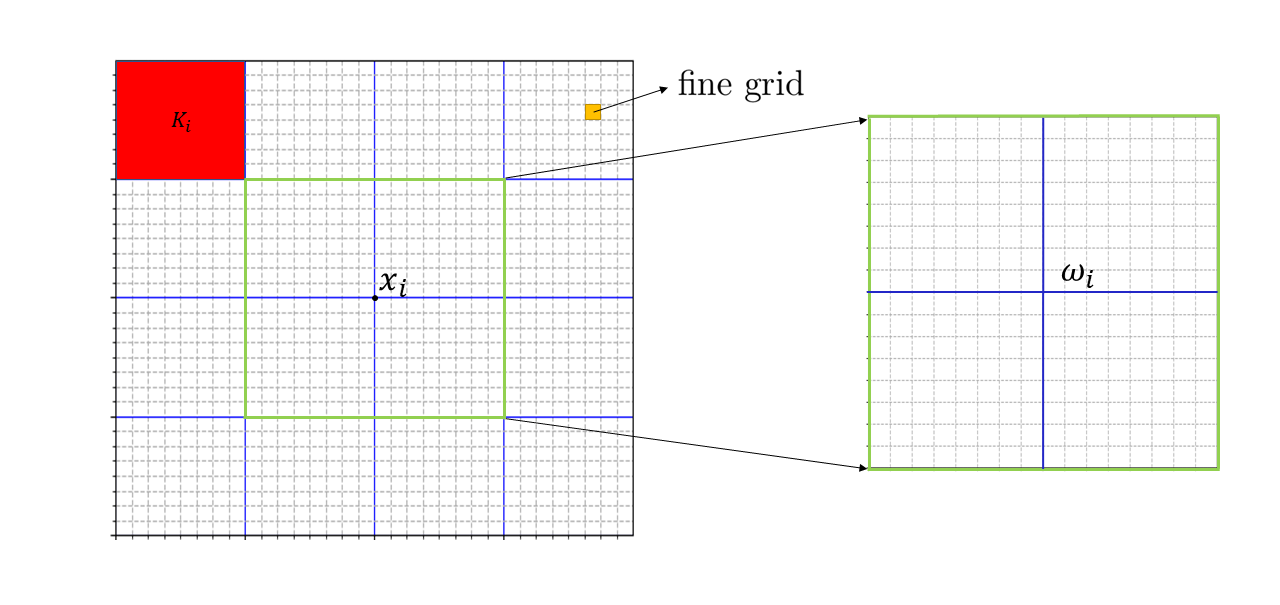}
  \caption{The fine grid $\mathcal{T}^h$,  the coarse grid $\mathcal{T}^H$, the coarse element $K_{i}$ and  neighborhood $\omega_{i}$ of the node $x_{i}$.}
  \label{sec2-fig1}
\end{figure}

Suppose the domain $\Omega$ is composed of a family of meshes $\mathcal{T}^{H}$, where $H=\underset{K_i\in\mathcal{T}^{H}}{\max}H_{K_i}$ is the coarse mesh size and $H_{K_i}$ is the diameter of coarse grid $K_i$. $\mathcal{T}^{H}$ is the conforming partition of $\Omega$ and is shape-regular. $\mathcal{T}^{h}$ is a conforming refinement of $\mathcal{T}^{H}$, and $h$ is the diameter of the fine grid. $N$ denotes the number of elements in $\mathcal{T}^{H}$, and $N_{v}$ denotes the number of vertices of the coarse grid. Let $\{x_{i}\}_{i=1}^{N_{v}}$ be the set of vertices in $\mathcal{T}^{H}$ and $\omega_{i}=\bigcup\{K_{j}\in \mathcal{T}^{H}|x_{i}\in\overline{K_{j}}\}$ the neighborhood of the node $x_{i}$.
Figure \ref{sec2-fig1} depicts the fine grid, the coarse element $K_i$, and the coarse neighborhood $\omega_{i}$ of the node $x_{i}$.

In order to construct coupling multiscale basis functions, the following {regularized} coupling spectral problem in each coarse neighborhood $\omega_{i}$ is given to find eigen-pairs
$\{\Lambda^{\omega_i},\psi^{\omega_i}\}$ such that
\begin{equation} 
  \begin{cases} \vspace{2mm}
    \ds    -\nabla\cdot\Big(\bsigma(\psi_{u}^{\omega_i})-\gamma_1\beta\psi_{\theta}^{\omega_i}\II\Big)=\frac{1}{H^2}\Lambda^{\omega_i}(\lambda+2\mu)\psi_{u}^{\omega_i}, & {\mathrm{in}\ \omega_i,}         \\
    \ds    -\nabla\cdot \kappa\nabla\psi_{\theta}^{\omega_i}+\gamma_2\beta\nabla\cdot\psi_{u}^{\omega_i}=\frac{1}{H^2}\Lambda^{\omega_i}\kappa\psi_{\theta}^{\omega_i},  & {\mathrm{in}\ \omega_i,}         \\
    \bsigma(\psi_{u}^{\omega_i})\cdot\nn=0,                                                                                                                   & {\mathrm{on}\ \partial\omega_i,} \\
    \kappa\nabla\psi_{\theta}^{\omega_i}\cdot\nn=0,                                                                                                           & {\mathrm{on}\ \partial\omega_i.}
  \end{cases}\label{eq:eigproblem}
\end{equation}
$\psi^{\omega_i}=(\psi_{u}^{\omega_i},\psi_{\theta}^{\omega_i})$, $\psi_{u}^{\omega_i}$ and $\psi_{\theta}^{\omega_i}$ are corresponding eigenfunctions for displacement and temperature fields in problem (\ref{eq:eigproblem}). 
$\gamma_1$ and $\gamma_2$ are the relaxation coefficients, the proper regularity of which is ensured by adjusting the value of relaxation coefficient. To be precise, it is always possible to choose appropriate $\gamma_1$ and $\gamma_2$ such that this spectral problem is positive definite. The introduction of $\frac{1}{H^2}$ in the right-hand side of this spectral problem is to ensure that $\Lambda^{\omega_i}$ is independent of $H$.
Here, it should be emphasized that the important coupling physical characteristics are naturally incorporated into multiscale basis functions, which can obviously reduce the orders of global stiffness matrix.

Let $V_{uh}(\Omega)$ and $V_{\theta h}(\Omega)$ be classical affine finite element spaces on $\mathcal{T}^{h}$ of $V_u(\Omega)$ and $V_{\theta}(\Omega)$ respectively, and the definitions of $V_{uh}(\omega_i)$ and $V_{\theta h}(\omega_i)$ are also similar.
The corresponding variational formulation of problem (\ref{eq:eigproblem}) is to find the weak solutions $\psi_{u}^{\omega_i}\in V_{uh}(\omega_i)$ and $\psi_{\theta}^{\omega_i}\in V_{\theta h}(\omega_i)$ such that
\begin{equation}
  \begin{cases}\vspace{2mm}
    \ds  a^{\omega_i}(\psi_{u}^{\omega_i},\vv_{uh})-\gamma_1 b^{\omega_i}(\vv_{uh},\psi_{\theta}^{\omega_i}) = \frac{1}{H^2}\Lambda^{\omega_i} \left\langle \psi_{u}^{\omega_i},\vv_{uh} \right\rangle_a^{\omega_i} ,                    & \forall \vv_{uh} \in V_{uh}(\omega_i),           \\
    \ds  d^{\omega_i}(\psi_{\theta}^{\omega_i},v_{\theta h}) + \gamma_2 b^{\omega_i}(\psi_{u}^{\omega_i},v_{\theta h}) = \frac{1}{H^2}\Lambda^{\omega_i} \left\langle \psi_{\theta}^{\omega_i}, v_{\theta h} \right\rangle_d^{\omega_i}, & \forall v_{\theta h} \in V_{\theta h}(\omega_i),
  \end{cases}\label{eq:varspectral}
\end{equation}
where
\begin{equation*}
  \begin{array}{lll}\vspace{2mm}
    \ds a^{\omega_i}(\psi_{u}^{\omega_i},\vv_{uh}) =\int_{\omega_i} \bsigma(\psi_{u}^{\omega_i}):\bepsilon(\vv_{uh}),\quad b^{\omega_i}(\vv_{uh},\psi_{\theta}^{\omega_i}) = \int_{\omega_i} \beta \psi_{\theta}^{\omega_i} \nabla \cdot \vv_{uh}, \\
    \ds  d^{\omega_i}(\psi_{\theta}^{\omega_i},v_{\theta h}) = \int_{\omega_i} \kappa \nabla \psi_{\theta}^{\omega_i} \cdot \nabla v_{\theta h},                                                                                                   \\
    \ds  \left\langle \psi_{u}^{\omega_i},\vv_{uh} \right\rangle_a^{\omega_i} = \int_{\omega_i} (\lambda+2\mu)\psi_{u}^{\omega_i} \cdot \vv_{uh}, \quad \left\langle \psi_{\theta}^{\omega_i}, v_{\theta h} \right\rangle_d^{\omega_i} = \int_{\omega_i} \kappa \psi_{\theta}^{\omega_i} v_{\theta h}.
  \end{array}
\end{equation*}
In the following text, substituting $\omega_i$ with $K$ signifies that it maintains the same definition within the region $K$. Thus, following the standard finite element discretization, we obtain the following algebraic eigenvalue problems,
\begin{equation}
  \bA^{\omega_i}\psi^{\omega_i}=\frac{1}{H^2}\Lambda^{\omega_i} \MM^{\omega_i}\psi^{\omega_i},
  \label{eq:eigMatrix2}
\end{equation}
where
\begin{equation*}
  \bA^{\omega_i}=
  \left[
    \begin{array}{cc}
      A_1^{\omega_i}         & -\gamma_1A_2^{\omega_i} \\
      \gamma_2A_3^{\omega_i} & A_4^{\omega_i}          \\
    \end{array}
    \right],
  \MM^{\omega_i}=
  \left[
    \begin{array}{cc}
      M_1^{\omega_i} & 0              \\
      0              & M_2^{\omega_i} \\
    \end{array}
    \right],
  \label{eq:eigKM2}
\end{equation*}
and
$$(A_1)^{\omega_i}_{j_1j_1}=a^{\omega_i}(\phi_{u,j_1},\phi_{u,j_1}),\
  (A_2)^{\omega_i}_{j_1j_2}=b^{\omega_i}(\phi_{u,j_1},\phi_{\theta,j_2}),\ A_3= A_2^T,$$
$$(A_4)^{\omega_i}_{j_2j_2}=d^{\omega_i}(\phi_{\theta,j_2},\phi_{\theta,j_2}),\ (M_1)^{\omega_i}_{j_1j_1}=\left\langle \phi_{u,j_1},\phi_{u,j_1} \right\rangle_a^{\omega_i},\ (M_2)^{\omega_i}_{j_2j_2}=\left\langle\phi_{\theta,j_2},\phi_{\theta,j_2} \right\rangle_d^{\omega_i}.$$
$\phi_{u,j}\in V_{uh}(\omega_i)$ and $\phi_{\theta,j}\in V_{\theta h}(\omega_i)$, $ \textcolor{black}{1\leq j_1\leq d N(\omega_i)},1\leq j_2\leq N(\omega_i)$, and $N(\omega_i)$ is the number of total fine grid nodes in $\omega_i$. 

Then the eigenfunctons $\{\psi^{\omega_i}_{l}\}_{l=1}^{L_i}$ related to the smallest $L_i$ eigenvalues are chosen as members of CGMsFEM space, which is defined by
\begin{equation} \nonumber
  V_{cgm}=\text{span}\Big\{\Phi^{\omega_i}_{l}|\Phi^{\omega_i}_{l}=\chi_i^T \II_{d+1}\psi^{\omega_i}_{l}, 1\leq i\leq N_v\ \mathrm{and}\ 1\leq l\leq L_i\Big\},
\end{equation} 
where $\{\chi_i\}_{i=1}^{N_v}$ is a set of partition of unity functions associated with the open cover $\{\omega_i\}_{i=1}^{N_v}$ of domain $\Omega$. The coupling multiscale basis function set
$\{\Phi^{\omega_i}_{l}\}$ can be placed in the following matrix
\begin{equation} \nonumber
  \RR=\Big[\Phi^{\omega_1}_1,\Phi^{\omega_1}_2...,\Phi^{\omega_1}_{L_1},\Phi^{\omega_2}_1,\Phi^{\omega_2}_2...,\Phi^{\omega_2}_{L_2},...,\Phi^{\omega_{N_v}}_1,\Phi^{\omega_{N_v}}_2...,\Phi^{\omega_{N_v}}_{L_{N_v}}\Big],\label{eq:base_func}
\end{equation}
Here, it should be noted that the matrix $\RR$ only needs to be constructed once, and it can be repeatedly used for computation.
\textcolor{black}{
\begin{remark}
  For multi-physics problems, it is more appropriate to select the partition of unit functions $\{\chi_i\}_{i=1}^{N_v}$ separately for each component. This paper chooses the standard Lagrangian basis functions as the partition of unit functions. The MsFEM basis functions of the elasticity operator and diffusion operator will result in slightly better performance.
\end{remark}
}

\subsection{Algorithm procedure of CGMsFEM}
\label{subsec:alg}
In the paper, the backward Euler scheme is used for temporal discretization. Let $\ww^n=(\uu_h^n,\theta_h^n)$ be the solutions at the n-th time level $t_n=\sum_{m=1}^{n}\tau_m,~n\in\{0,1,...,N_T\}$, where $\tau_m$ is the time step. Then we have the following weak formulation for the Eqs. (\ref{eq:thermovariation}),
\begin{equation}
  \begin{cases}\vspace{2mm}
    \ds \textcolor{black}{ a(\uu_h^n,\vv_{uh}) - b(\vv_{uh},\theta_h^n) } = \left\langle \ff^n,\vv_{uh} \right\rangle ,   & \forall \vv_{uh} \in V_{uh}(\Omega),       \\
    \ds  c\bigg(\frac{\theta^n_h-\theta^{n-1}_h}{\tau_n}, v_{\theta h}\bigg) + d(\theta_h^n,v_{\theta h}) + b\bigg(\frac{\uu^n_h-\uu^{n-1}_h}{\tau_n},v_{\theta h}\bigg) =  \left\langle g^n, v_\theta \right\rangle, & \forall v_{\theta h} \in V_{\theta h}(\Omega),
  \end{cases}\label{eq:discrevariation}
\end{equation}
Thus, following the standard finite element discretization, Eqs. (\ref{eq:discrevariation}) can be rewritten as follows
\begin{equation} 
  \bA^n \ww^n=\BB \ww^{n-1}+\FF^n,\label{eq:FEmatrix}
\end{equation}
where
\begin{equation*}
  \bA^n =
  \left[
    \begin{array}{cc}
      A_1 & -A_2 \\A_3&M'+\tau_n A_4\\
    \end{array}
    \right],
  \BB=
  \left[
    \begin{array}{cc}
      0 & 0 \\A_3&M'\\
    \end{array}
    \right],
  \FF^n=
  \left[
    \begin{array}{cc}
      F^n \\\tau_n G^n
    \end{array}
    \right],
\end{equation*}
and
$$(A_1)_{j_1j_1}=a(\phi_{u,j_1},\phi_{u,j_1}),\
  (A_2)_{j_1j_2}=b(\phi_{u,j_1},\phi_{\theta,j_2}),\ A_3=A_2^T,$$
$$(A_4)_{j_2j_2}=d(\phi_{\theta,j_2},\phi_{\theta,j_2}),\ (M')_{j_2j_2}=c(\phi_{\theta,j_2},\phi_{\theta,j_2}),$$
$$(F^n)_{j_1}=(\ff^n,\phi_{u,j_1}),\ (G^n)_{j_2}=(g^n,\phi_{\theta,j_2}).$$
Define
$$\bA^n_c=\RR^T\bA^n\RR,~\BB^c=\RR^T\BB\RR,~\ww^n_c=(\RR^T\RR)^{-1}\RR^T\ww^n,~\FF^n_c=(\RR^T\RR)^{-1}\RR^T\FF^n,\ $$
and combining with Eq. (\ref{eq:FEmatrix}), the following algebraic system can be given,
\begin{equation}
  \bA^n_c\ww^n_c=\BB^c \ww^{n-1}_c+\FF^n_c. \label{eq:CGMsFE}
\end{equation}
Then the solution $\ww^n_c$ can be calculated by iteration, and the solutions in fine grid can be obtained by use of the coupling multiscale basis functions,
\begin{equation} \nonumber 
  \ww^n=\RR\ww^n_c.
\end{equation}
Here, the key difference between CGMsFEM and GMsFEM is the construction of more accurate coupling multiscale basis functions. The matrix $\RR$ can be computed offline, and \textcolor{black}{it} can be reused in all time steps. Thus, a much smaller system is solved with CGMsFEM, which can significantly improve the computational efficiency.

\section{Convergence analysis of CGMsFEM}
\label{sec:analysis}
\subsection{Interpolation error}
Before proving the convergence and error estimates for CGMsFEM, some definitions and Lemmas are first given.
Let $\ww = (\uu,\theta), \vv = (\vv_{u},v_{\theta}) \in V_u(\Omega) \times V_{\theta}(\Omega)$. 
The $ b(\cdot,\cdot) $ is continuous in $ V_{u}(\Omega) \times L^2(\Omega) $, and it can be assumed that there exists a constant $C_0$ such that for all $\uu \in V_{u}(\Omega), \theta \in L_2(\Omega)$,
we have
\begin{equation} \nonumber
  b^{K} \left( \uu,\theta \right) \leq C_0 \left\| \uu \right\|_{a,K} \left\| \theta \right\|_{L_d,K}, \quad \forall~K \in \mathcal{T}_H,
\end{equation}
where
\begin{equation*}
  \ds \left\| \theta \right\|_{L_d,K}^2 = \int_{K} \kappa \theta^2.
\end{equation*}
Also define
\begin{equation*}
 \left\| \uu \right\|_{L_a,K}^2 = \int_{K} \left( \lambda+2\mu \right) \uu \cdot \uu.
\end{equation*}
And for all $\uu \in V_{u}(\Omega), \theta \in V_{\theta}(\Omega)$, the energy norms can be defined as follows
\begin{equation} \nonumber
  \|\uu\|_{a,K}^2 = a^{K}(\uu,\uu),\quad \|\theta\|_{d,K}^2 = d^{K}(\theta,\theta).
\end{equation}
In the following text, substituting $K$ with $\omega_i$ signifies that it maintains the same definition within the region $\omega_i$. Meanwhile, the local bilinear functionals can be defined as follows
\begin{equation}
  \begin{aligned}
    \mathcal{A}^{\omega_i}(\ww,\vv) & = a^{\omega_i}(\uu,\vv_{u}) - \gamma_1 b^{\omega_i}(\vv_{u},\theta) + \gamma_2 b^{\omega_i}(\uu, v_{\theta}) +d^{\omega_i}(\theta,v_{\theta}), \\
    \mathcal{M}^{\omega_i}(\ww,\vv) & = \int_{\omega_i} (\lambda+2\mu)\uu \cdot \vv_{u}+\kappa \theta v_{\theta}.
  \end{aligned}
\end{equation}
Then the local eigenvalue problem (\ref{eq:eigproblem}) can be written as
\begin{equation}
  \mathcal{A}^{\omega_i}(\bpsi^{\omega_i},\vv) = \frac{1}{H^2} \Lambda^{\omega_i} \mathcal{M}^{\omega_i}(\bpsi^{\omega_i},\vv), \quad  \forall  \vv \in V_u(\omega_i) \times V_{\theta}(\omega_i).
\end{equation}
The obtained eigenvalues can be arranged in ascending order
\begin{equation} \nonumber
  \Lambda^{\omega_i}_1 \leq \Lambda^{\omega_i}_2 \leq \cdots \leq \Lambda^{\omega_i}_L \leq \cdots .
\end{equation}
Assume that except for $d+1$ zero eigenvalues, all eigenvalues are real numbers and greater than 0, and define $\widetilde{V}\left(\omega_i\right):= V_u(\omega_i) \times V_{\theta}(\omega_i)$.

We have
\begin{equation} \nonumber
  \vv = \sum_{l=1}^{\infty} \mathcal{M}^{\omega_i}(\vv,\bpsi_l^{\omega_i}) \bpsi_l^{\omega_i},\quad \forall \vv \in \widetilde{V}\left(\omega_i\right),
\end{equation}
where $\{ \bpsi_l^{\omega_i} \}_{l=1}^{\infty}$ are a set of complete orthogonal basis functions in $\widetilde{V}\left(\omega_i\right)$ with the inner product $M^{\omega_i}(\cdot,\cdot)$.
It is easily obtained that
\begin{equation} 
  \mathcal{A}^{\omega_i}(\vv,\vv) = \frac{1}{H^2} \sum_{l=1}^{\infty} \mathcal{M}^{\omega_i}(\vv,\bpsi_l^{\omega_i})^2 \Lambda^{\omega_i}_l,\quad \mathcal{M}^{\omega_i}(\vv,\vv) = \sum_{l=1}^{\infty} \mathcal{M}^{\omega_i}(\vv,\bpsi_l^{\omega_i})^2.
\end{equation}
Inspired by Reference \cite{ADG:CAfGMs:2019}, we define the semi-norm by
\begin{equation}
  \interleave  \vv  \interleave_{s,\omega_i}^2 = \sum_{l=1}^{\infty} \left( \frac{\Lambda^{\omega_i}_l}{H^2} \right)^s \mathcal{M}^{\omega_i}(\vv,\bpsi_l^{\omega_i})^2  \text{ for }0 \leq s \leq 2,
\end{equation}
and it follows that
\begin{equation} \nonumber
  \interleave  \vv  \interleave_{0,\omega_i}^2 = \mathcal{M}^{\omega_i}(\vv,\vv),  \quad \interleave  \vv  \interleave_{1,\omega_i}^2 = \mathcal{A}^{\omega_i}(\vv,\vv).
\end{equation}
It's important to emphasize that the eigenvalues of the local spectral problem (\ref{eq:eigproblem}) are $\left\{ \frac{\Lambda^{\omega_i}_l}{H^2} \right\}$. Therefore, the seminorm is defined in this manner.
Following the above definitions, we give
\begin{mylem}
  Assume that $\ww \in V_u(\omega_i) \times V_{\theta}(\omega_i)$ and \textcolor{black}{$ \interleave  \ww  \interleave_{s,\omega_i} < \infty $}.
  For $0\leq t \leq s \leq 2 , L_i \geq d+1$, there holds
  \begin{equation}
    \interleave  \ww-\mathcal{I}_{L_i}^{\omega_i} \ww  \interleave_{t,\omega_i}^2 \leq \left( \frac{\Lambda_{L_i+1}^{\omega_i}}{H^2} \right)^{t-s} \interleave  \ww  \interleave_{s,\omega_i}^2,
  \end{equation}
  where $\mathcal{I}_{L_i}^{\omega_i}$ is the local interpolation operator defined by
  \begin{equation} \nonumber
    \mathcal{I}_{L_i}^{\omega_i} \ww = \sum_{l=1}^{L_i} \mathcal{M}^{\omega_i}(\ww,\bpsi_l^{\omega_i}) \bpsi_l^{\omega_i}. \label{eq:interop}
  \end{equation}
\end{mylem}

\begin{mypf}
  By definition (\ref{eq:interop}), it can be obtained that
  \begin{equation*}
    \begin{aligned}
      \interleave  \ww-\mathcal{I}_{L_i}^{\omega_i} \ww  \interleave_{t,\omega_i}^2 = & \sum_{l=L_i+1}^{\infty} \left( \frac{\Lambda^{\omega_i}_l}{H^2} \right)^t \mathcal{M}^{\omega_i}(\ww,\bpsi_l^{\omega_i})^2                                                             \\
      =                                                                               & \sum_{l=L_i+1}^{\infty} \left( \frac{\Lambda^{\omega_i}_l}{H^2} \right)^{t-s} \left( \frac{\Lambda^{\omega_i}_l}{H^2} \right)^s \mathcal{M}^{\omega_i}(\ww,\bpsi_l^{\omega_i})^2       \\
      \leq                                                                            & \left( \frac{\Lambda_{L_i+1}^{\omega_i}}{H^2} \right)^{t-s} \sum_{l=L_i+1}^{\infty} \left( \frac{\Lambda^{\omega_i}_l}{H^2} \right)^s \mathcal{M}^{\omega_i}(\ww,\bpsi_l^{\omega_i})^2 \\
      \leq                                                                            & \left( \frac{\Lambda_{L_i+1}^{\omega_i}}{H^2} \right)^{t-s} \interleave  \ww -\mathcal{I}_{L_i}^{\omega_i} \ww \interleave_{s,\omega_i}^2                                    \\
      \leq                                                                            & \left( \frac{\Lambda_{L_i+1}^{\omega_i}}{H^2} \right)^{t-s} \interleave  \ww  \interleave_{s,\omega_i}^2.
    \end{aligned}
  \end{equation*}
\end{mypf}
Then the following corollary can be given by this lemma with $t =0, 1$.
\begin{myinfer}
  For all $\ww \in V_u(\omega_i) \times V_{\theta}(\omega_i), 1 \leq s \leq 2 , L_i \geq d+1$,
  \begin{equation} 
    \begin{aligned}
      \mathcal{M}^{\omega_i}(\ww-\mathcal{I}_{L_i}^{\omega_i} \ww,\ww-\mathcal{I}_{L_i}^{\omega_i} \ww)
       & \leq \frac{H^{2s}}{\left( \Lambda_{L_i+1}^{\omega_i} \right)^{s}} \interleave  \ww  \interleave_{s,\omega_i}^2,     \\
      \mathcal{A}^{\omega_i}(\ww-\mathcal{I}_{L_i}^{\omega_i} \ww,\ww-\mathcal{I}_{L_i}^{\omega_i} \ww)
       & \leq \frac{H^{2s-2}}{\left( \Lambda_{L_i+1}^{\omega_i} \right)^{s-1}} \interleave  \ww  \interleave_{s,\omega_i}^2.
    \end{aligned}
  \end{equation}
\end{myinfer}
For interpolation operator $\mathcal{I}$: $ V_u(\Omega) \times V_{\theta}(\Omega) \rightarrow V_u(\Omega) \times V_{\theta}(\Omega)$ and all  $\ww = (\uu, \theta) \in V_u(\Omega) \times V_{\theta}(\Omega) $, we define $\left(\mathcal{I} \ww \right)_u$ and $\left(\mathcal{I} \ww \right)_{\theta}$, s.t. $\mathcal{I} \ww = ( \left(\mathcal{I} \ww \right)_u, \left(\mathcal{I} \ww \right)_{\theta} )$.

\begin{mylem}
  For all $\ww = (\uu, \theta) \in V_u(\Omega) \times V_{\theta}(\Omega) $, there holds
  \begin{equation}
    \begin{aligned}
      \int_{K} (\lambda+2\mu) \left( \uu-\left(\mathcal{I}_{\text{ms}} \ww \right)_u \right) \cdot \left( \uu-\left(\mathcal{I}_{\text{ms}} \ww \right)_u \right)+\kappa \left( \theta-\left(\mathcal{I}_{\text{ms}} \ww \right)_{\theta} \right)^2 \\
      \leq N_K \frac{H^4}{ \Lambda_{K,L+1} } \sum_{y_i \in K}  \interleave  \ww  \interleave_{2,\omega_i}^2,
    \end{aligned}\label{eq:wuest}
  \end{equation}
  where $N_K = \text{card}~\{ y_i: y_i \in K\} $  and $\Lambda_{K, L+1}=\min_{y_i \in K} \Lambda_{L_i+1}^{\omega_i}$. $\mathcal{I}_{\text{ms}}: V_u(\Omega) \times V_{\theta}(\Omega) \rightarrow V_{cgm}$ is the global interpolation operator, which is defined by
  \begin{equation} \nonumber
    \mathcal{I}_{\text{ms}} \ww = \sum_{i=1}^{N_v} \sum_{l=1}^{L_i} \mathcal{M}^{\omega_i}(\ww,\psi_l^{\omega_i}) \Phi_l^{\omega_i}.
  \end{equation}
\end{mylem}
\begin{mypf}
  Based on the definition of $\mathcal{I}_{\text{ms}}$, we have
  \begin{equation}
    \begin{aligned}
      \ww - \mathcal{I}_{\text{ms}} \ww 
       & =\ww -  \sum_{i=1}^{N_v} \sum_{l=1}^{L_i} \mathcal{M}^{\omega_i}(\ww,\bpsi_l^{\omega_i}) \chi_i^T(x) I_{d+1} \bpsi_l^{\omega_i} \\
       & =\ww -  \sum_{i=1}^{N_v}  \chi_i^T I_{d+1} \mathcal{I}_{L_i}^{\omega_i} \ww                                           \\
       & = \sum_{i=1}^{N_v} \chi_i^T(x) I_{d+1} \left( \ww-\mathcal{I}_{L_i}^{\omega_i} \ww \right).
    \end{aligned}\label{eq:interp1}
  \end{equation}
  Combining with Eq. (\ref{eq:interp1}), the left part of Eq. (\ref{eq:wuest}) gives
  \begin{equation} \nonumber
    \begin{aligned}
       & \int_{K} (\lambda+2\mu) \left(  \sum_{y_i \in K} \chi_{i,u}^T I_{d} \left( \uu-\left(\mathcal{I}_{L_i}^{\omega_i} \ww \right)_u \right) \right) \cdot \left(  \sum_{y_i \in K} \chi_{i,u}^T I_{d} \left( \uu-\left(\mathcal{I}_{L_i}^{\omega_i} \ww \right)_u \right) \right)                         \\
       & \quad + \kappa \left(  \sum_{y_i \in K} \chi_{i,\theta} \left( \theta-\left(\mathcal{I}_{L_i}^{\omega_i} \ww \right)_{\theta} \right) \right)^2                                                                                                                                                       \\
       & \leq N_K  \sum_{y_i \in K} \int_{K} (\lambda+2\mu)  \left( \chi_{i,u}^T I_{d} \left( \uu-\left(\mathcal{I}_{L_i}^{\omega_i} \ww \right)_u \right) \right) \cdot \left(  \chi_{i,u}^T I_{d} \left( \uu-\left(\mathcal{I}_{L_i}^{\omega_i} \ww \right)_u \right) \right)                                \\
       & \quad + \kappa \left( \chi_{i,\theta} \left( \theta-\left(\mathcal{I}_{L_i}^{\omega_i} \ww \right)_{\theta} \right) \right)^2                                                                                                                                                                         \\
       & \leq N_K \sum_{y_i \in K} \int_{\omega_i} (\lambda+2\mu)  \left( \uu-\left(\mathcal{I}_{L_i}^{\omega_i} \ww \right)_u \right)  \cdot  \left( \uu-\left(\mathcal{I}_{L_i}^{\omega_i} \ww \right)_u \right)  + \kappa  \left( \theta-\left(\mathcal{I}_{L_i}^{\omega_i} \ww \right)_{\theta}  \right)^2 \\
       & = N_K \sum_{y_i \in K} 	\mathcal{M}^{\omega_i}(\ww-\mathcal{I}_{L_i}^{\omega_i} \ww,\ww-\mathcal{I}_{L_i}^{\omega_i} \ww)                                                                                                                                                                                        \\
       & \leq N_K \sum_{y_i \in K}  \frac{H^4}{\left( \Lambda_{L_i+1}^{\omega_i} \right)^{2}} \interleave  \ww  \interleave_{2,\omega_i}^2                                                                                                                                                                     \\
       & \leq  N_K \frac{H^4}{ \left( \Lambda_{K,L+1} \right)^2 } \sum_{y_i \in K}  \interleave  \ww  \interleave_{2,\omega_i}^2.
    \end{aligned}
  \end{equation}
\end{mypf}

\begin{mylem}\label{Lemma:localinter}
  For all $\ww  = (\uu, \theta)  \in V_u\left( \Omega \right) \times V_{\theta}\left( \Omega \right) $ and the definition of $ \mathcal{I}_{\text{ms}} $, the local interpolation error estimate can be given as
  \begin{equation}
    \begin{aligned}
      \left\| \uu-\left(\mathcal{I}_{\text{ms}} \ww \right)_u  \right\|_{a,K} +\left\| \theta-\left(\mathcal{I}_{\text{ms}} \ww \right)_{\theta}   \right\|_{d,K} \leq 2N_K \sum_{y_i \in K} \left(  C_2   \mathcal{A}^K\left( \ww-\mathcal{I}_{\text{ms}} \ww,\ww-\mathcal{I}_{\text{ms}}  \ww\right) \right. \\
      \left. + \left( \frac{2C_1^2}{H^2}+C_3 \right)  \left(  \left\| \uu-\left(\mathcal{I}_{\text{ms}} \ww \right)_u \right\|_{L_a,K} +  \left\| \theta-\left(\mathcal{I}_{\text{ms}} \ww \right)_{\theta}   \right\|_{L_d,K}  \right) \right),
    \end{aligned}\label{eq:local:interp}
  \end{equation}
  {where $\ds C_2= \frac{2}{2-| \gamma_1 - \gamma_2| C_0}$, $\ds C_3 = \frac{| \gamma_1 - \gamma_2| C_0}{2-| \gamma_1 - \gamma_2| C_0}$, and $C_1$ is a constant such that}
  \begin{equation} \nonumber
    \mathop{\max}   \left| \nabla \chi_i \right| \leq \frac{C_1}{H},\quad \forall \ i \leq N_v.
  \end{equation}
\end{mylem}
The proof of this lemma can be found in \ref{sec:appendix:lemma}.

  

Lemma \ref{Lemma:localinter} directly leads to the following estimation of global interpolation error.
\begin{mythm}\label{theom:globalinter}
  For all $\ww \in V_u\left( \Omega \right) \times V_{\theta}\left( \Omega \right) $ and definition of $ \mathcal{I}_{\text{ms}} $, the global interpolation error can be given as
  \begin{equation}
    \begin{aligned}
       & \left\| \uu-\left(\mathcal{I}_{\text{ms}} \ww \right)_u  \right\|_{a,\Omega}^2 +\left\| \theta-\left(\mathcal{I}_{\text{ms}} \ww \right)_{\theta}  \right\|_{d,\Omega}^2                                                        \\
       & \quad \quad \leq 2N_{\text{max}} \left( C_2 \frac{H^2}{ \Lambda_{L+1} } +  \left( 2C_1^2 H^2+C_3 H^4  \right) \frac{1}{\left( \Lambda_{L+1} \right)^{2}}  \right) \sum_{\omega_i} \interleave  \ww  \interleave_{2,\omega_i}^2,
    \end{aligned}
  \end{equation}
  where $ \Lambda_{L+1}=\min _{\omega_i} \Lambda_{L_i+1}^{\omega_i}  $ and $ N_{\text{max}}= \max_{K} N_K $.
\end{mythm}
\begin{mypf}
  From Eq. (\ref{eq:local:interp}), the following estimations can be given as
  \begin{equation}
    \begin{aligned}
      \sum_{K \in \mathcal{T}^H} \sum_{y_i \in K} \mathcal{A}^K\left( \ww-\mathcal{I}_{L_i}^{\omega_i} \ww,\ww-\mathcal{I}_{L_i}^{\omega_i}  \ww\right) & = \sum_{\omega_i} \mathcal{A}^{\omega_i} \left( \ww-\mathcal{I}_{L_i}^{\omega_i} \ww,\ww-\mathcal{I}_{L_i}^{\omega_i}  \ww\right) \\
                                                                                                                                              & \leq \sum_{\omega_i}  \frac{H^2}{ \Lambda_{L_i+1}^{\omega_i}  }\interleave  \ww  \interleave_{2,\omega_i}^2,
    \end{aligned}\label{eq:akest}
  \end{equation}
  \begin{equation}
    \begin{aligned}
      \sum_{K \in \mathcal{T}^H} \sum_{y_i \in K} & \left(  \left\| \uu-\left(\mathcal{I}_{L_i}^{\omega_i} \ww \right)_u  \right\|_{L_a,K} +  \left\| \theta-\left(\mathcal{I}_{L_i}^{\omega_i} \ww \right)_{\theta}  \right\|_{L_d,K}  \right)                                \\
                                                  & = \sum_{\omega_i} \left(  \left\| \uu-\left(\mathcal{I}_{L_i}^{\omega_i} \ww \right)_u  \right\|_{L_a,\omega_i} +  \left\| \theta-\left(\mathcal{I}_{L_i}^{\omega_i} \ww \right)_{\theta}  \right\|_{L_d,\omega_i} \right) \\
                                                  & = \sum_{\omega_i} \interleave \ww-\mathcal{I}_{L_i}^{\omega_i} \ww \interleave_{0,\omega_i}^2 \leq  \sum_{\omega_i} \frac{H^4}{\left( \Lambda_{L+1}^{\omega_i} \right)^{2}} \interleave  \ww  \interleave_{2,\omega_i}^2.
    \end{aligned}\label{eq:ukest}
  \end{equation}
  Then combining the summation of Eq. $\left( \ref{eq:local:interp} \right)$ for all $K \in \mathcal{T}_H$ with Eqs. (\ref{eq:akest}) and (\ref{eq:ukest}), the proof is complete.
\end{mypf}

\subsection{Steady state case}
Some results for the prior error estimate of the problem (\ref{eq:govequs}) will be first given in the steady state case. For $ \bar{\uu} \in V_u(\Omega)$ and $ \bar{\theta} \in V_{\theta}(\Omega)$, it follows that
\begin{equation}
  \begin{aligned}
    a(\bar{\uu}, \vv_u)-b(\vv_u, \bar{\theta}) & =\langle\bar{\ff}, \vv_u\rangle_a,     & \forall \vv_u \in V_u(\Omega),             \\
    d(\bar{\theta}, v_{\theta})                & =\langle\bar{g}, v_{\theta} \rangle_d, & \forall v_{\theta} \in V_{\theta}(\Omega).
  \end{aligned}\label{eq:steadyproblem}
\end{equation}
Define $V_{uH} = \text{span} \left\{\bPhi_{lu}^{\omega_i}\right\}$ and $V_{\theta H} = \text{span} \left\{  \Phi_{l\theta}^{\omega_i}\right\}$, where $ \bPhi_l^{\omega_i}(x) = \left( \bPhi_{lu}^{\omega_i}(x),\Phi_{l\theta}^{\omega_i}(x) \right)$, ($1\leq i \leq N_v, 1 \leq l \leq L_i $). The discretization form of Eq. (\ref{eq:steadyproblem}) is given as
\begin{equation}
  \begin{aligned}
    a(\bar{\uu}_H, \vv_{aH})-b(\vv_{uH}, \bar{\theta}_H) & =\langle\bar{\ff}, \vv_{uH}\rangle_a,    & \forall \vv_{uH} \in V_{uH},           \\
    d(\bar{\theta}_H, v_{\theta H})                      & =\langle\bar{g}, v_{\theta H} \rangle_d, & \forall v_{\theta H} \in V_{\theta H}.
  \end{aligned}
\end{equation}
Then for all $\ww = (\uu,\theta) \in V_u(\Omega) \times V_{\theta}(\Omega)$, define the Riesz projection operator $ \mathscr{R}_H = \left(  \mathscr{R}_{Hu}\left( u,\theta \right),  \mathscr{R}_{H \theta}\left(\theta \right) \right):  V_u(\Omega) \times V_{\theta}(\Omega) \rightarrow V_{uH} \times V_{\theta H}$, such that
\begin{equation} \nonumber
  \begin{aligned}
    a(\uu-\mathscr{R}_{Hu}\left( \uu,\theta \right), \vv_{uH})-b(\vv_{uH}, \theta-\mathscr{R}_{H \theta}\left(\theta \right)) & =0, & \forall \ \vv_{uH} \in V_{uH},           \\
    d(\theta-\mathscr{R}_{H \theta}\left(\theta \right), v_{\theta H})                                                        & =0, & \forall \ v_{\theta H} \in V_{\theta H}.
  \end{aligned}
\end{equation}
The following Lemma can be given with reference to \cite{EM:EMMNA-MMAN:2009}.
\begin{mylem}
  For all $\ww = (\uu,\theta) \in V_u(\Omega) \times V_{\theta}(\Omega) $, and the definition of $\mathscr{R}_H $, we have
  \begin{equation}
  \textcolor{black}{  \left\|  \uu-\mathscr{R}_{H u}(\uu,\theta)  \right\|_a } \leq \inf _{\vv_{uH} \in V_{uH}} \left\|  \uu-\vv_{uH}   \right\|_a + C_0 \left\|  \theta-\mathscr{R}_{H \theta} (\theta)  \right\|_c,\label{eq:cea_a}
  \end{equation}
  \begin{equation}
    \left\|  \theta-\mathscr{R}_{H \theta} (\theta)  \right\|_d \leq \inf _{v_{\theta H} \in V_{\theta H}} \left\|  \theta-v_{\theta H}   \right\|_d.
  \end{equation}
\end{mylem}

\begin{mylem}
  For all $r \in L_d$, and let $ \phi \in V_{\theta}(\Omega) $ be the solution of the dual problem $ d(\phi, v_{\theta})=c(r, v_{\theta}), \forall\ v_{\theta} \in V_{\theta}(\Omega)$, and  $\phi_H \in V_{\theta H}$ be the solution of the discrete problem $ d(\phi_H, v_{\theta})=c(r, v_{\theta}), \forall v_{\theta} \in V_{\theta H}$. There holds
  \begin{equation}
    \left\|  \phi -  \phi_H \right\|_d \leq C_4 \left\|  r \right\|_c.
  \end{equation}
  where $ C_4=C_p^{\frac{1}{2}} \kappa_{\text{min} }^{-\frac{1}{2}} $, and $ C_p $ is the Poincaré constant.
\end{mylem}
\begin{mypf}
  The following equation can first be estimated,
  \begin{equation}\label{eq:phierr1}
    \left\|  \phi -  \phi_H \right\|_d^2 =   d\left( \phi, \phi -  \phi_H \right) = c \left( r, \phi -  \phi_H \right) \leq 	\left\|  r \right\|_c 	\left\|  \phi -  \phi_H \right\|_c.
  \end{equation}
  \begin{equation}
    C_p^{-1}\kappa_{\text{min} }	\left\|  \phi -  \phi_H \right\|_c^2 \leq \kappa_{\text{min} }	\left\|  \nabla \left(\phi -  \phi_H \right) \right\|_c^2 \leq 	\left\|  \phi -  \phi_H \right\|_d^2.\label{eq:est:C4}
  \end{equation}
  Combining Eq. (\ref{eq:phierr1}) with Eq. (\ref{eq:est:C4}), it follows that
  \begin{equation} \nonumber
    \left\|  \phi -  \phi_H \right\|_d^2 \leq C_p^{\frac{1}{2}}\kappa_{\text{min} }^{-\frac{1}{2}}	\left\|  r \right\|_c 	\left\|  \phi -  \phi_H \right\|_d.
  \end{equation}
\end{mypf}

\begin{mylem}
  For all $\ww = (\uu,\theta) \in V_u(\Omega) \times V_{\theta}(\Omega) $, and the definition of $\mathscr{R}_H $, we have
  \begin{equation}
    \left\|  \uu-\mathscr{R}_{Hu}(\uu,\theta)  \right\|_a \leq \text{max} \{ 1,C_4 C_0\} C\left(H,\Lambda^{\omega_i}_{L+1}\right) 	\interleave  \ww \interleave_{2,\Omega},\label{eq:interpest_a}
  \end{equation}
  \begin{equation} \label{eq_interpest_d}
    \left\|  \theta-\mathscr{R}_{H \theta} (\theta)  \right\|_d \leq  C\left(H,\Lambda^{\omega_i}_{L+1}\right) \interleave \ww \interleave_{2,\Omega},
  \end{equation}
  \begin{equation}
    \left\|  \theta-\mathscr{R}_{H \theta} (\theta)  \right\|_c \leq  C_4 C\left(H,\Lambda_{L+1}\right) \interleave \ww \interleave_{2,\Omega},\label{eq:interpest_c}
  \end{equation}
  where
  \begin{equation} \label{eq_constant_HLambda}
    C\left(H,\Lambda_{L+1}\right) =  \left[ C_2 \frac{H^2}{ \Lambda_{L+1} } +  \left( 2C_1^2 H^2+C_3 H^4  \right) \frac{1}{\left( \Lambda_{L+1} \right)^{2}} \right]^{\frac{1}{2}},
  \end{equation}
  and
  $\ds\interleave  \ww  \interleave_{s,\Omega}^2 =  \sum_{\omega_i \in \Omega} \interleave \ww \interleave_{s,\omega_i}^2$. 
\end{mylem}
\begin{mypf}
  For all $\theta \in V_{\theta}(\Omega) $, and let $ r:= \theta-\mathscr{R}_{Hd} (\theta) $, it follows
  \begin{equation} \nonumber
    \begin{aligned}
      \left\|  \theta-\mathscr{R}_{H \theta} (\theta)  \right\|_c^2 & = c(r,r) = d(\phi,r) = d(\phi - \phi_H,r)                                                                            \\
        & \leq  \left\| \phi -  \phi_H \right\|_d \left\|  r \right\|_d \leq C_4 \left\|  r \right\|_c  \left\|  r \right\|_d.
    \end{aligned}
  \end{equation}
  Then we have
  \textcolor{black}{
  \begin{equation}
    \left\|\theta-\mathscr{R}_{H \theta}(\theta)\right\|_c \leq C_4\left\|\theta-\mathscr{R}_{H \theta}(\theta)\right\|_d.
    \label{eq:interp_cd}
  \end{equation}
  }
  Moreover, combining Eq. (\ref{eq:cea_a}) with Eq. (\ref{eq:interp_cd}), it follows
  \begin{equation} \nonumber
    \left\|  \uu-\mathscr{R}_{Hu}(\uu,\theta)  \right\|_a \leq \inf _{\vv_{uH} \in V_{uH}} \left\| \uu-\vv_{uH}   \right\|_a + C_4 C_0 \inf _{v_{\theta H} \in V_{\theta H}} \left\|  \theta-v_{\theta H}   \right\|_d.
  \end{equation}
  By the definition of $\ds\interleave  \ww  \interleave_{s,\Omega}^2$, the proof is complete.
\end{mypf}

\subsection{The prior error estimate of CGMsFEM}
For $ \ff $ and $ g $ in problem (\ref{eq:govequs}), we have
\begin{equation} \nonumber
  \langle \tilde{\ff}, v_u\rangle_a = \langle \ff, v_u\rangle, \forall \ v_u \in V_u(\Omega), \quad \langle \tilde{g}, v_{\theta}\rangle_d = \langle g, v_{\theta}\rangle , \forall \ v_{\theta} \in V_{\theta}(\Omega),
\end{equation}
where $ \tilde{\ff} \in V_u(\Omega)$ and $ \tilde{g} \in V_{\theta}(\Omega)$. Define
\begin{equation}  \nonumber
  C^n(\ff,g) = \frac{1}{2} \left\| \tilde{\ff}^n-\tilde{\ff}^n_H \right\|_{a}^2 + \tau_n  \left\| \tilde{g}^n-\tilde{g}^n_H \right\|_{d}^2,\quad \forall \ n \in {1,\cdots,N_T}.
\end{equation}
For simplicity of notation, let
\begin{equation}  \nonumber
  C_1^n(\ww) = 4 C_4^2\left(C_4^2+C_0^2 \text{max}\{ 1,C_4^2 C_0^2\} \right) \left\| \partial_t \ww  \right\|_{L^{\infty}\left(T_n, \interleave \cdot \interleave_{2,\Omega} \right)},
\end{equation}
\begin{equation}  \nonumber
  C_2^n(\ww) = 2 C_4^2\left( \left\| \partial_{tt} \theta \right\|_{L^{\infty}\left(T_n, \| \cdot \|_{c} \right)} + C_0^2 \left\| \partial_{tt} \uu \right\|_{L^{\infty}\left(T_n, \| \cdot \|_{a} \right)} \right).
\end{equation}
Then we have the following prior error estimate of the CGMsFEM.
\begin{mythm}\label{theorem:energyerror}
  Let $ \ww = \left(\uu,\theta\right) $ and $ \ww_H = \left(\uu_H,\theta_H\right) $ be the unique solution and CGMsFEM solution of problem ($ \ref{eq:govequs} $), there holds
  \begin{equation}
    \begin{aligned}
      \frac{1}{4} \left\| \uu^n - \uu^n_H \right\|_a^2 + \frac{1}{4} \left\| \theta^n - \theta^n_H \right\|_c^2 \leq  \sum_{m=1}^{n} \left[ C^m\left(\ff,g\right) + \tau_m C^2\left(H,\Lambda_{L+1}\right) C_1^m(\ww) \right. \\
        + \left. \tau_m^3 C_2^m(\ww) \right]
      + \frac{1}{2} \left(C_4^2 + \text{max} \left\{ 1,C_4^2 C_0^2 \right\}\right) C^2\left(H,\Lambda_{L+1}\right) \interleave \ww^n \interleave_{2,\Omega}^2,
    \end{aligned}
  \end{equation}
  and
  \begin{equation}
    \begin{aligned}
      \sum_{m=1}^{n} \frac{1}{8} \tau_m  \left\| \theta^m - \theta^m_H \right\|_d^2 \leq  \sum_{m=1}^{n} \left[ C^m\left(\ff,g\right) + \tau_m C^2\left(H,\Lambda_{L+1}\right) C_1^m(\ww) \right. \\
        \left. + \tau_m^3 C_2^m(\ww) \right]
      + \sum_{m=1}^{n} \frac{1}{4} \tau_m C^2\left(H,\Lambda_{L+1}\right) \interleave \ww^n \interleave_{2,\Omega}^2.
    \end{aligned}
  \end{equation}
  where $n \in \left\{1,2,\cdots,N_T\right\} $.
\end{mythm}
The proof of this theorem is similar to Theorem 3.1 in Reference \cite{EM:EMMNA-MMAN:2009}. Those interested can refer to the original source or find a brief imitation in  \ref{sec:appendixB}. 
Excluding terms unrelated to $H$ and $\Lambda_{L+1}$, it can be concluded that the energy error of the CGMsFEM solution is primarily controlled by $C^2\left(H,\Lambda_{L+1}\right)$ defined in Eq. (\ref{eq_constant_HLambda}).
When a sufficient number of basis functions are selected for each subdomain $\omega_i$, which implies that $\Lambda_{L+1}$ is sufficiently large, we have $C^2\left(H,\Lambda_{L+1}\right) = O(\frac{H^2}{ \Lambda_{L+1} })$.

\section{Numerical experiments}
\label{sec:numexp}
In this section, we present several numerical examples to evaluate the performance of the proposed CGMsFEM for the thermomechanical problem (\ref{eq:govequs}) and mainly focus on the verification of the \textcolor{black}{accuracy} and efficiency of the CGMsFEM. The computation domain is $\Omega=[0,1]^2$, and the time interval is $(0,1]$. The Dirichlet condition and Neumann boundary conditions for the displacement and temperature fields are defined as 0 in $\Gamma_D^u=\Gamma_D^{\theta}=[0,1]\times0$ and $\Gamma_N^u=\Gamma_N^{\theta}=\partial\Omega\backslash \Gamma_D^u$, respectively. The solutions $(\uu^{cgm},\theta^{cgm})$ of CGMsFEM will be compared with the reference solutions $(\uu^{ref},\theta^{ref})$ of the standard finite element method and the solutions $(\uu^{gm},\theta^{gm})$ of GMsFEM. Then, based on Theorem \ref{theorem:energyerror}, the relative energy errors of each solution and the total relative energy errors are defined as follows
\begin{equation}
  \begin{array}{lll}\vspace{2mm}
    ||E_u||_{e}=\frac{\bigg(\ds \int_{\Omega}\sigma(E_u):\epsilon(E_u) dx\bigg)^\frac{1}{2}}{\bigg(\ds \int_{\Omega}\sigma(\uu):\epsilon(\uu) dx\bigg)^\frac{1}{2}},\quad
    ||E_\theta||_{e}=\frac{\bigg(\ds \int_{\Omega}\kappa\nabla E_\theta\cdot\nabla E_\theta dx\bigg)^\frac{1}{2}}{\bigg(\ds\int_{\Omega}\kappa\nabla \theta\cdot\nabla \theta dx\bigg)^\frac{1}{2}}, \\
    ||E_w||_{e}=\frac{\bigg(\ds \int_{\Omega}\sigma(E_u):\epsilon(E_u)+\kappa\nabla E_\theta\cdot\nabla E_\theta dx\bigg)^\frac{1}{2}}{\bigg(\ds\int_{\Omega}\sigma(\uu):\epsilon(\uu)+\kappa\nabla \theta\cdot\nabla \theta dx\bigg)^\frac{1}{2}}.
  \end{array}\label{eq:energyerror}
\end{equation}
$E_u$ represents the \textcolor{black}{displacement} error $E_u^{cgm}=\uu^{cgm}-\uu^{ref}$ or $E_u^{gm}=\uu^{gm}-\uu^{ref}$, and $E_\theta$ represents the \textcolor{black}{temperature} error $E_\theta^{cgm}=\theta^{cgm}-\theta^{ref}$ or $E_\theta^{gm}=\theta^{gm}-\theta^{ref}$.
Our numerical experiments are performed on a desktop workstation with 16 GB of memory and a 3.4 GHz Core i7 CPU.

\subsection{Verification of CGMsFEM with periodic microstructure}
\label{subsec:verification}
In this example, the efficiency and accuracy of the CGMsFEM are verified for heterogeneous media with periodic microstructure. The body force $\ff$ and heat source $g$ are chosen as
$$\ff(x,y)=0,\quad g(x,y)=10.$$
The initial boundary condition is defined as $$ \theta_0(x,y)=500x(1-x)y(1-y). $$
\begin{figure}[h]
  \centering
  {\tiny(a)}\includegraphics[width=0.45\textwidth]{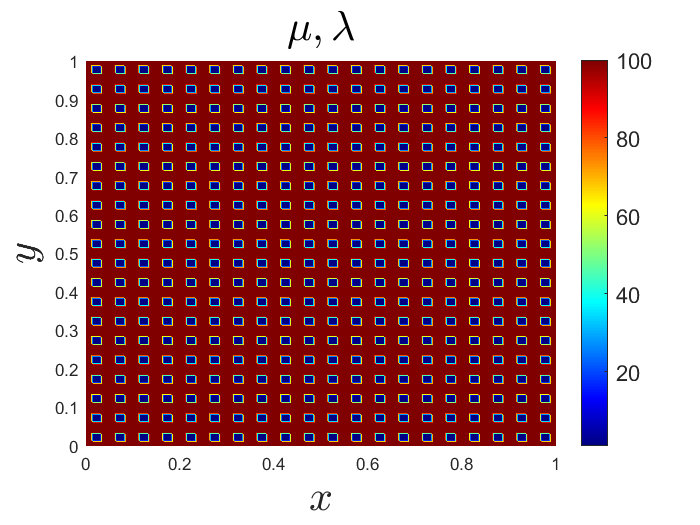}~
  {\tiny(b)}\includegraphics[width=0.45\textwidth]{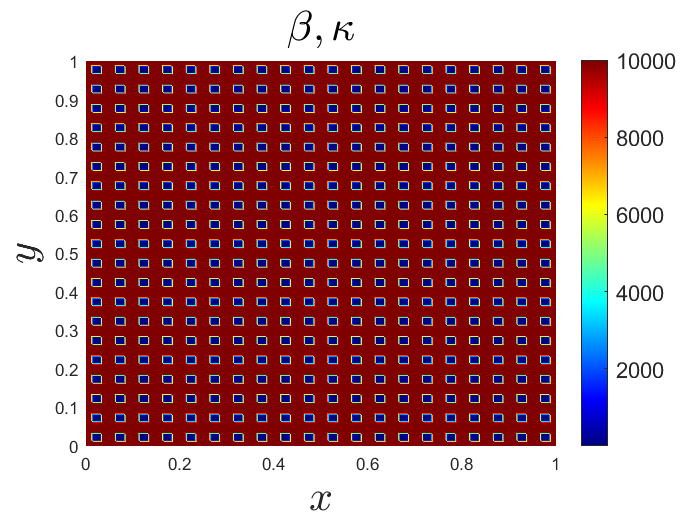}
  \caption{
    Contour plots of the material coefficients in periodic microstructure.
    (a) Lam$\acute{e}$ coefficients $\mu$, and $\lambda$;
    (b) Thermal conductivity coefficient $\kappa$, and expansion coefficient $\beta$.
  }\label{sec5-fig1}
\end{figure}

Then the Lam$\acute{e}$ coefficients $\mu$, $\lambda$, thermal conductivity coefficient $\kappa$, and expansion coefficient $\beta$ are shown in Figure \ref{sec5-fig1}, where the contrasts are chosen as $\lambda _{\text{max}}:\lambda_ {\text{min}}=10^2:1$, $\mu_{\text{max}}:\mu_{\text{min}}=10^2:1$, $\kappa_{\text{max}}:\kappa_{\text{min}}=10^4:1$,$\beta_{\text{max}}:\beta_{\text{min}}=10^4:1$. The corresponding relaxation coefficients $\gamma_1=0.4,\gamma_2=0.04$ are given, and the time step is $\tau=0.02$. Here, the $200\times200$ fine grid is used for the reference solution, and the $20\times 20$ coarse grid is used for the proposed CGMsFEM and GMsFEM. The number of local coupling multiscale basis functions for CGMsFEM is fixed at 8, and the total number of GMsFEM multiscal basis functions for displacement $\uu$ and $\theta$ is also chosen at 8.

\begin{figure}[htbp]
  \centering
  {\tiny(1)}\includegraphics[width=0.15\textwidth]{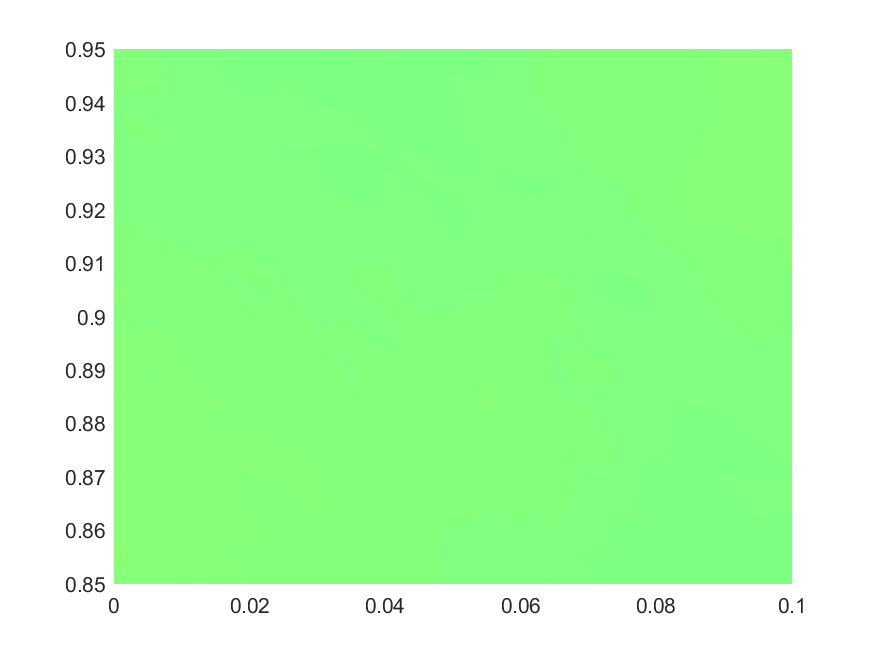}
  \includegraphics[width=0.15\textwidth]{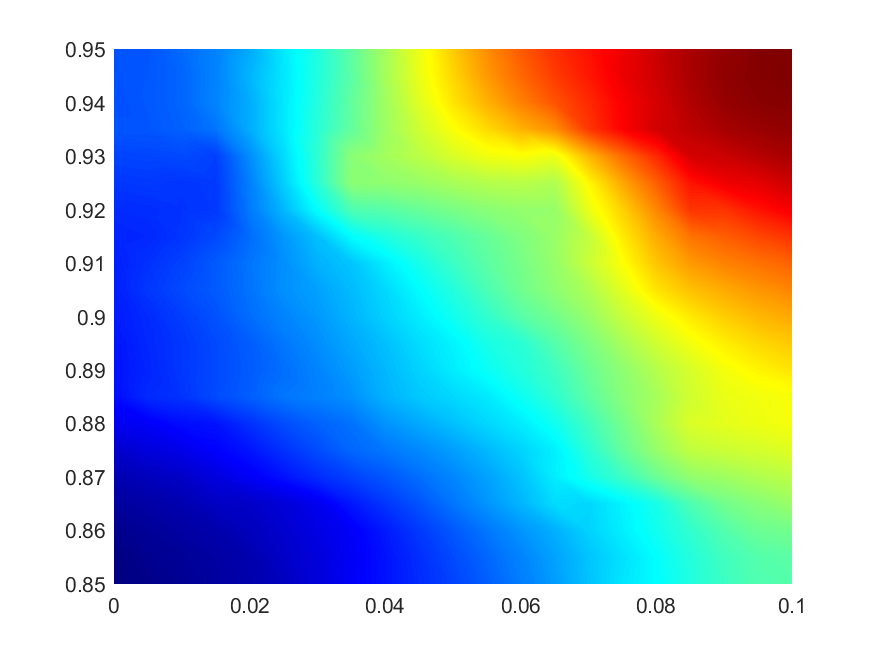}
  \includegraphics[width=0.15\textwidth]{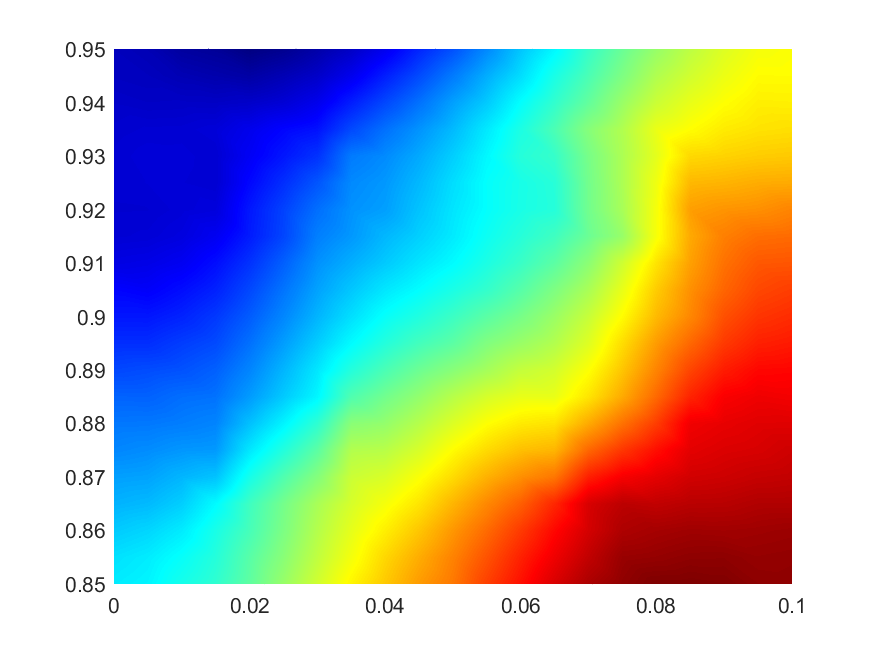}
  \includegraphics[width=0.15\textwidth]{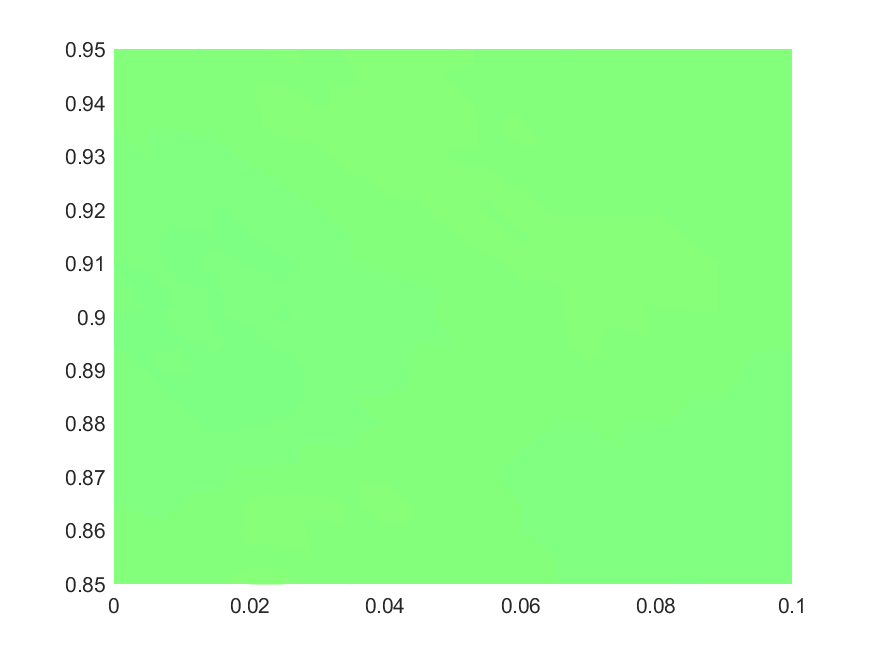}
  \includegraphics[width=0.15\textwidth]{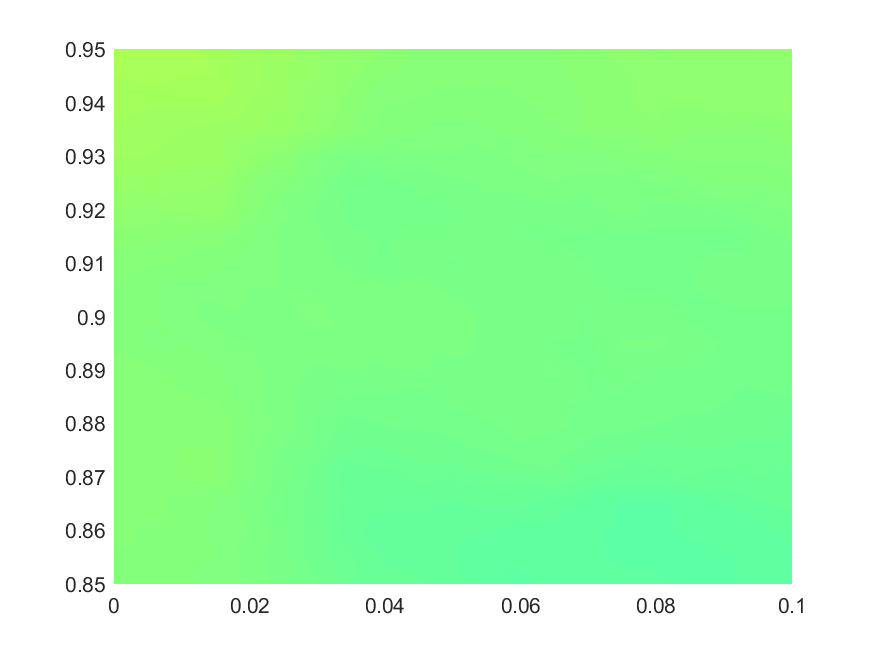}
  \includegraphics[width=0.15\textwidth]{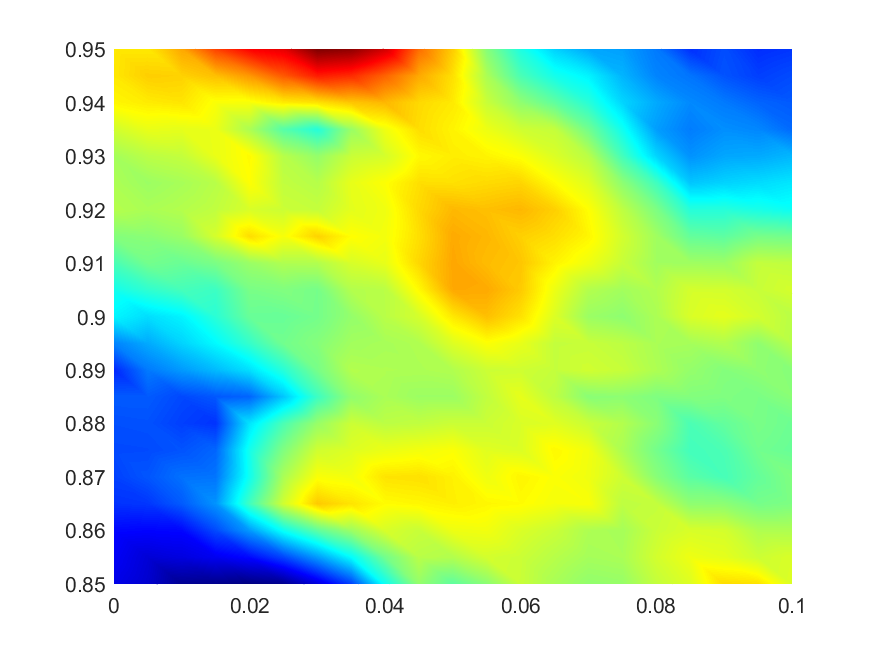}
  {\tiny(2)}\includegraphics[width=0.15\textwidth]{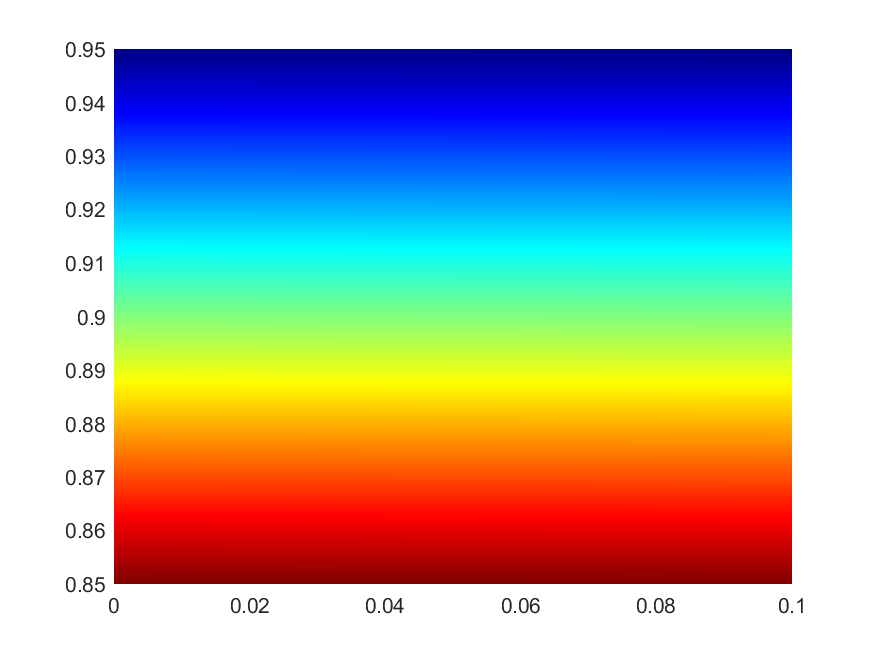}
  \includegraphics[width=0.15\textwidth]{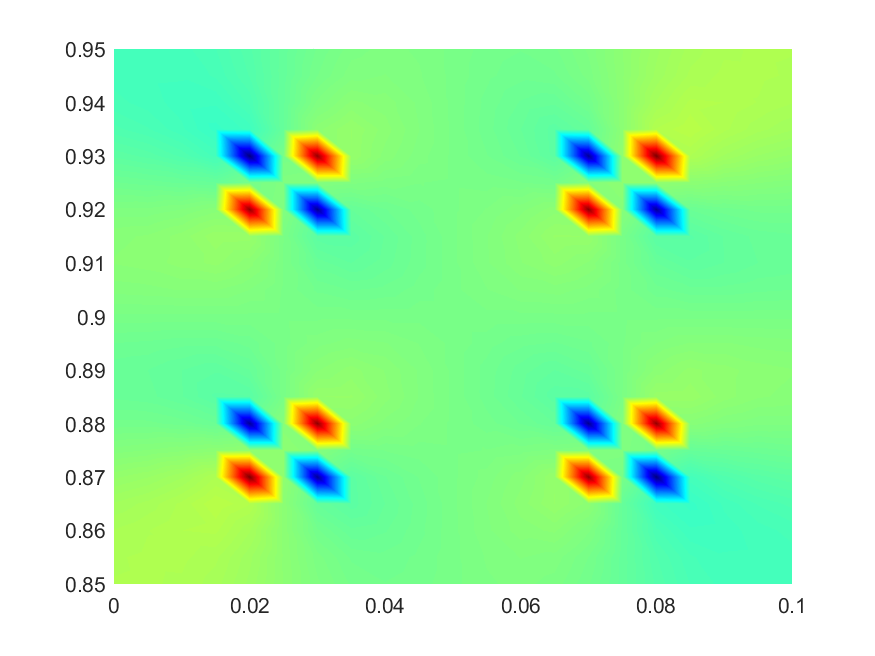}
  \includegraphics[width=0.15\textwidth]{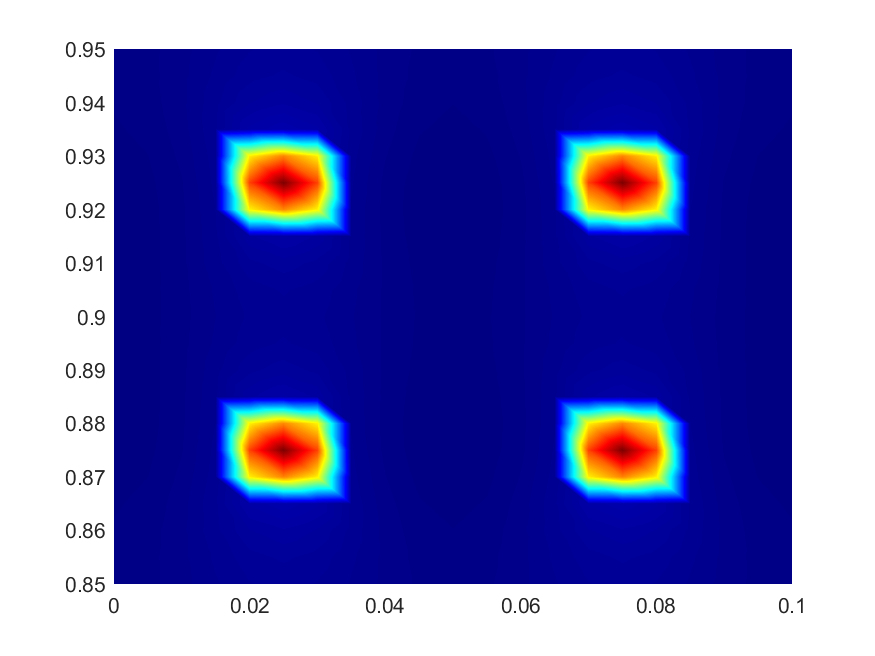}
  \includegraphics[width=0.15\textwidth]{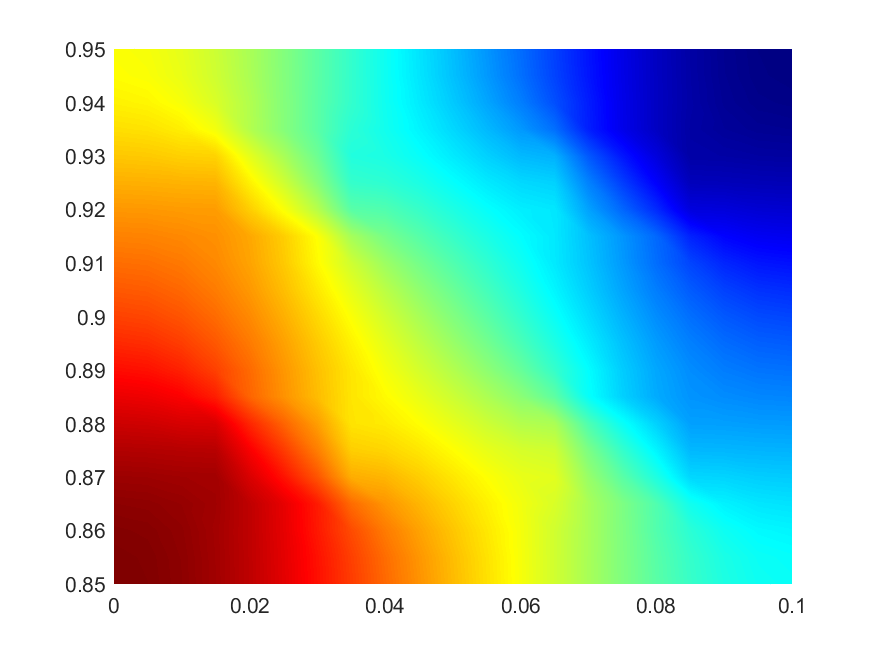}
  \includegraphics[width=0.15\textwidth]{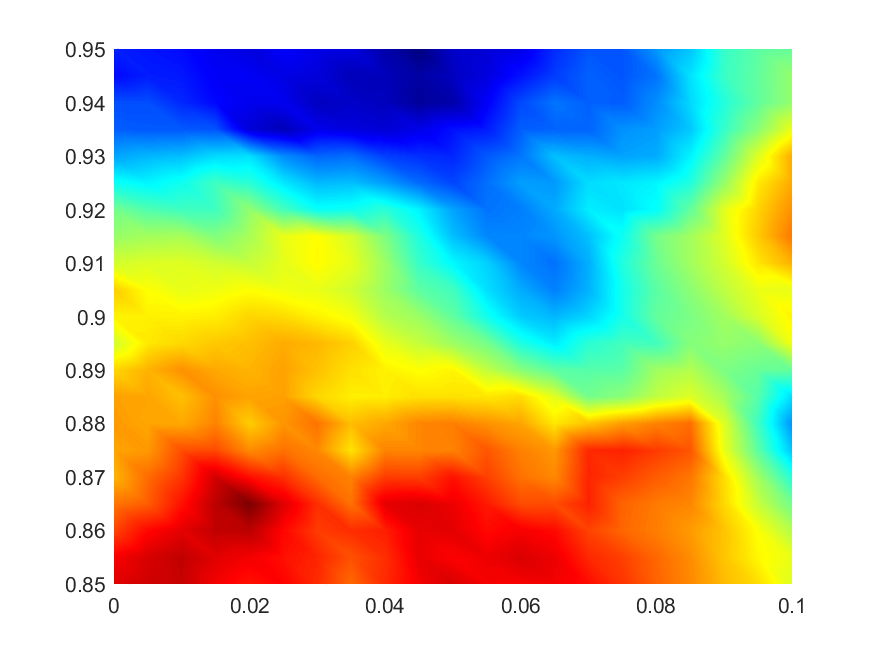}
  \includegraphics[width=0.15\textwidth]{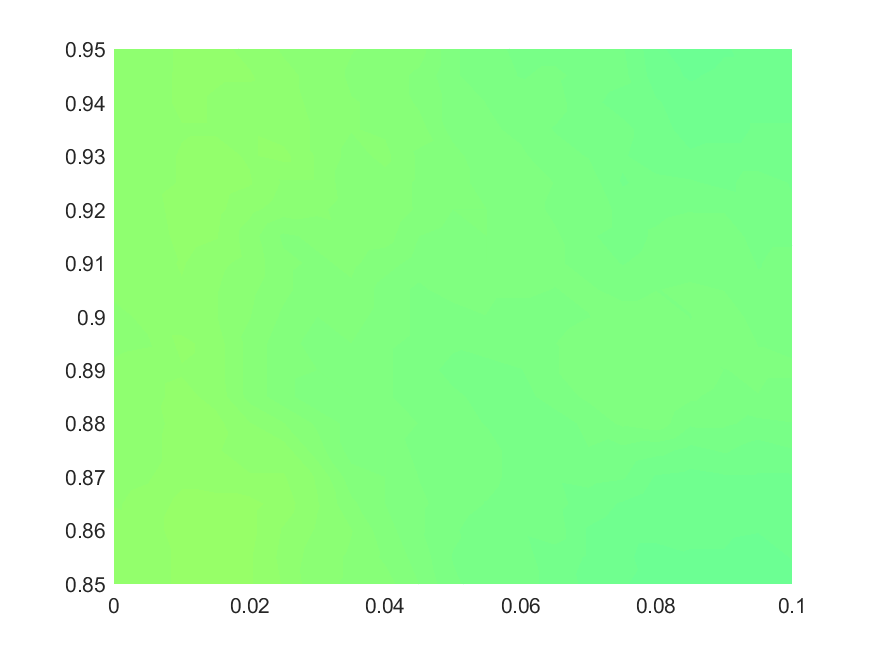}
  {\tiny(3)}\includegraphics[width=0.15\textwidth]{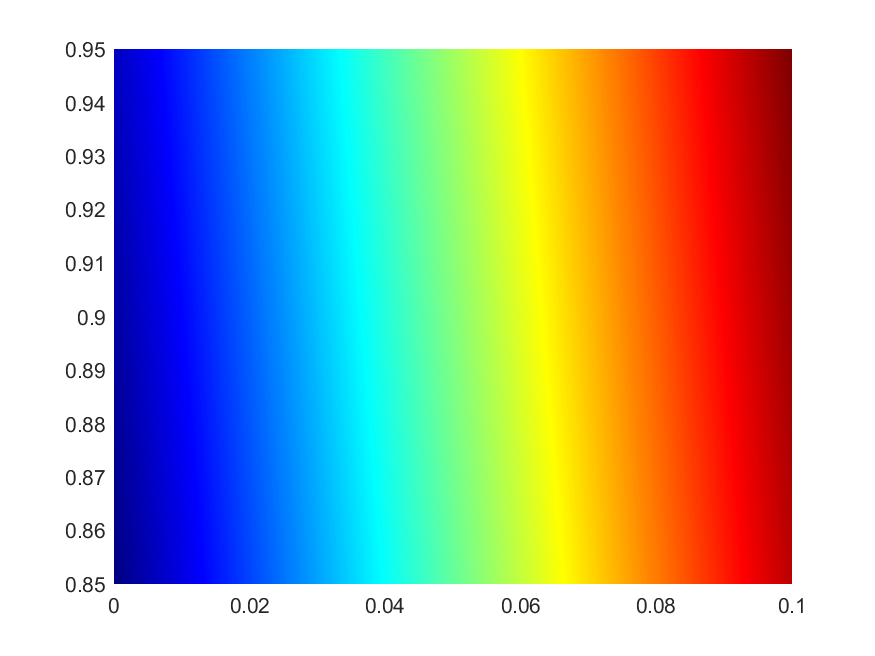}
  \includegraphics[width=0.15\textwidth]{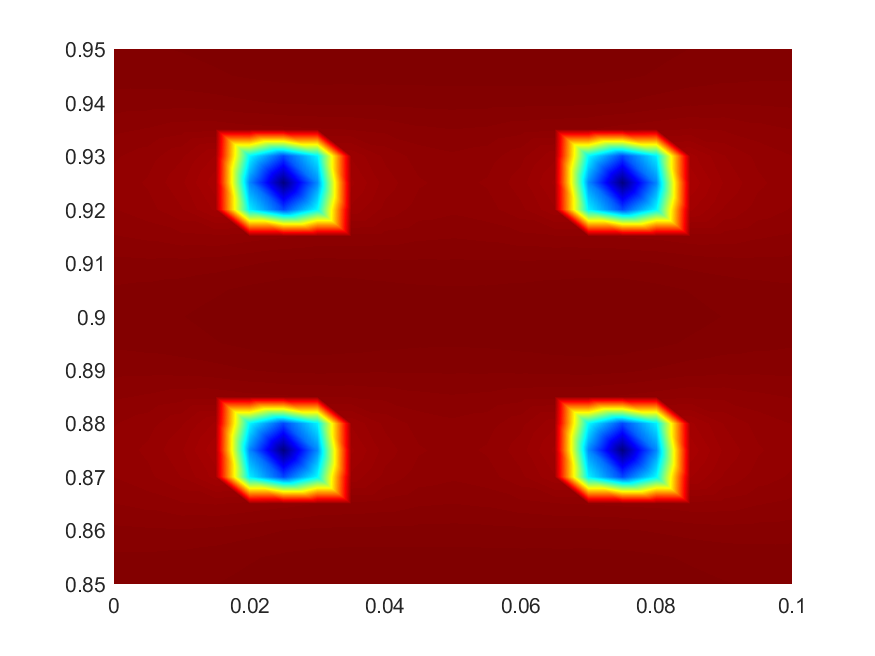}
  \includegraphics[width=0.15\textwidth]{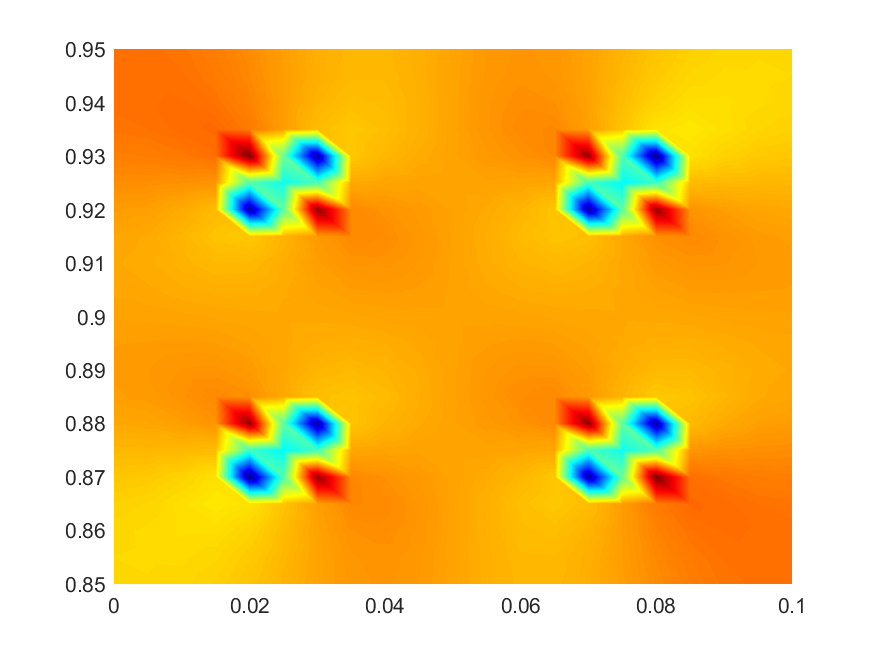}
  \includegraphics[width=0.15\textwidth]{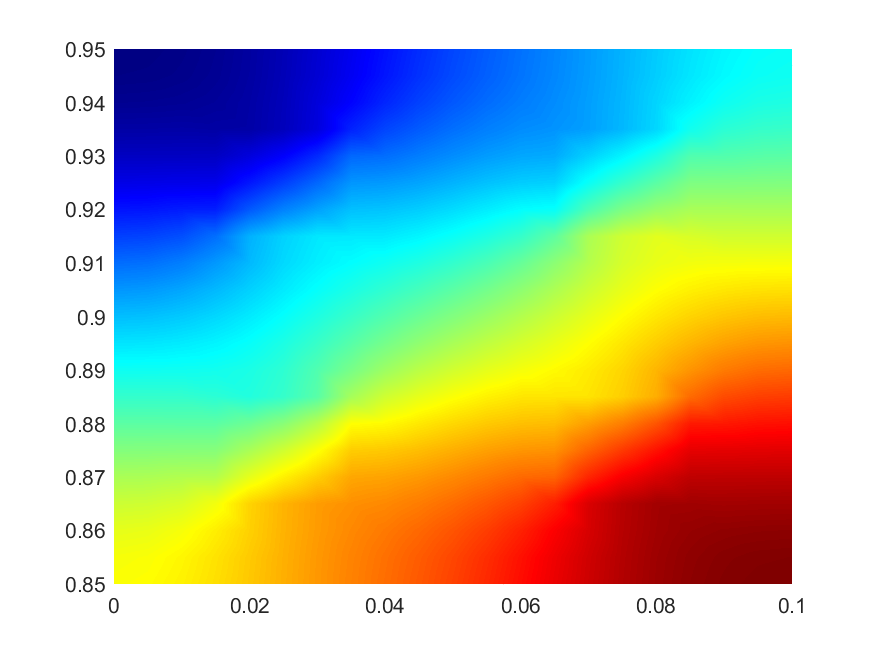}
  \includegraphics[width=0.15\textwidth]{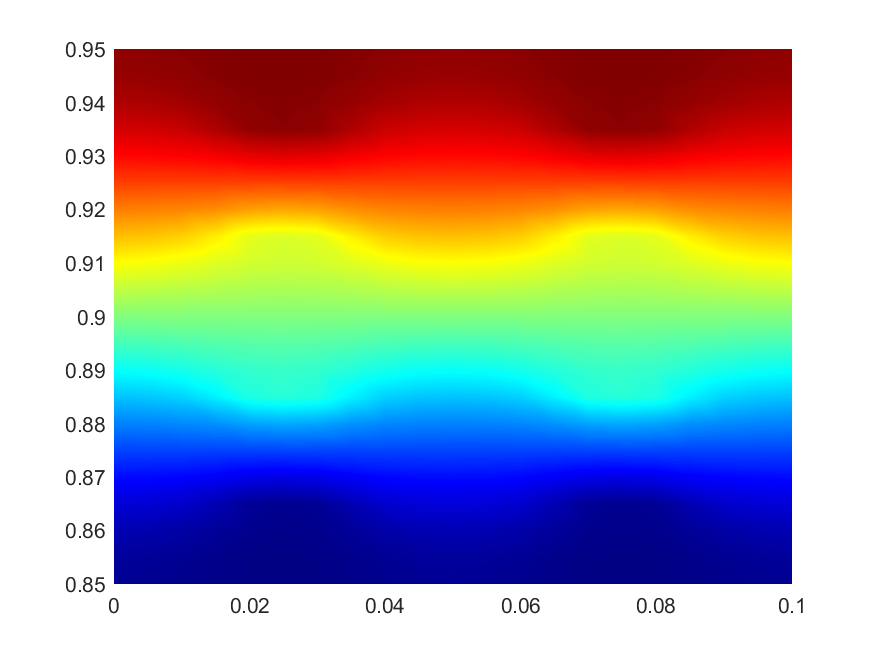}
  \includegraphics[width=0.15\textwidth]{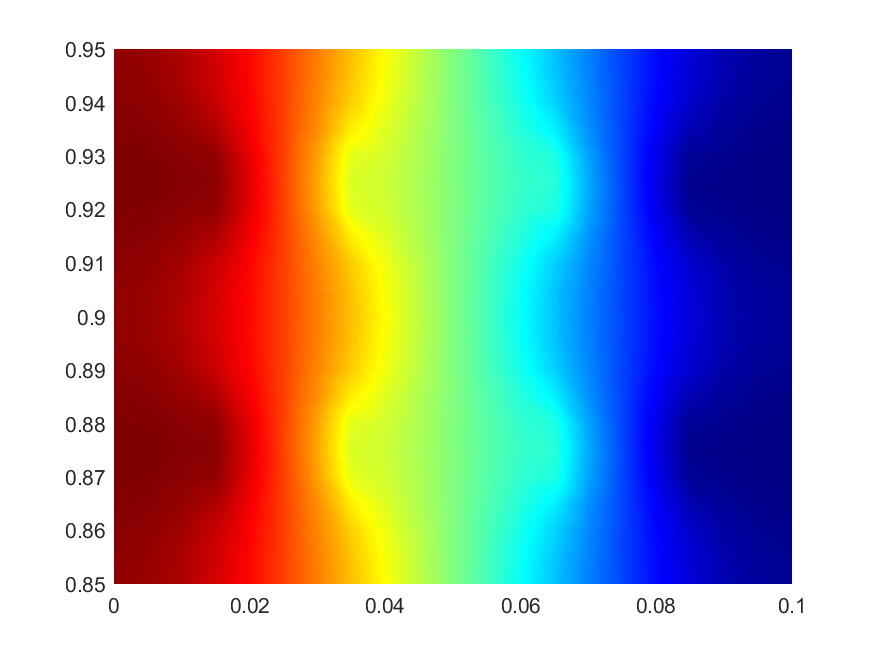}
  {\tiny(4)}\includegraphics[width=0.15\textwidth]{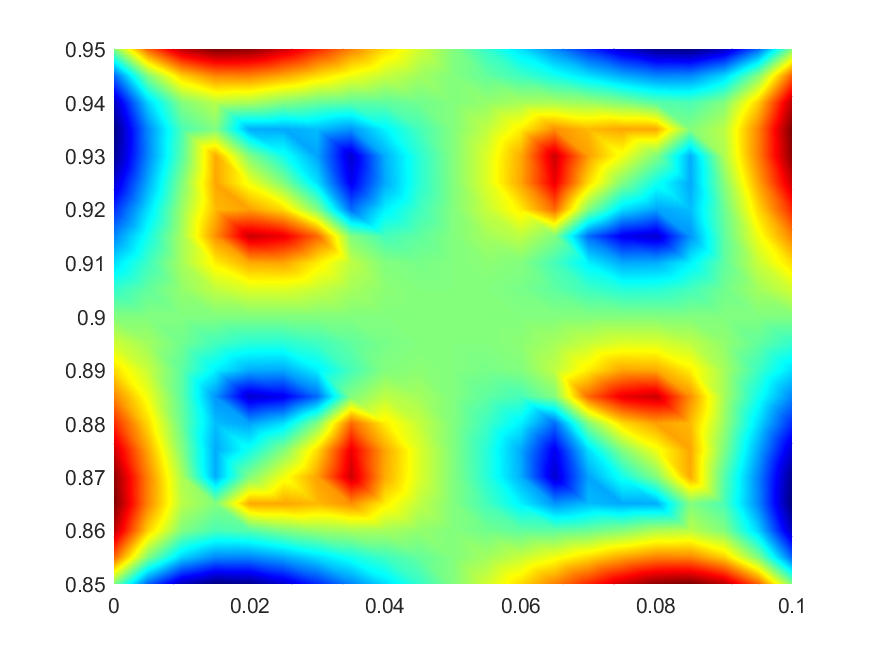}
  \includegraphics[width=0.15\textwidth]{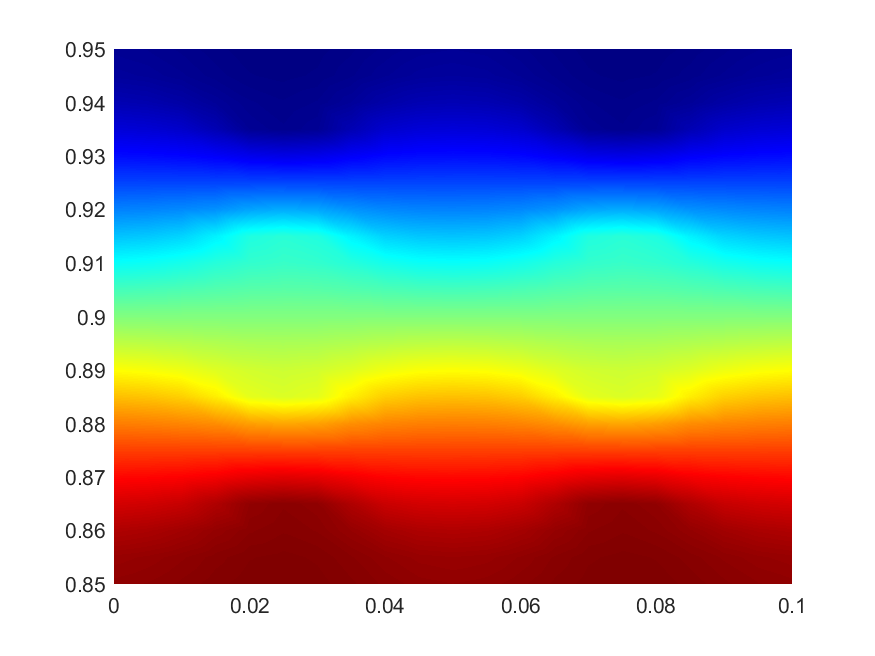}
  \includegraphics[width=0.15\textwidth]{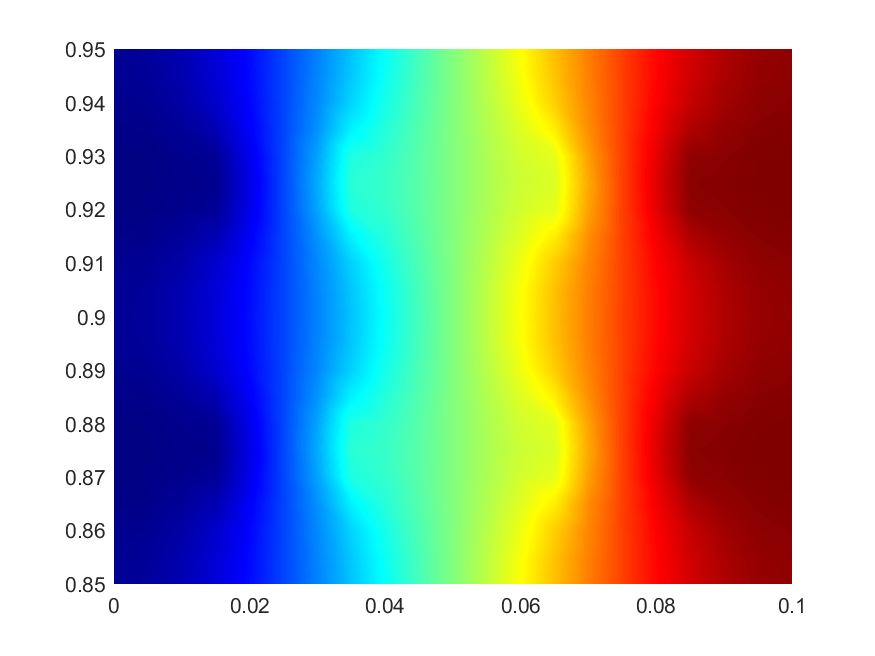}
  \includegraphics[width=0.15\textwidth]{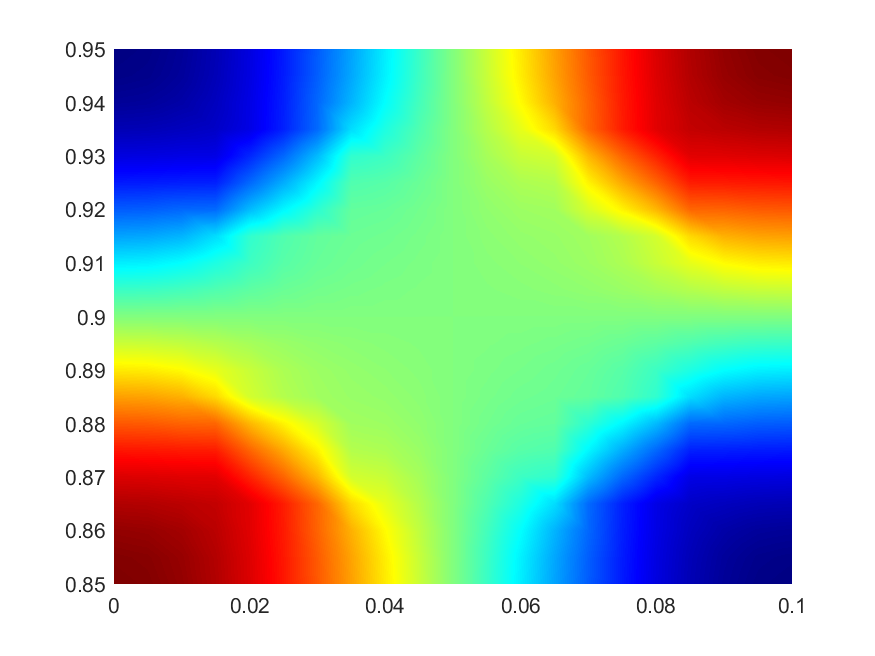}
  \includegraphics[width=0.15\textwidth]{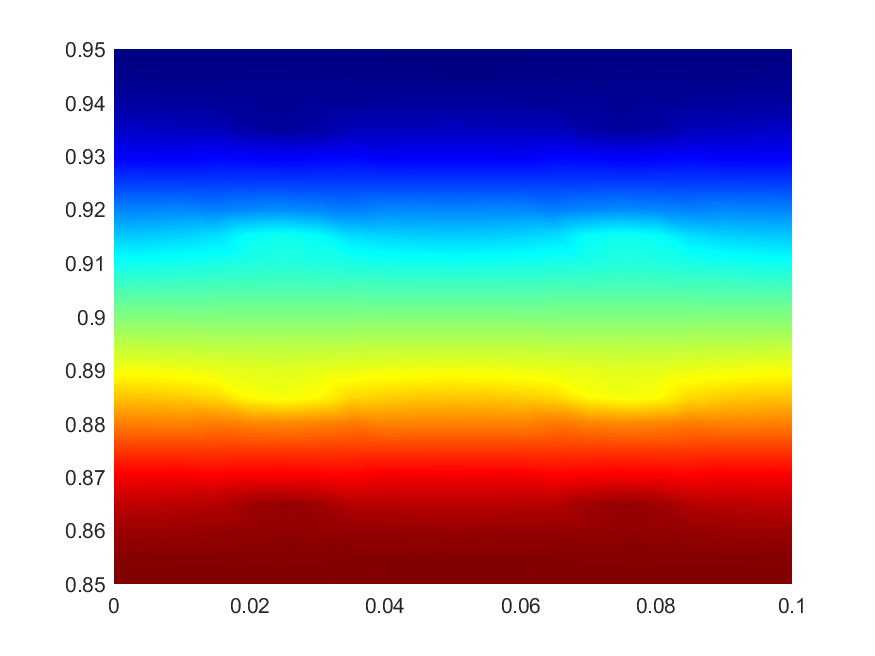}
  \includegraphics[width=0.15\textwidth]{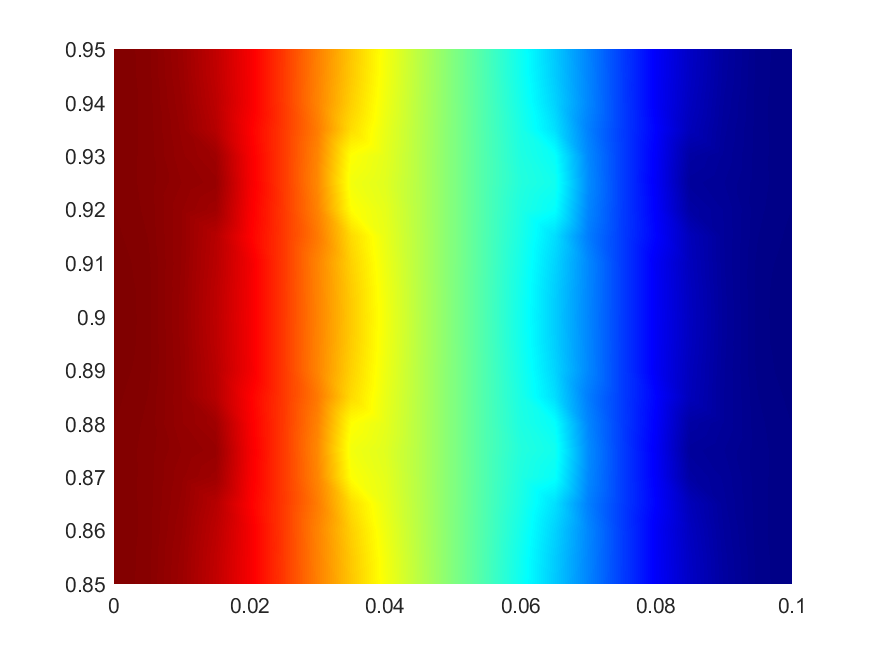}
  {\tiny(5)}\includegraphics[width=0.15\textwidth]{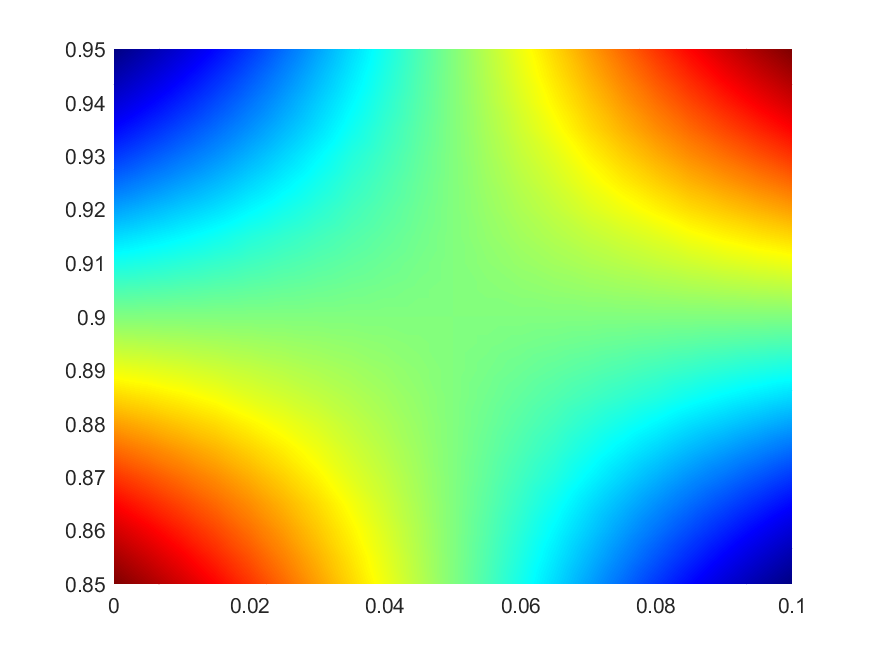}
  \includegraphics[width=0.15\textwidth]{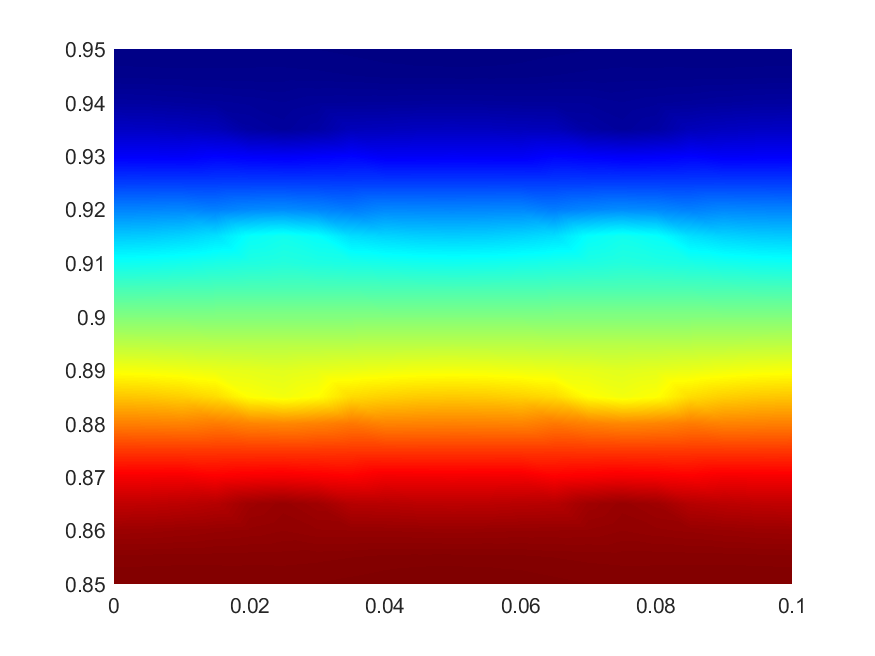}
  \includegraphics[width=0.15\textwidth]{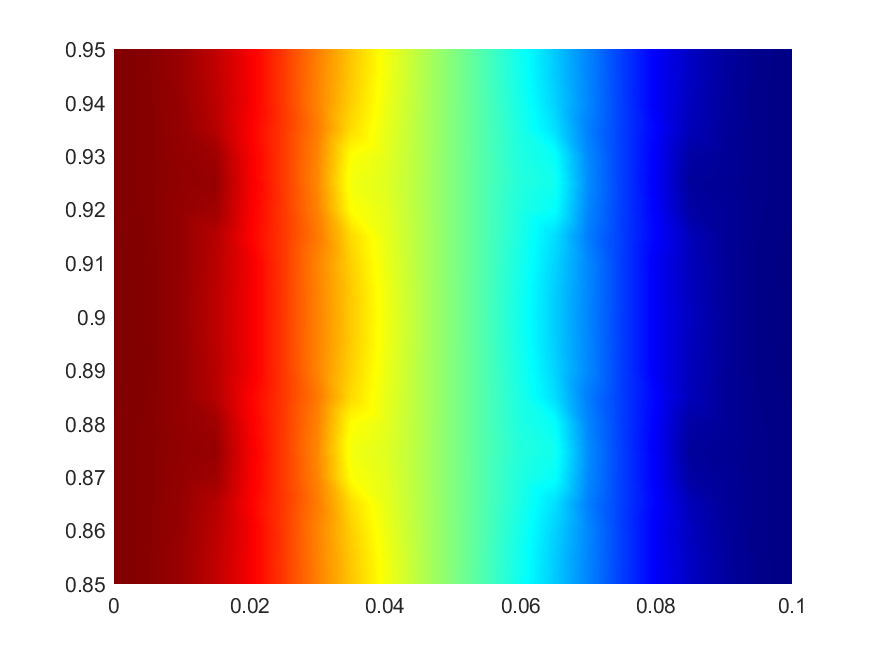}
  \includegraphics[width=0.15\textwidth]{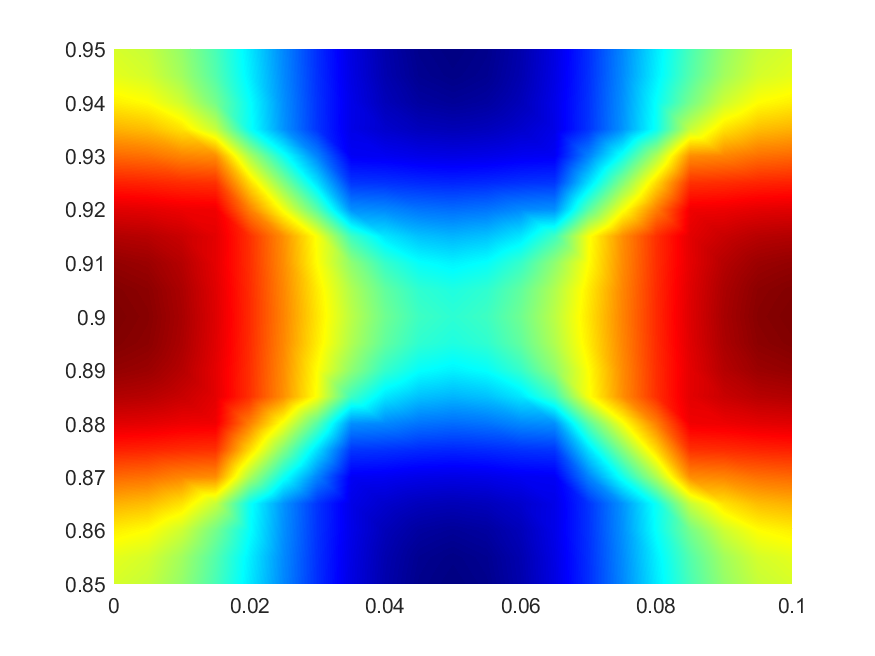}
  \includegraphics[width=0.15\textwidth]{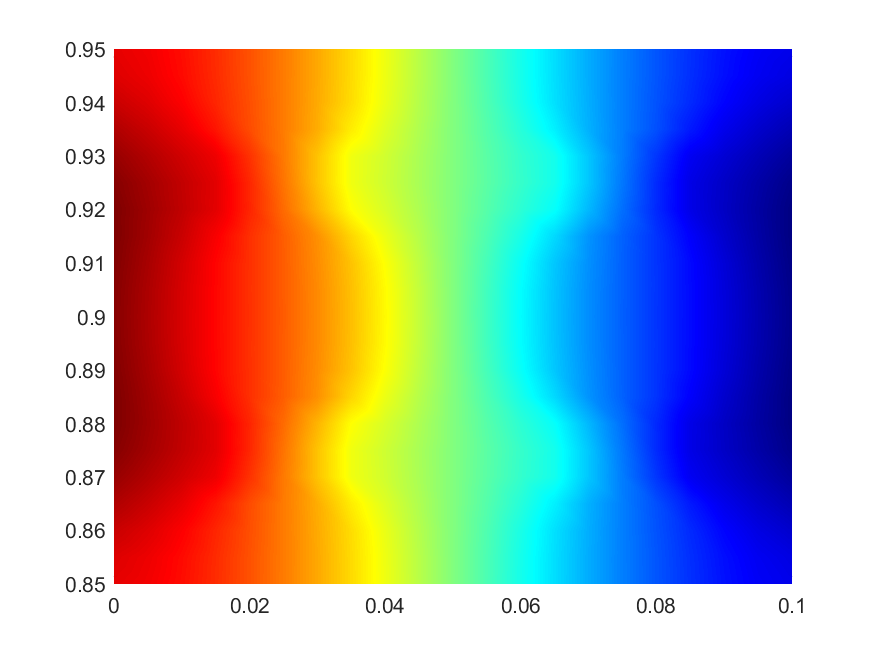}
  \includegraphics[width=0.15\textwidth]{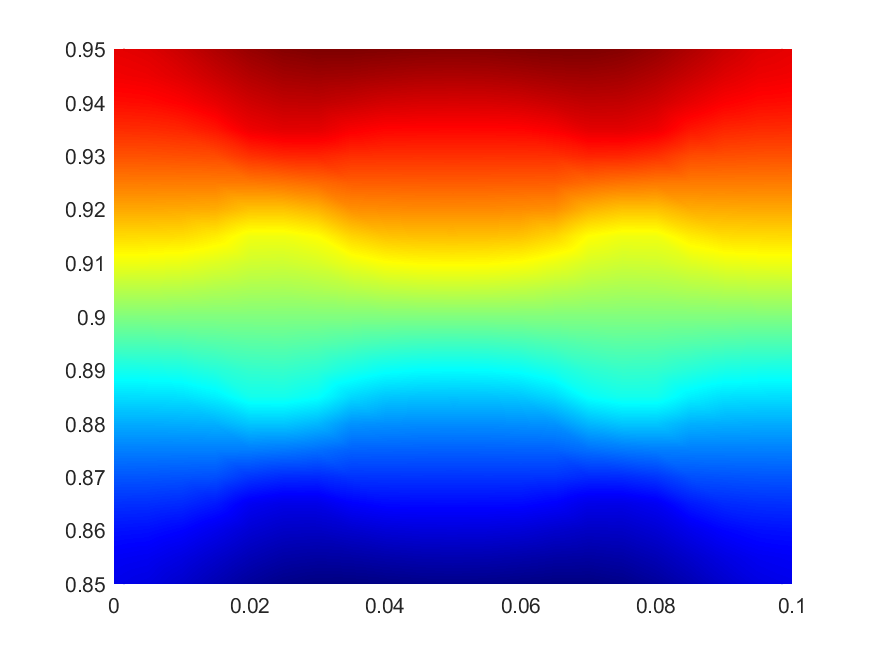}
  {\tiny(6)}\includegraphics[width=0.15\textwidth]{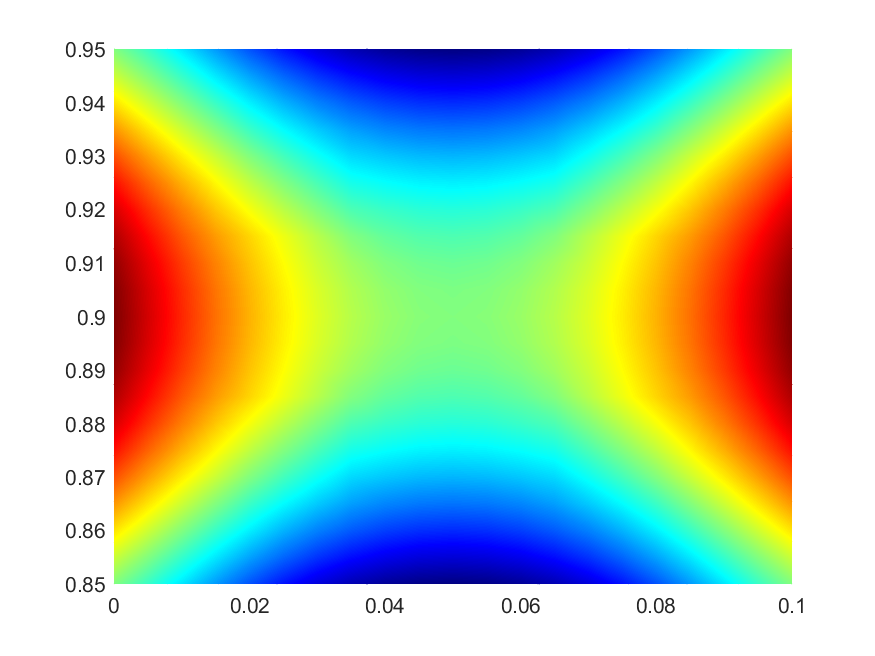}
  \includegraphics[width=0.15\textwidth]{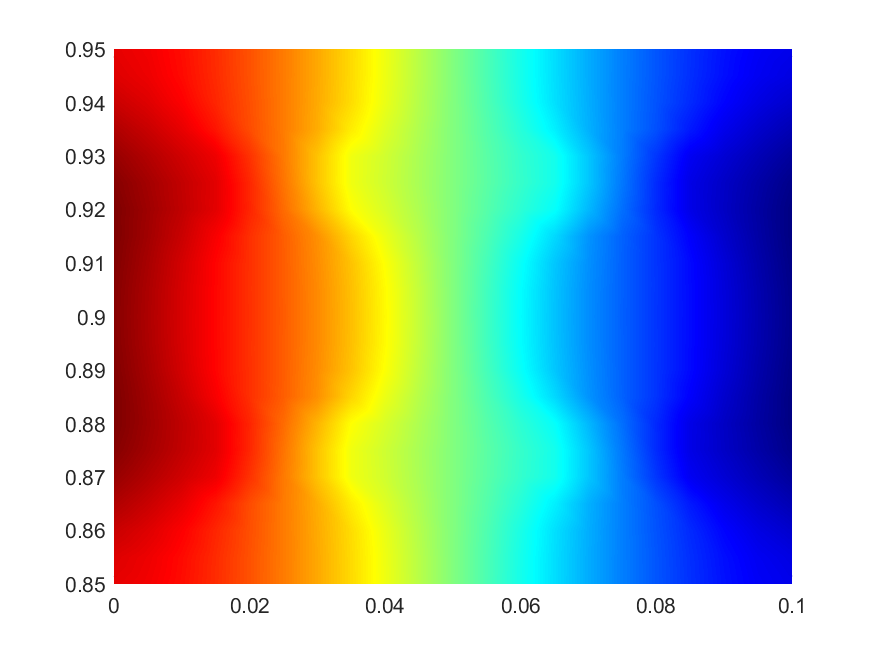}
  \includegraphics[width=0.15\textwidth]{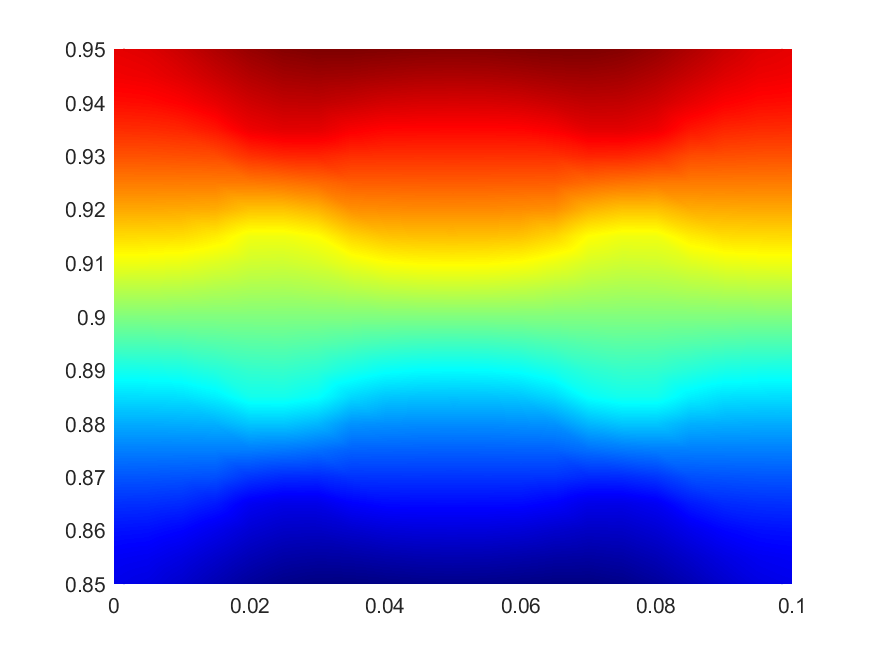}
  \includegraphics[width=0.15\textwidth]{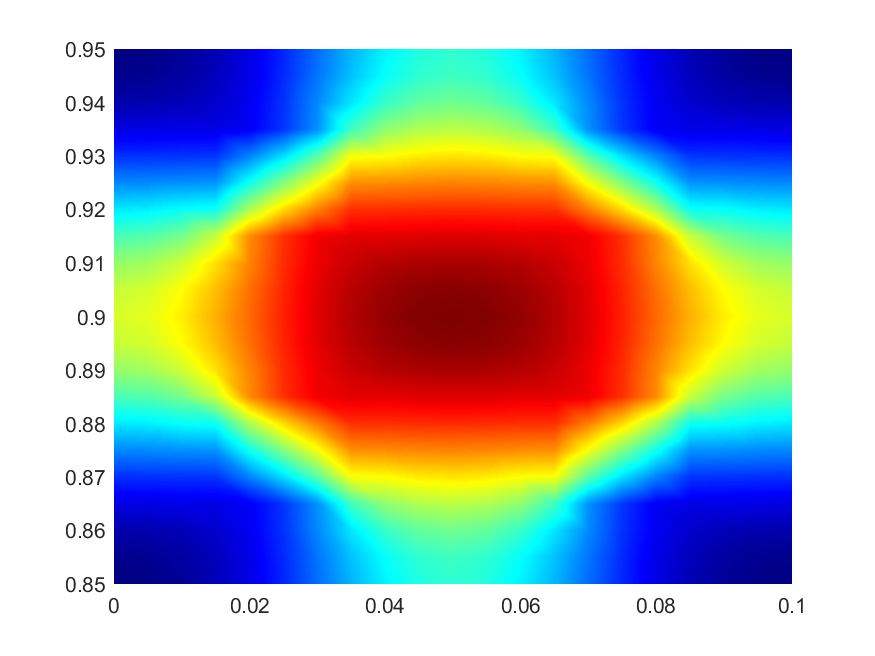}
  \includegraphics[width=0.15\textwidth]{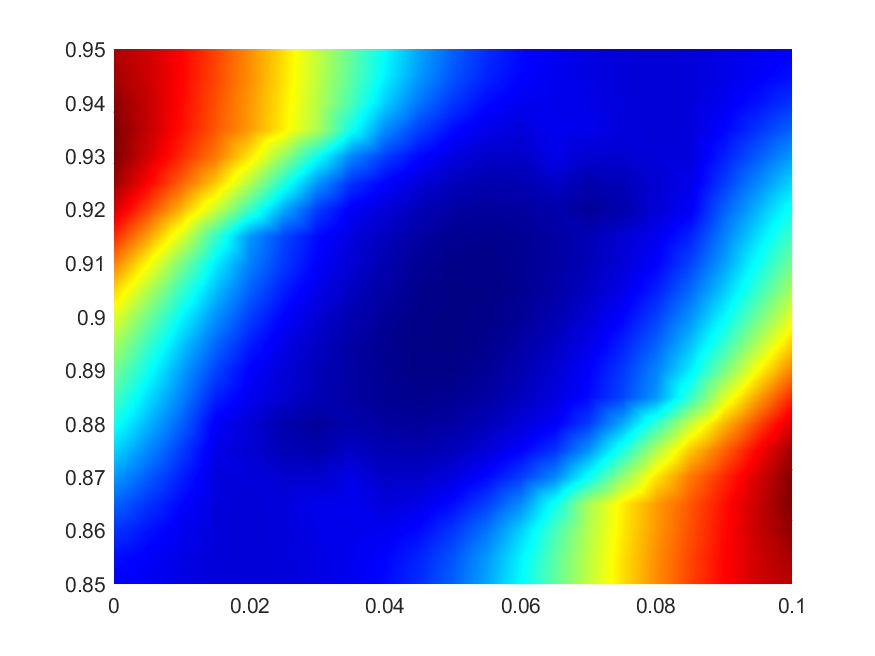}
  \includegraphics[width=0.15\textwidth]{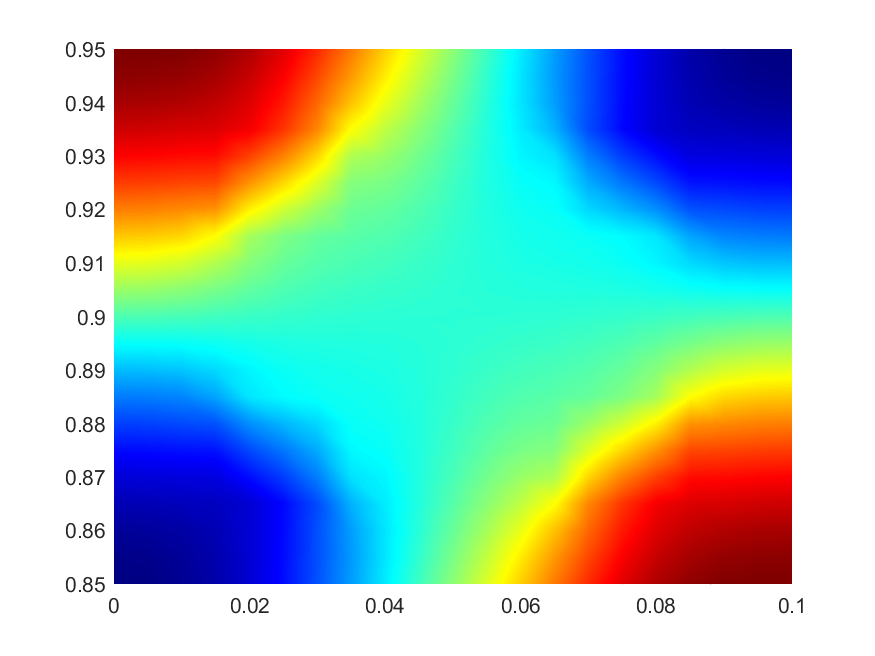}
  {\tiny(7)}\includegraphics[width=0.15\textwidth]{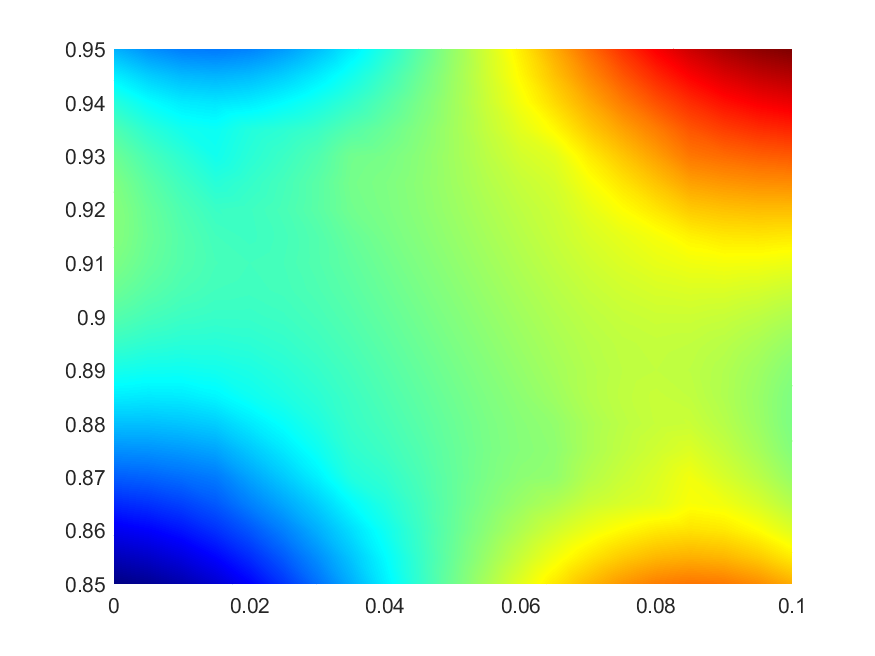}
  \includegraphics[width=0.15\textwidth]{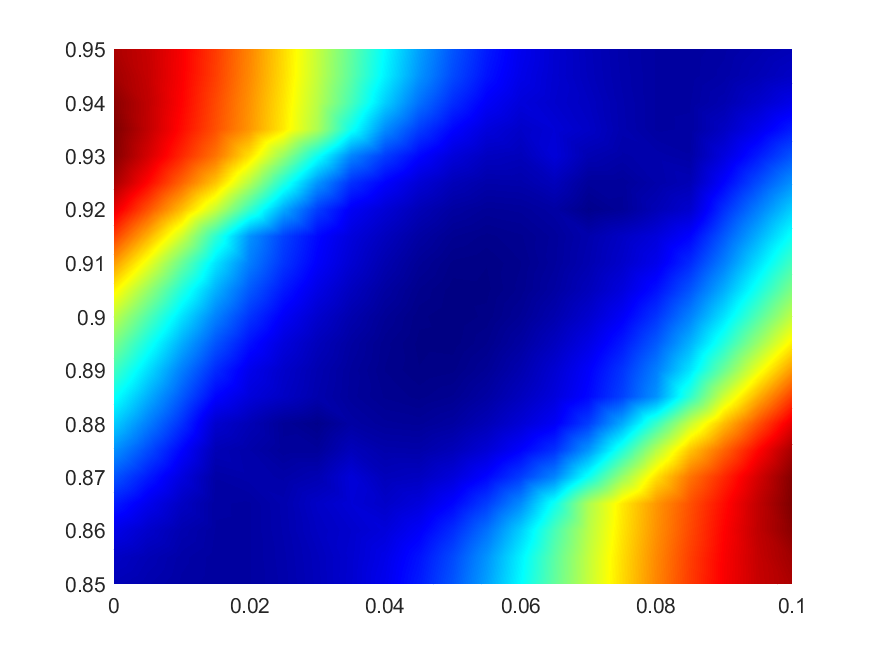}
  \includegraphics[width=0.15\textwidth]{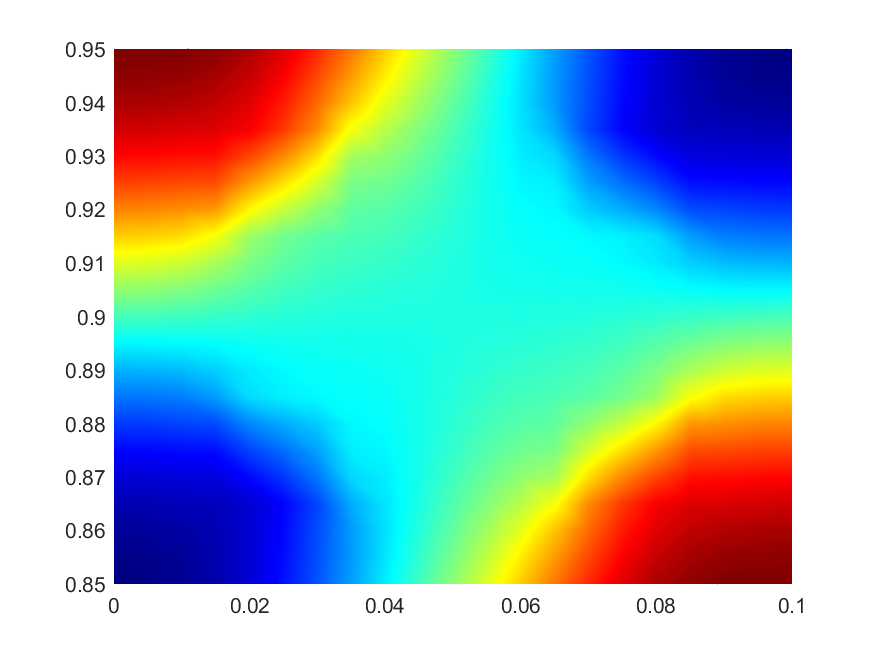}
  \includegraphics[width=0.15\textwidth]{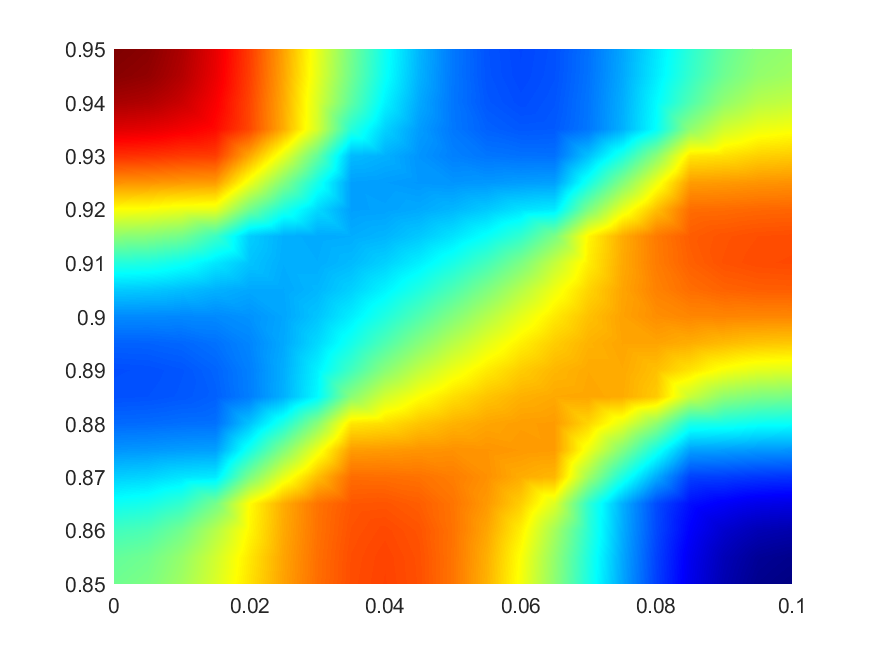}
  \includegraphics[width=0.15\textwidth]{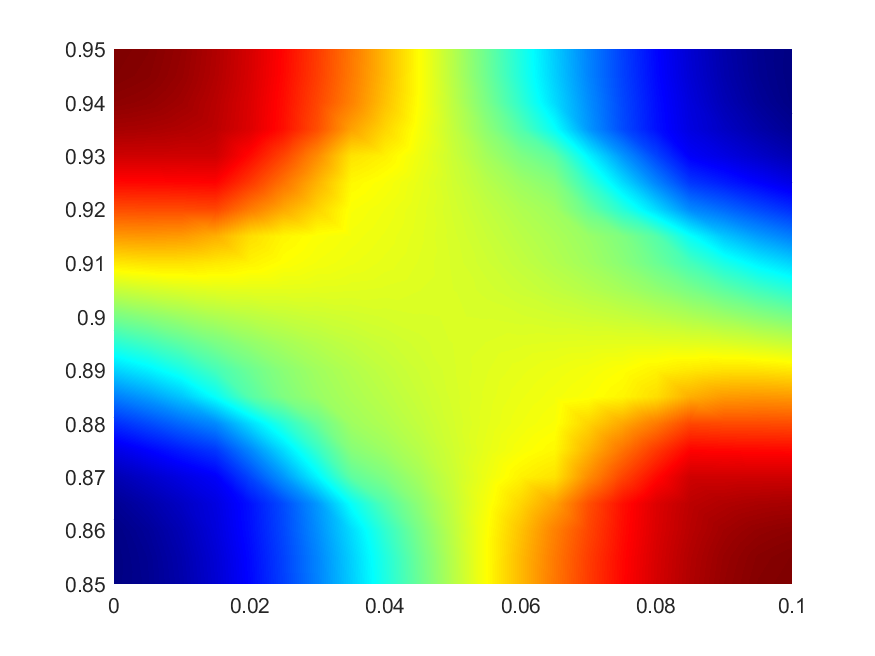}
  \includegraphics[width=0.15\textwidth]{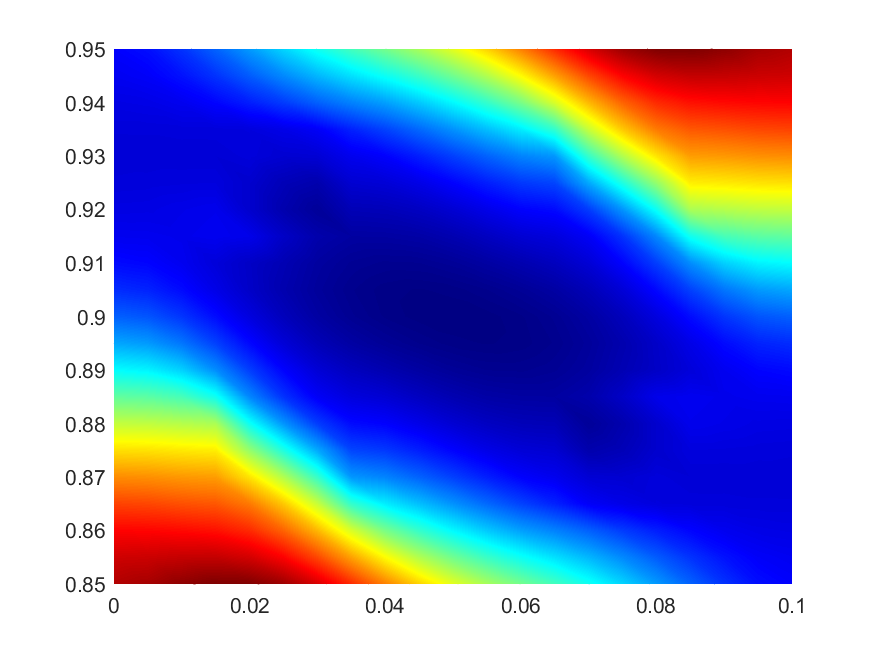}
  {\tiny(8)}\includegraphics[width=0.15\textwidth]{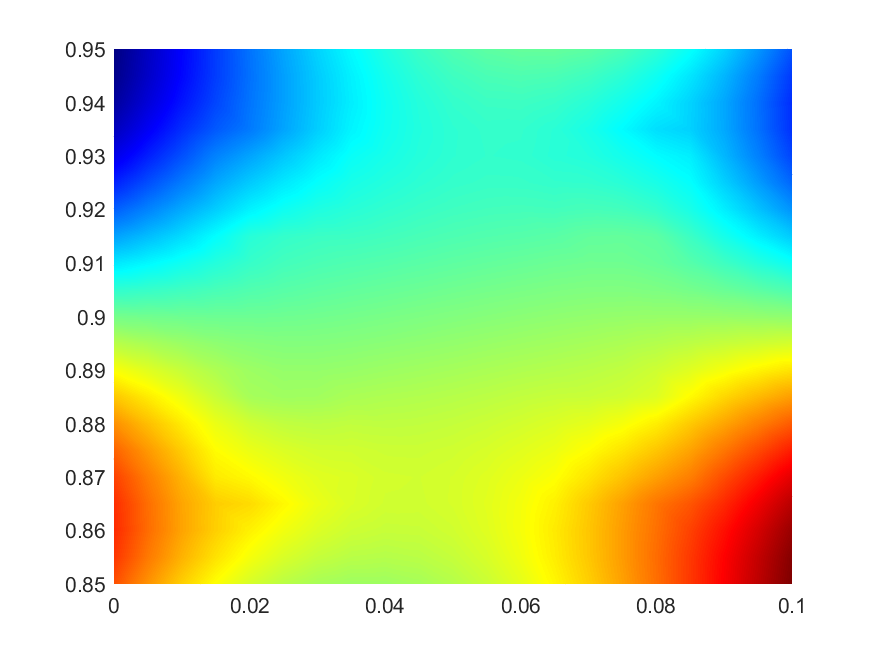}
  \includegraphics[width=0.15\textwidth]{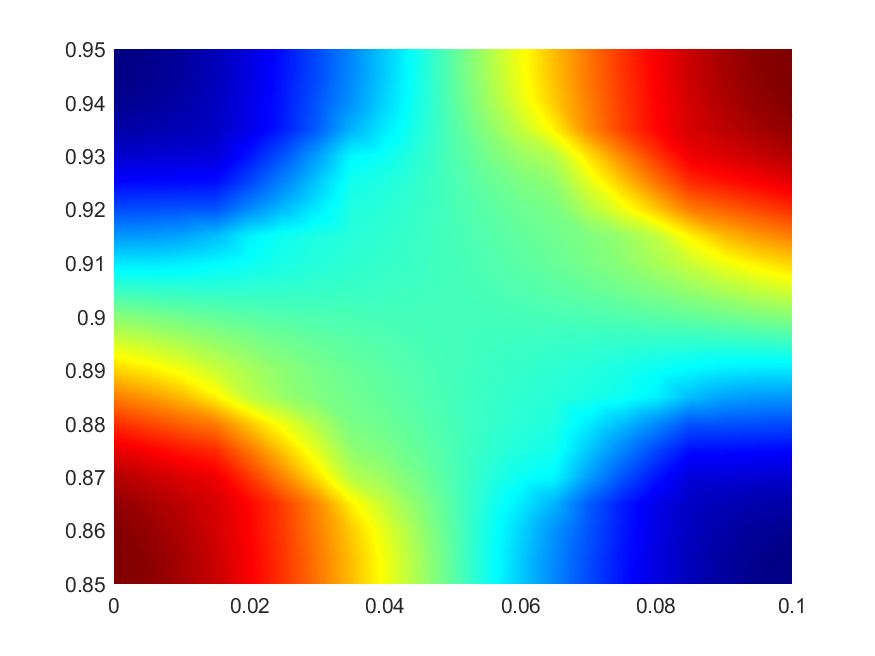}
  \includegraphics[width=0.15\textwidth]{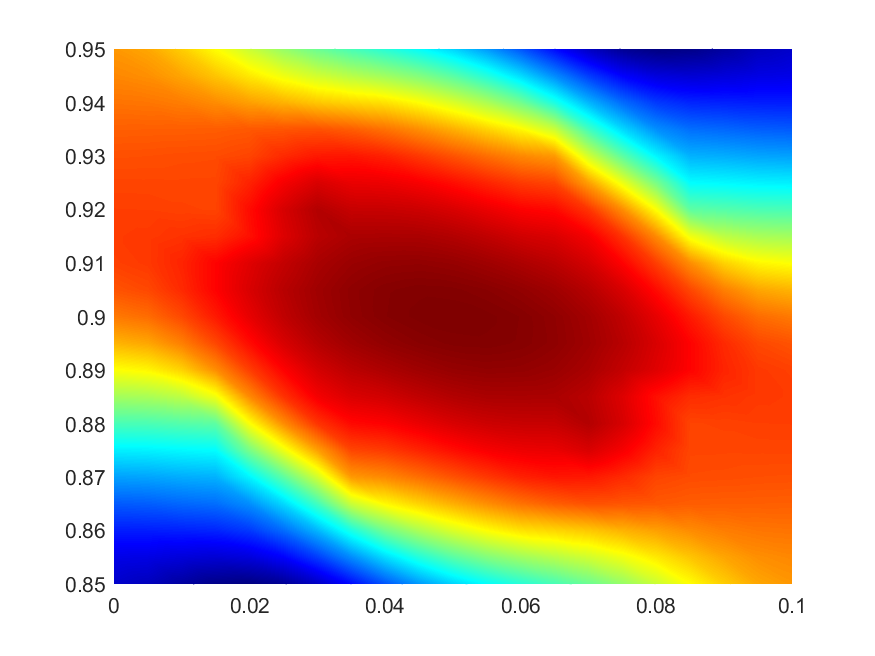}
  \includegraphics[width=0.15\textwidth]{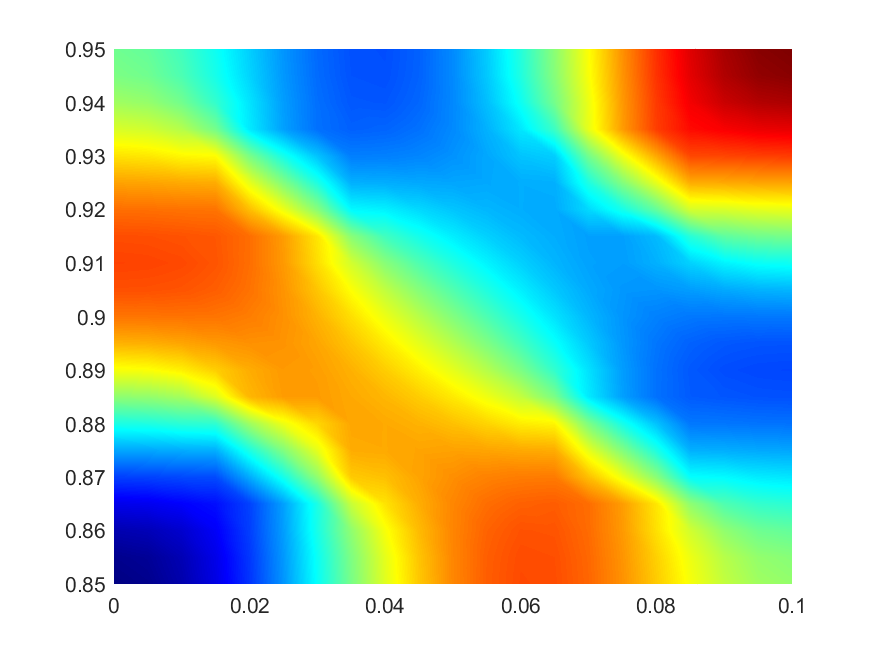}
  \includegraphics[width=0.15\textwidth]{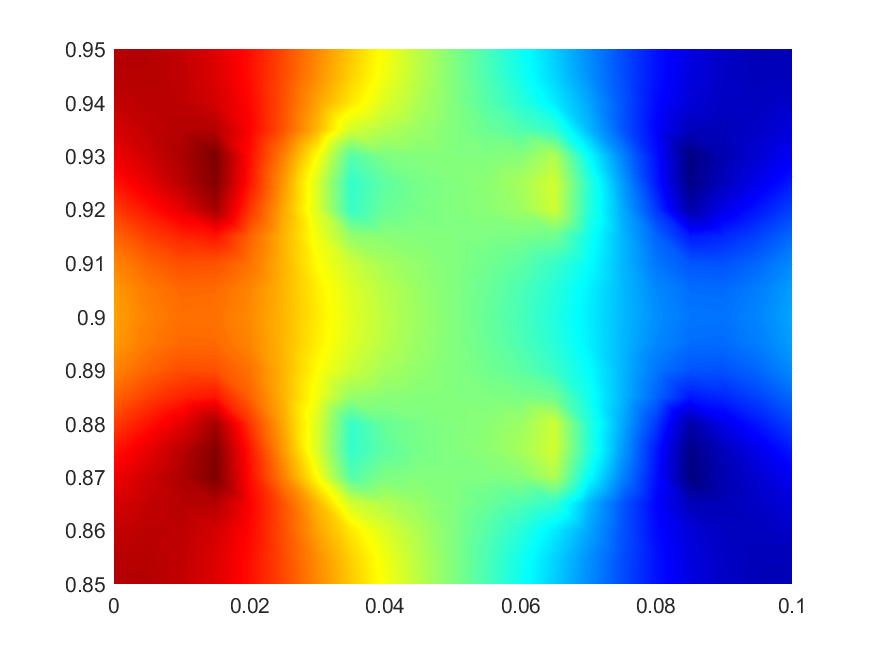}
  \includegraphics[width=0.15\textwidth]{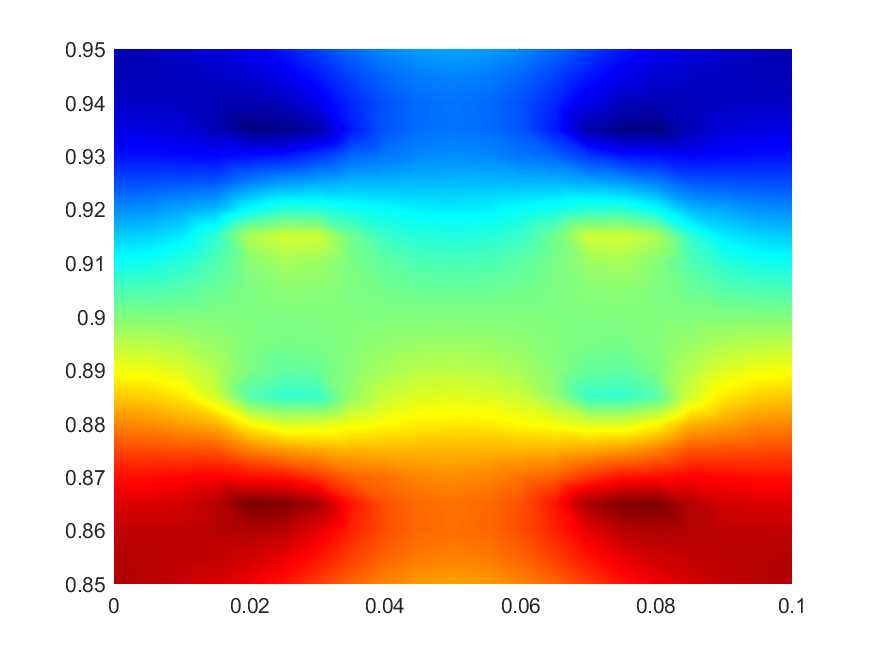}\\
  \hspace{0.6cm}{\tiny(a)$\psi_{\theta}^{cgm}$}\hspace{1cm}
  {\tiny(b)~$\psi_{u_1}^{cgm}$}\hspace{1cm}
  {\tiny(c)~$\psi_{u_2}^{cgm}$}\hspace{1cm}
  {\tiny(d)~$\psi_{\theta}^{gm}$}\hspace{1cm}
  {\tiny(e)~$\psi_{u_1}^{gm}$}\hspace{1.1cm}
  {\tiny(f)~$\psi_{u_2}^{gm}$}
  \caption{
    Contour plots of 8 eigenfunctions of the CGMsFEM and GMsFEM used in Section \ref{subsec:verification}.
    CGMsFEM: (a) $\psi_{\theta}^{cgm}$ (b) $\psi_{u_1}^{cgm}$ and (c) $\psi_{u_2}^{cgm}$;
    GMsFEM: (d) $\psi_{\theta}^{gm}$ (e) $\psi_{u_1}^{gm}$ and (f) $\psi_{u_2}^{gm}$;
  } \label{fig:sec6-fig1}
\end{figure}

Figure \ref{fig:sec6-fig1} demonstrates the difference in eigenfunctions for the CGMsFEM and GMsFEM. It can be seen that the first three or four eigenfunctions from the local coupling spectral problem of the CGMsFEM can better represent the microscopic oscillating information than those from the GMsFEM. This is because for the CGMsFEM, fewer eigenfunctions are needed to represent the physical properties of the microstructure.

\begin{figure}[h]
  \centering
  {\tiny(a)}\includegraphics[width=0.3\textwidth]{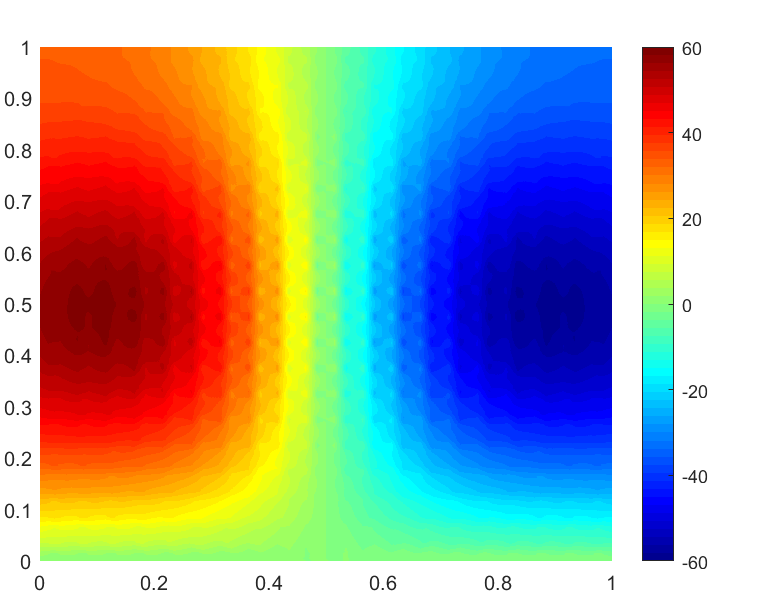}~
  {\tiny(d)}\includegraphics[width=0.3\textwidth]{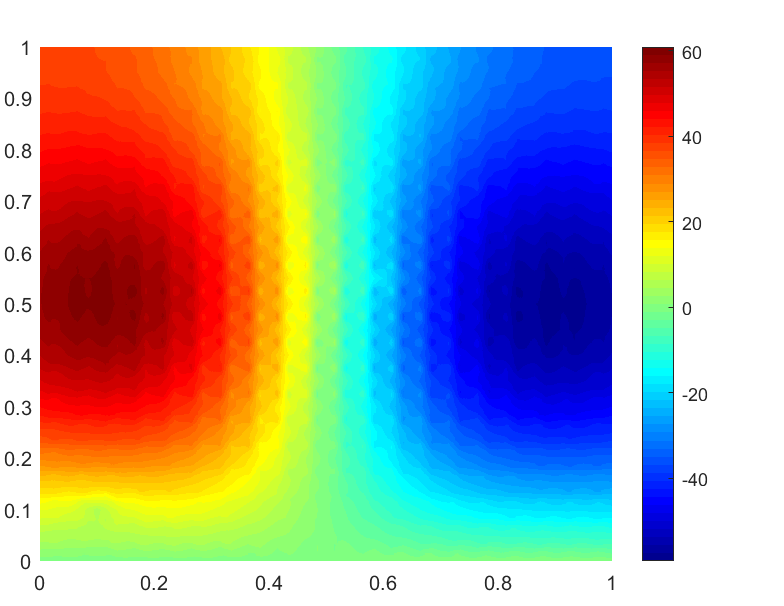}~
  {\tiny(g)}\includegraphics[width=0.3\textwidth]{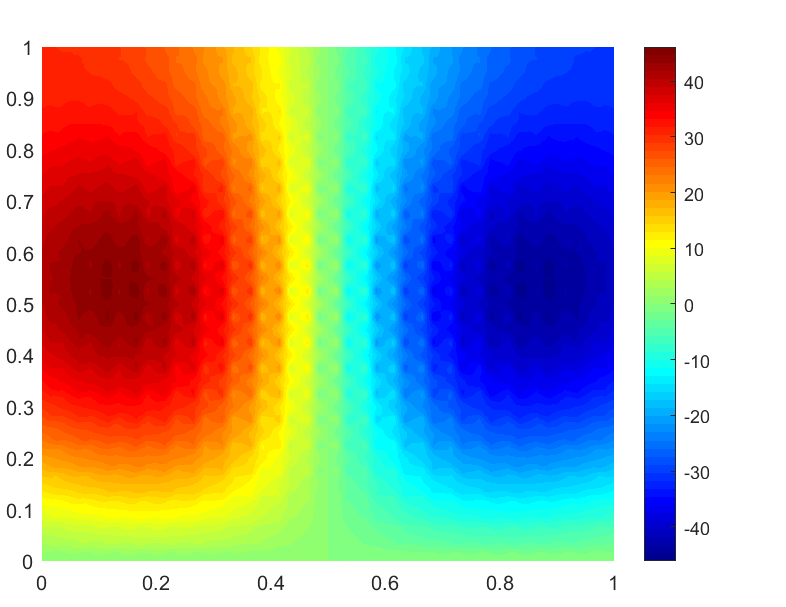}
  {\tiny(b)}\includegraphics[width=0.3\textwidth]{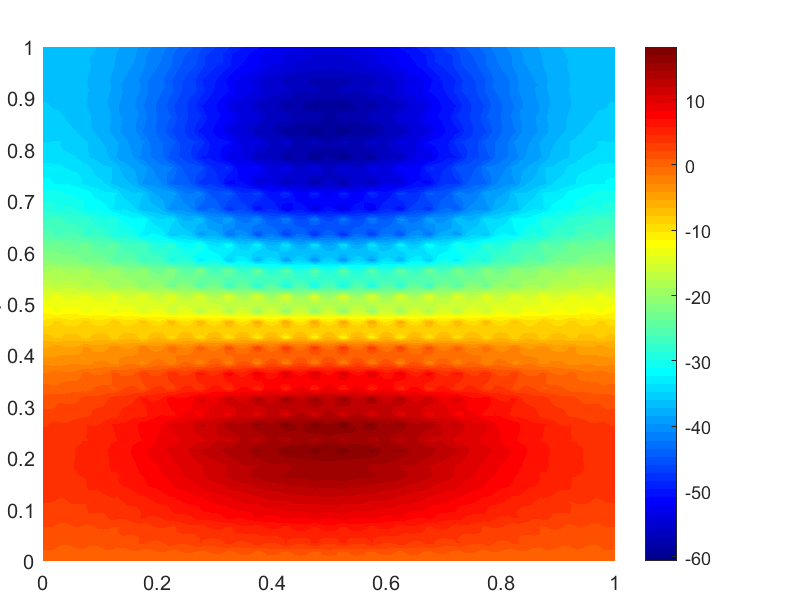}~
  {\tiny(e)}\includegraphics[width=0.3\textwidth]{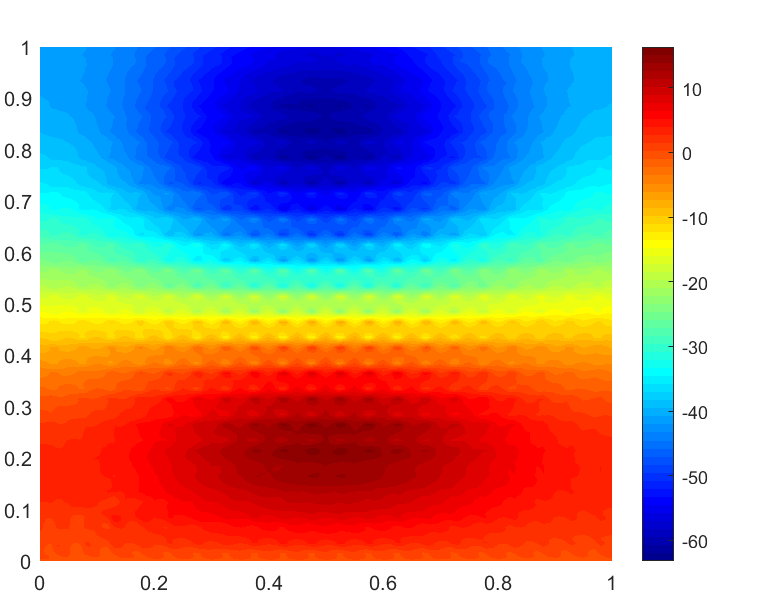}~
  {\tiny(h)}\includegraphics[width=0.3\textwidth]{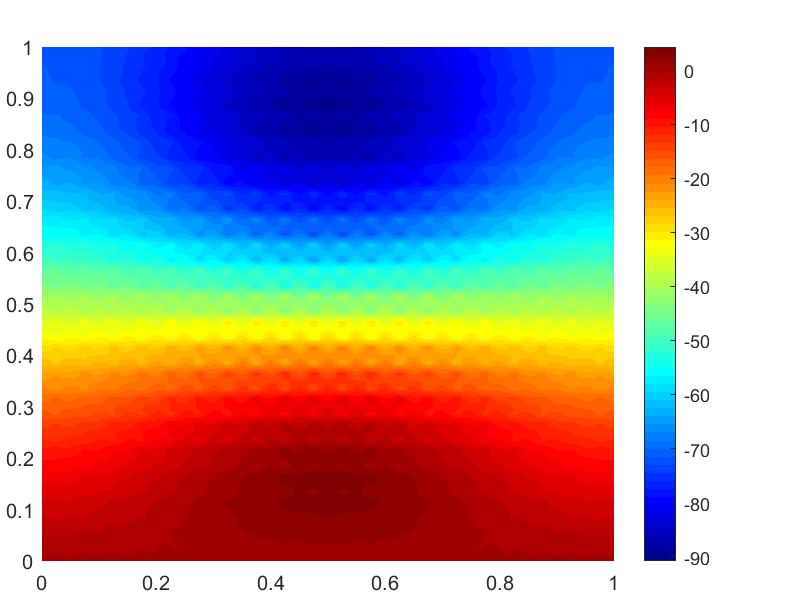}
  {\tiny(c)}\includegraphics[width=0.3\textwidth]{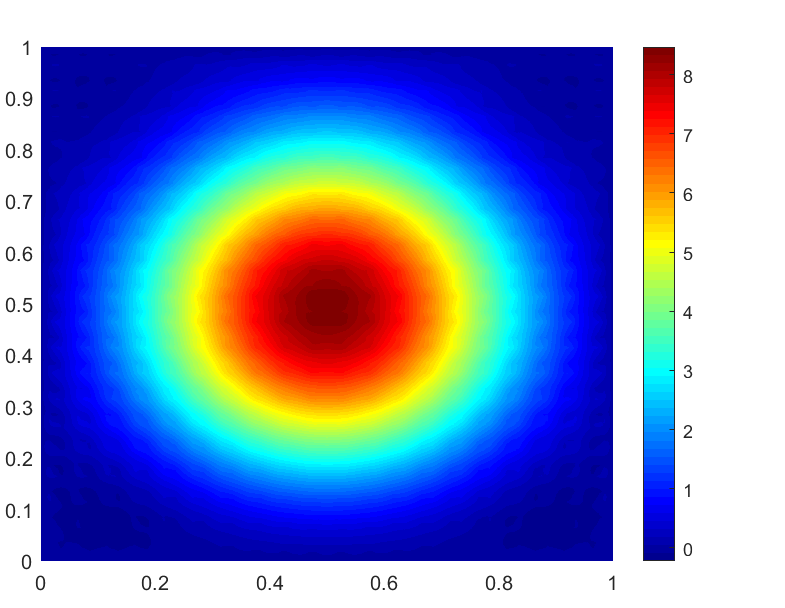}~
  {\tiny(f)}\includegraphics[width=0.3\textwidth]{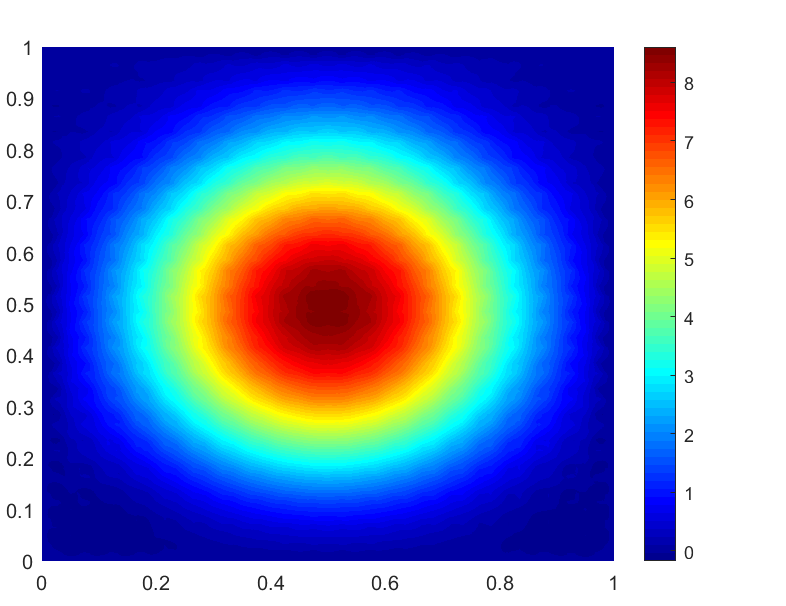}~
  {\tiny(i)}\includegraphics[width=0.3\textwidth]{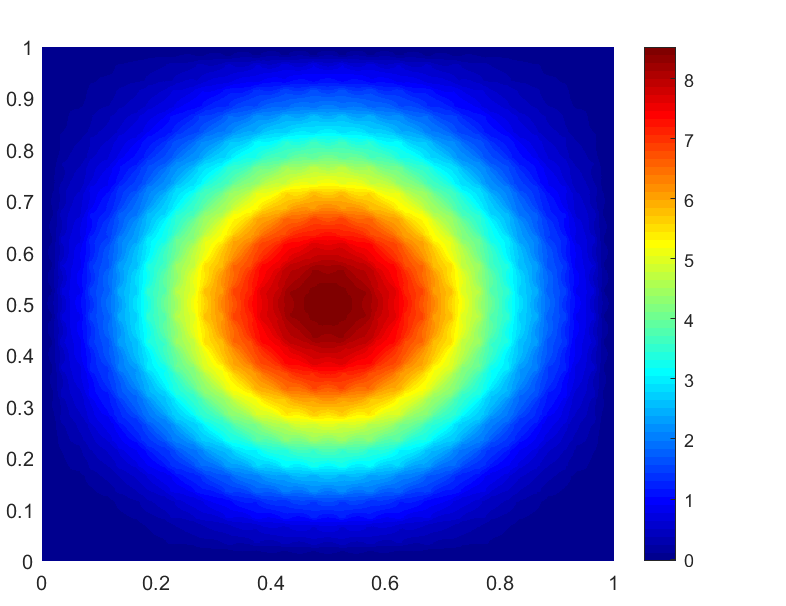}
  \caption{
  Contour plots of solutions for periodic microstructure.
  The reference solutions: (a) $u_1^{ref}$ (b) $u_2^{ref}$ and (c) $\theta^{ref}$;
  The CGMsFEM solutions: (d) $u_1^{cgm}$ (e) $u_2^{cgm}$ and (f) $\theta^{cgm}$;
  The GMsFEM  solutions: (g) $u_1^{gm}$ (h) $u_2^{gm}$ and (i) $\theta^{gm}$.
  }\label{sec5-fig2-solu}
\end{figure}

\begin{table}[htbp]
  \footnotesize
  \caption{Comparison of energy errors for CGMsFEM and CGMsFEM}\label{sec5-tab1}
  \begin{center}
    \begin{tabular}{ccccccccccc}
      \toprule
      $L$ & $Eigvalue$ & $||E_\theta^{cgm}||_e$ & $||E_\theta^{gm}||_e$ & $||E_u^{cgm}||_e$ & $||E_u^{gm}||_e$ & $||E_w^{cgm}||_e$ & $||E_w^{gm}||_e$ \\
      \midrule
      4   & 14.80      & 0.1630                 & 0.4611                & 0.2048            & 0.3226           & 0.1922            & 0.3734           \\
      6   & 35.56      & 0.1316                 & 0.4221                & 0.1687            & 0.2914           & 0.1576            & 0.3395           \\
      8   & 37.54      & 0.1275                 & 0.4212                & 0.1420            & 0.2897           & 0.1374            & 0.3382           \\
      10  & 40.86      & 0.0793                 & 0.4045                & 0.0790            & 0.2325           & 0.0791            & 0.2995           \\
      12  & 43.65      & 0.0687                 & 0.3972                & 0.0728            & 0.2134           & 0.0715            & 0.2865           \\
      14  & 70.87      & 0.0634                 & 0.3880                & 0.0543            & 0.2046           & 0.0574            & 0.2779           \\
      16  & 77.68      & 0.0590                 & 0.3785                & 0.0383            & 0.1768           & 0.0461            & 0.2602           \\
      \bottomrule
    \end{tabular}
  \end{center}
\end{table}
\begin{figure}[h]
  \centering
  {\tiny(a)}\includegraphics[width=0.3\textwidth]{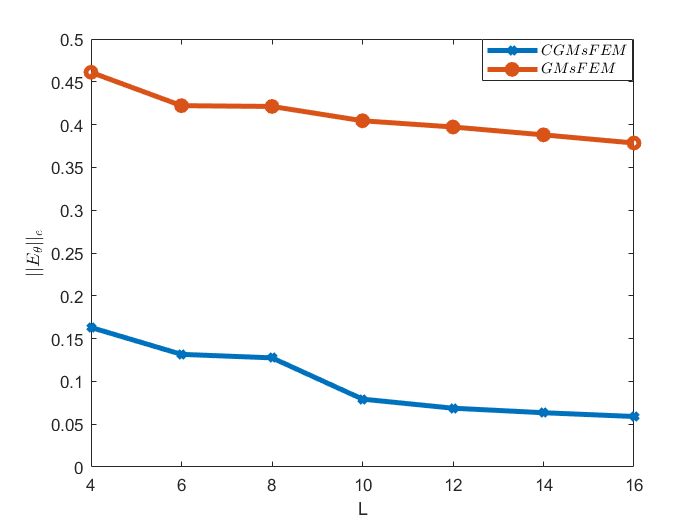}
  {\tiny(b)}\includegraphics[width=0.3\textwidth]{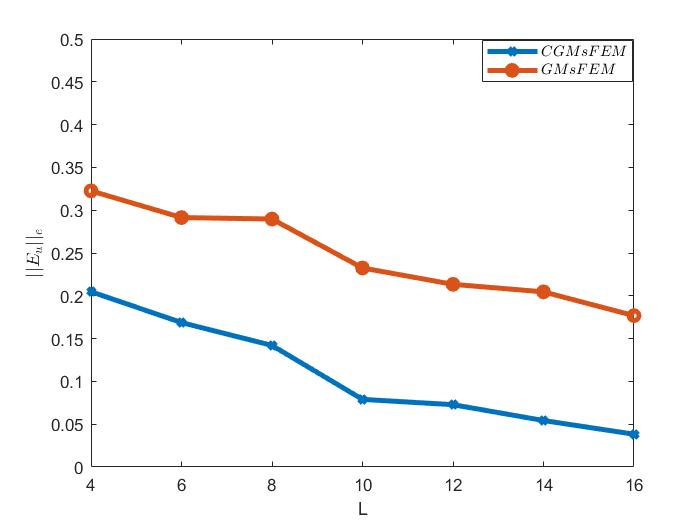}
  {\tiny(c)}\includegraphics[width=0.3\textwidth]{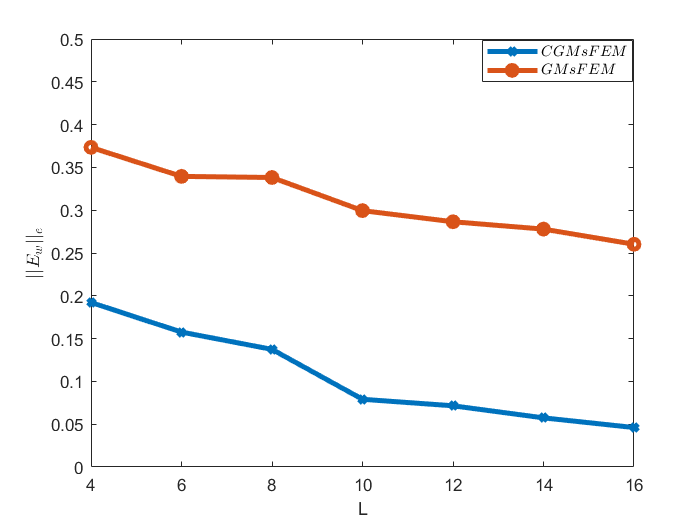}
  \caption{
    Comparison of relative energy errors of the CGMsFEM and GMsFEM in periodic microstructure.
    (a) $E_\theta$ (b) $E_u$, and (c) $E_w$;
  }\label{sec5-fig3-err}
\end{figure}
Figure  \ref{sec5-fig2-solu} demonstrates the contour plots of reference solutions $(\uu^{ref},~\theta^{ref})$, CGMsFEM solutions $(\uu^{cgm},~\theta^{cgm})$, and GMsFEM solutions $(\uu^{gm},~\theta^{gm})$. It can be concluded that the proposed CGMsFEM has higher accuracy than the GMsFEM with the same number of multiscale basis functions. For comparison purposes, we \textcolor{black}{calculate} the energy errors of CGMsFEM and CGMsFEM as defined by the Eq. (\ref{eq:energyerror}). The numercal results for different number $L=4,6,\cdots,16$ of multiscale basis functions are given in Table \ref{sec5-tab1}, \textcolor{black}{where the column labeled 'Eigenvalue' refers to $\Lambda_{L+1}$ of the CGMsFEM}. FFrom the table, we observe that the total relative energy error $||E_w^{cgm}||_e$ for CGMsFEM decreases as $L$ increases, and similar results can also be shown for the energy errors $||E_\theta^{cgm}||_e$ and $||E_u^{cgm}||_e$. Moreover, it can also clearly be found that the energy errors of CGMsFEM with $L=4$ are obviously smaller than those of GMsFEM with $L=16$ in Figure \ref{sec5-fig3-err}, which demonstrate that the CGMsFEM is more efficient than the GMsFEM.

\subsection{Application to random microstructure material coefficients}
In order to validate the good applicability of CGMsFEM, two kinds of tests (Test A, B) are performed in this subsection. In test A, heterogeneous media have random microstructure and deterministic coefficients, and in test B,  both the microstructure and coefficients of materials are random.
\begin{figure}[h]
  \centering
  {\tiny(a)}\includegraphics[width=0.45\textwidth]{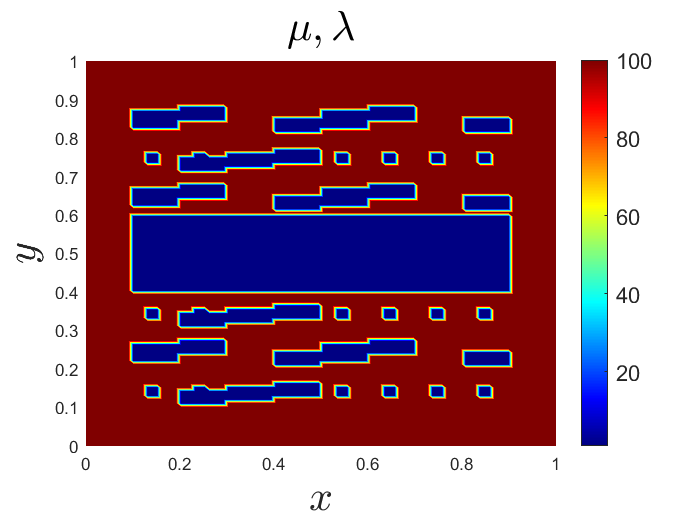}~
  {\tiny(c)}\includegraphics[width=0.45\textwidth]{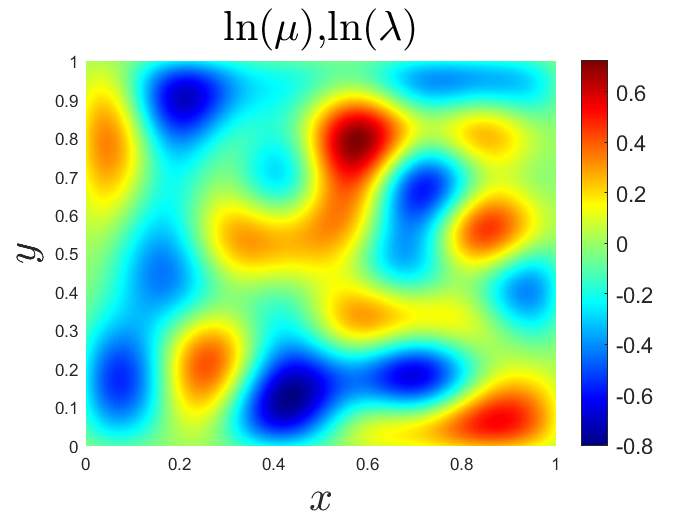}
  {\tiny(b)}\includegraphics[width=0.45\textwidth]{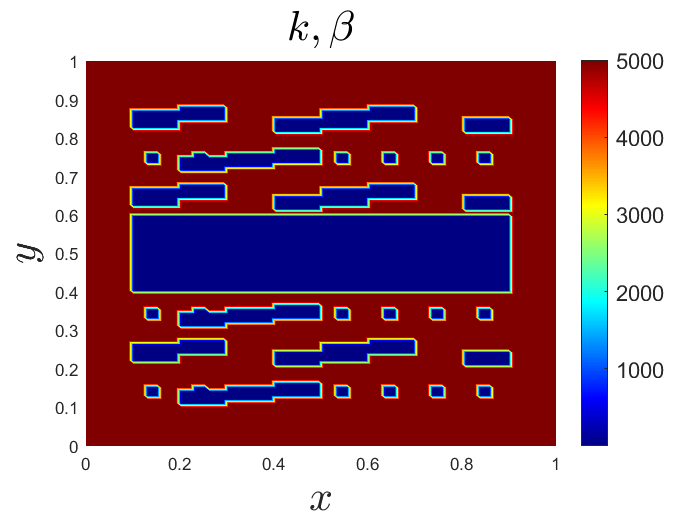}~
  {\tiny(d)}\includegraphics[width=0.45\textwidth]{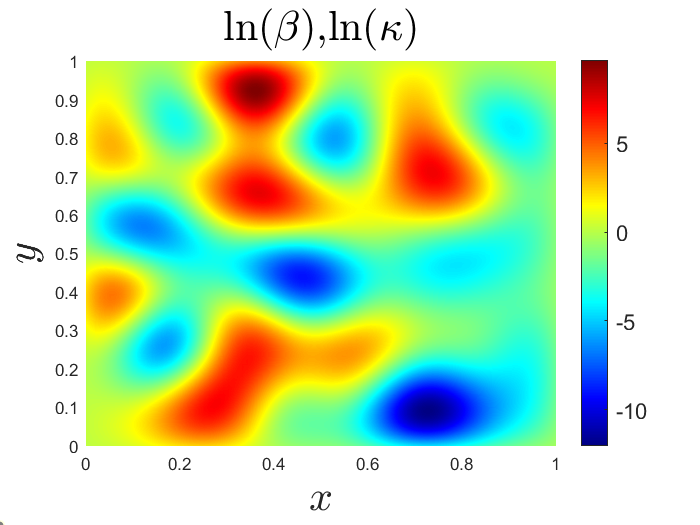}
  \caption{
    Contour plots of one sample for the material coefficients.
    Lam$\acute{e}$ coefficients $\mu$, and $\lambda$: (a) Test A and (c) Test B ;
    Thermal conductivity coefficient $\kappa$, and expansion coefficient $\beta$: (b) Test A and (d) Test B .
  }\label{fig:caseAB}
\end{figure}

\textbf{Test A:}\textit{ Heterogeneous media with random microstructure and deterministic coefficients.} In this simulation, the body force $\ff$ and heat source $g$ are chosen as
$$\ff(x,y)=\mathbf{0},\quad g(x,y)=10\times exp\Big(-\frac{(x-0.2)^2+(y-0.4)^2}{2*0.2^2}\Big).$$
The initial boundary \textcolor{black}{condition} is defined as
$$\theta_0(x,y)=cos(\pi x)cos(\pi y)+1.5.$$
Then the Lam$\acute{e}$ coefficients $\mu$, $\lambda$, Conductivity coefficient $\kappa$, and thermal expansion coefficient $\beta$ are showed in the left of Figure \ref{fig:caseAB}, where the contrasts are chosen as $\lambda _{\text{max}}:\lambda_ {\text{min}}=10^2:1$, $\mu_{\text{max}}:\mu_{\text{min}}=10^2:1$, $\kappa_{\text{max}}:\kappa_{\text{min}}=10^3:1$,$\beta_{\text{max}}:\beta_{\text{min}}=5.0\times10^4:1$. The corresponding relaxation coefficients $\gamma_1=0.75,\gamma_2=7.0\times10^{-2}$ are given, and the time step is $\tau=0.01$. Here, the $100\times100$ fine grid is used for reference solution, and $10\times 10$ coarse grid for the proposed CGMsFEM and GMsFEM. The number of local coupling multiscale basis functions for CGMsFEM are fixed to 10, and the total number of local GMsFEM multiscale basis functions for displacement $\uu$ and $\theta$ are also chosen as 10.

\textbf{Test B:}\textit{ Heterogeneous media with random microstructure and coefficients.}
In this test, the initial conditions, boundary conditions, and source terms are the same as in Section \ref{subsec:verification}. The material coefficients $\kappa(\xx;\xi)$, $\lambda(\xx;\xi) $, $\mu(\xx;\xi)$, and $\beta(\xx;\xi)$ satisfy the following logarithmic Gaussian random field
\begin{equation}  \nonumber
  \text{exp}\Big(\mathcal{GP} \big(b_0(\xx),\text{Cov}(\xx_1,\xx_2)\big)\Big),
\end{equation}
where $\text{Cov}(\xx_1,\xx_2)=\sigma^2\text{exp}(-||\xx_1-\xx_2||^2/l^2)$, $x_1$ and $x_2$ are spatial coordinates in $\Omega$, $\sigma^2$ is the overal variance, and $l$ is the length scale. Then we take $l=0.01,\sigma=10,b_0(x)=0$ for $\kappa(\xx;\xi)$, $\beta(\xx;\xi)$, and $l=0.01,\sigma=10,b_0(x)=0$ for $\lambda(\xx;\xi)$, $\mu(\xx;\xi)$. For simplicity, the samples are obtained with the help of Karhunen-Lo$\acute{e}$ve expansion (KLE) with 50 truncated terms. The right parts of Figure \ref{fig:caseAB} show one sample of $\kappa(\xx;\xi)$, $\lambda(\xx;\xi)$, $\mu(\xx;\xi)$ and $\beta(\xx;\xi)$. Let the time step size be $\tau_n=0.01$. Fine grid size $h=\frac{1}{100}$, coarse grid size $H=\frac{1}{10}$, and relaxation coefficients $\gamma_1=0.7,\gamma_2=8.0\times10^{-6}$ are used for simulation. The number of local coupling multiscale basis functions for CGMsFEM is fixed at 8, and the total number of local GMsFEM multiscal basis functions for displacement $\uu$ and $\theta$ is also chosen at 8.

\begin{figure}[h]
  \centering
  {\tiny(a)}\includegraphics[width=0.3\textwidth]{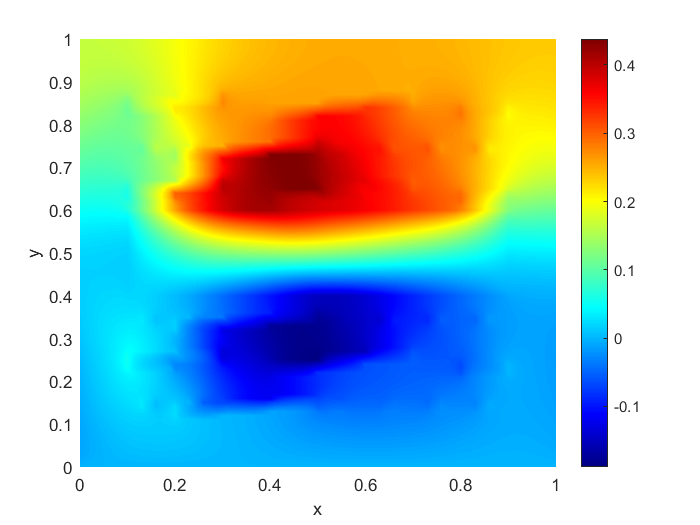}~
  {\tiny(d)}\includegraphics[width=0.3\textwidth]{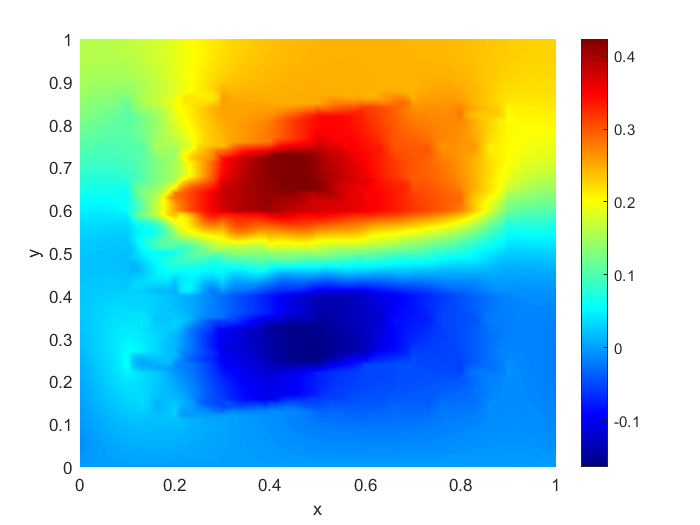}~
  {\tiny(g)}\includegraphics[width=0.3\textwidth]{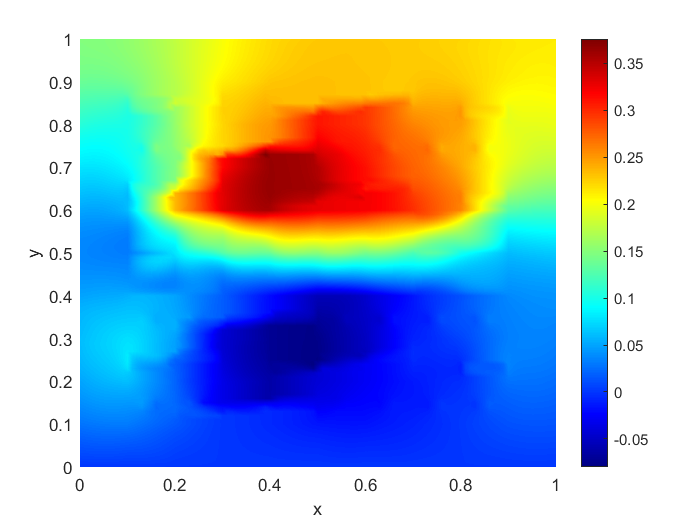}
  {\tiny(b)}\includegraphics[width=0.3\textwidth]{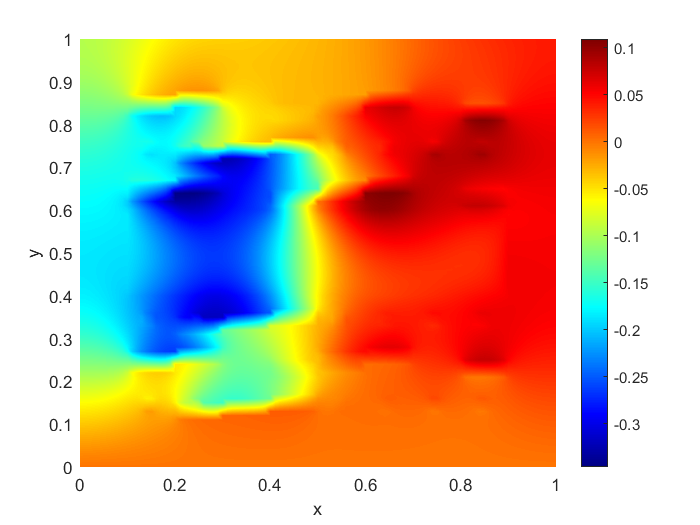}~
  {\tiny(e)}\includegraphics[width=0.3\textwidth]{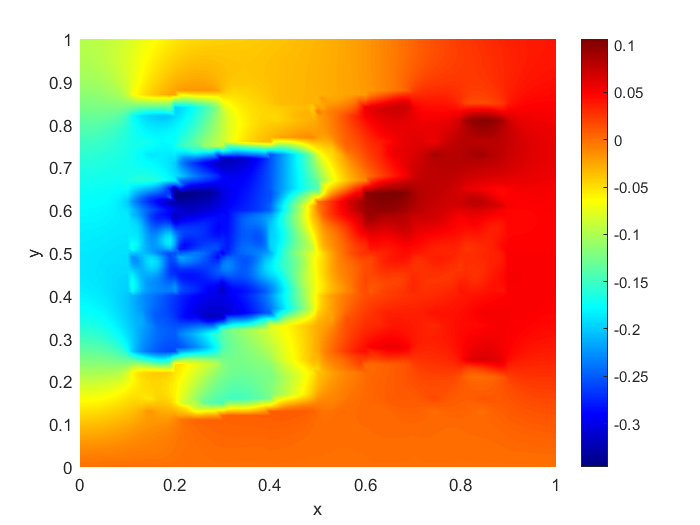}~
  {\tiny(h)}\includegraphics[width=0.3\textwidth]{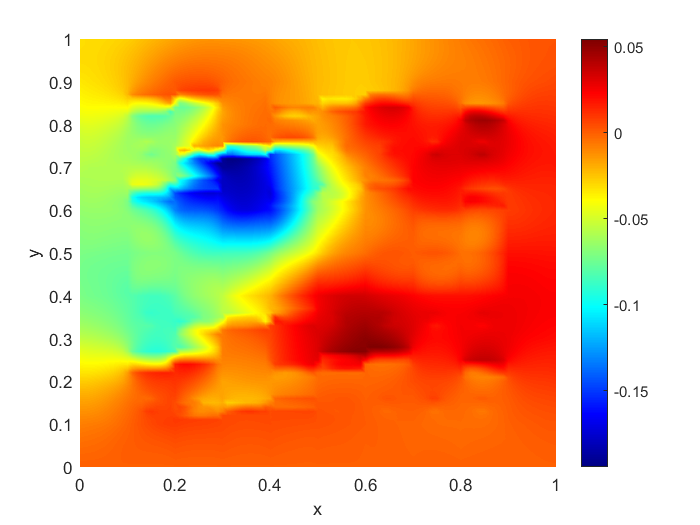}
  {\tiny(c)}\includegraphics[width=0.3\textwidth]{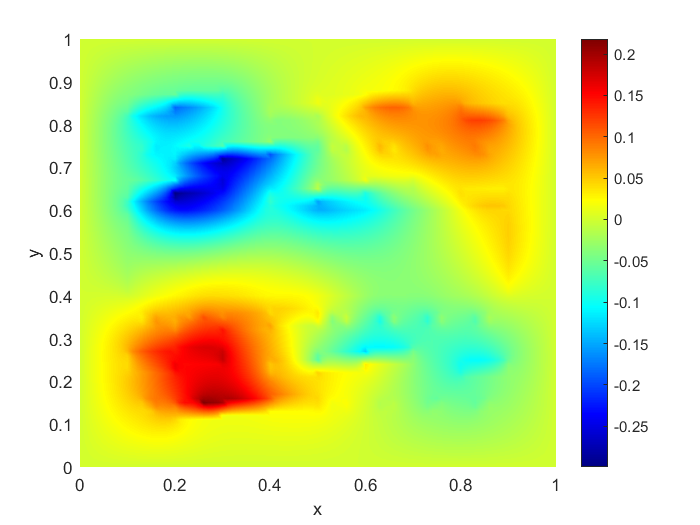}~
  {\tiny(f)}\includegraphics[width=0.3\textwidth]{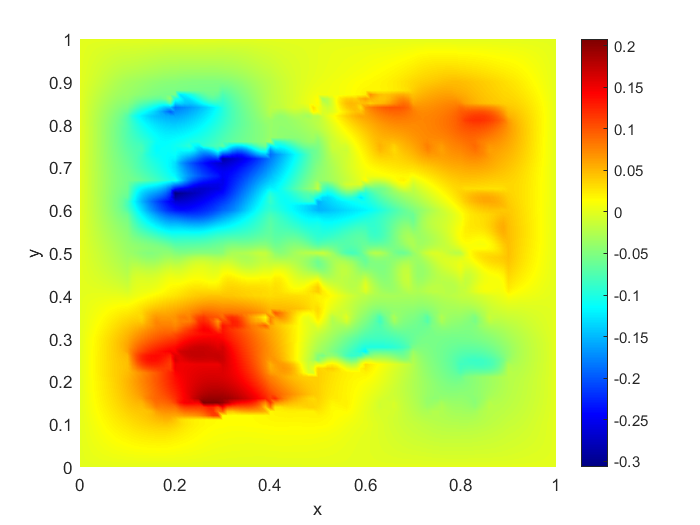}~
  {\tiny(i)}\includegraphics[width=0.3\textwidth]{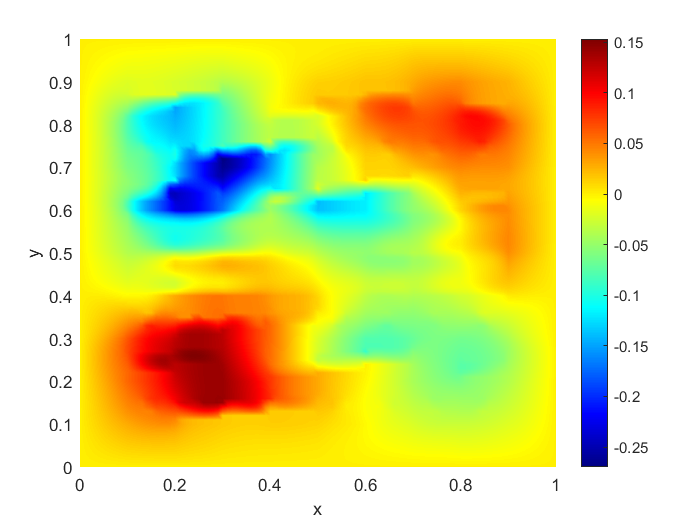}
  \caption{
  Contour plots of solutions for Test A.
  The reference solutions: (a) $u_1^{ref}$ (b) $u_2^{ref}$ and (c) $\theta^{ref}$;
  The CGMsFEM solutions: (d) $u_1^{cgm}$ (e) $u_2^{cgm}$ and (f) $\theta^{cgm}$;
  The GMsFEM  solutions: (g) $u_1^{gm}$ (h) $u_2^{gm}$ and (i) $\theta^{gm}$.
  }\label{fig:caseAsolution}
\end{figure}

\begin{figure}[h]
  \centering
  {\tiny(a)}\includegraphics[width=0.3\textwidth]{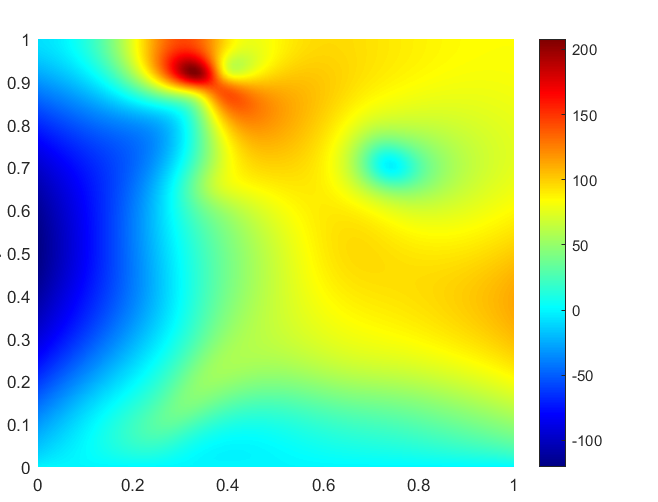}~
  {\tiny(d)}\includegraphics[width=0.3\textwidth]{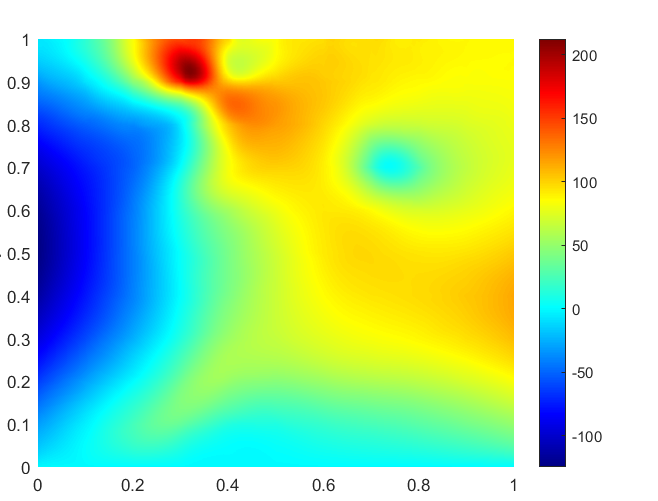}~
  {\tiny(g)}\includegraphics[width=0.3\textwidth]{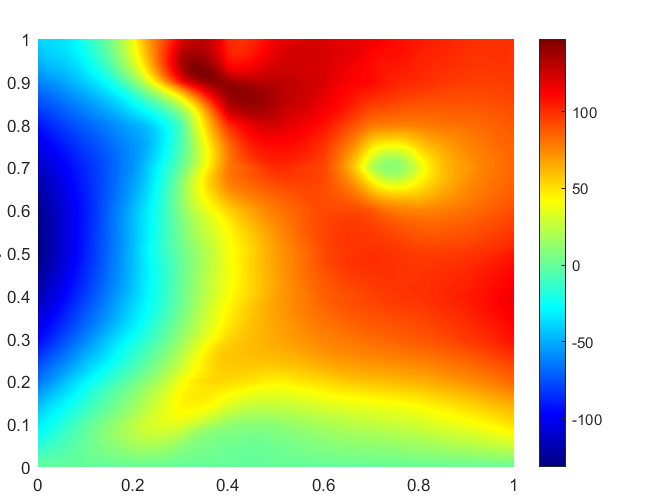}
  {\tiny(b)}\includegraphics[width=0.3\textwidth]{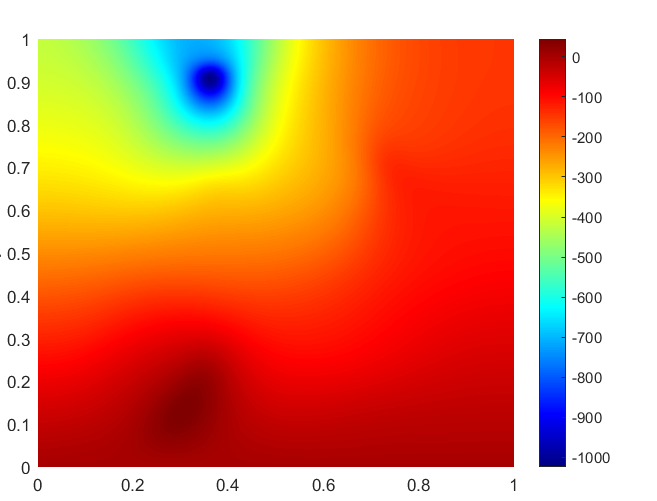}~
  {\tiny(e)}\includegraphics[width=0.3\textwidth]{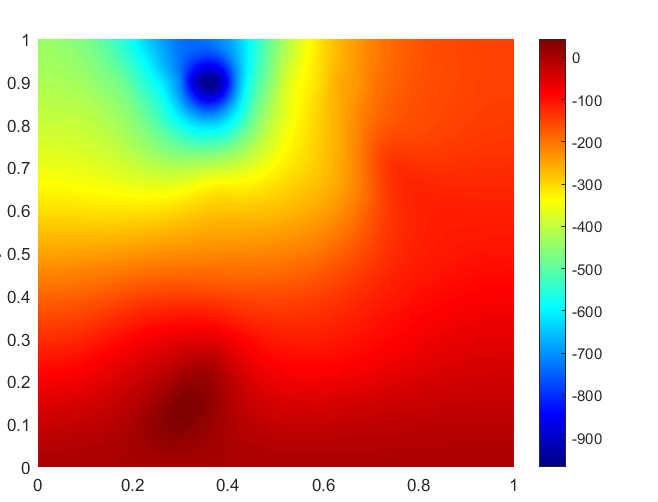}~
  {\tiny(h)}\includegraphics[width=0.3\textwidth]{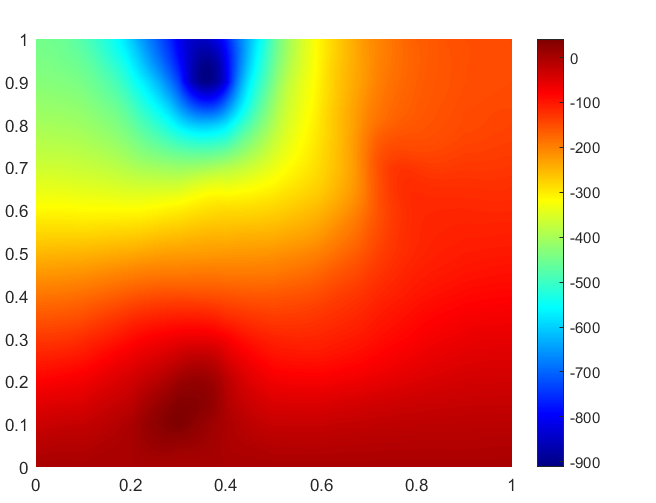}
  {\tiny(c)}\includegraphics[width=0.3\textwidth]{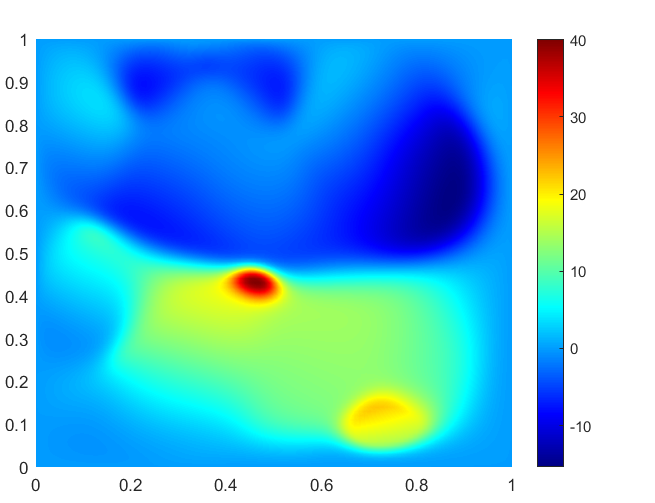}~
  {\tiny(f)}\includegraphics[width=0.3\textwidth]{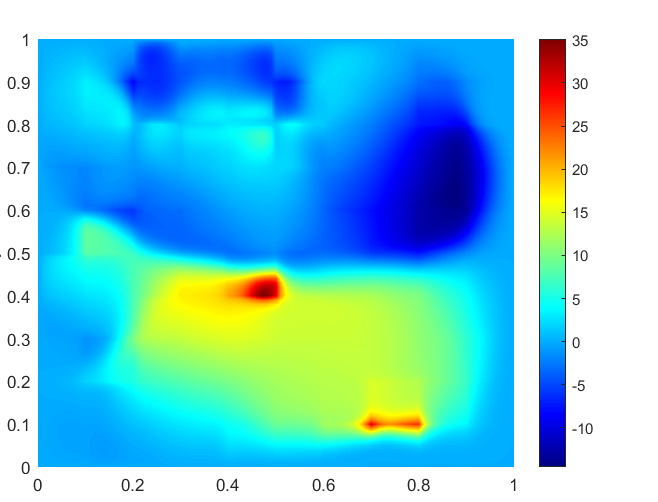}~
  {\tiny(i)}\includegraphics[width=0.3\textwidth]{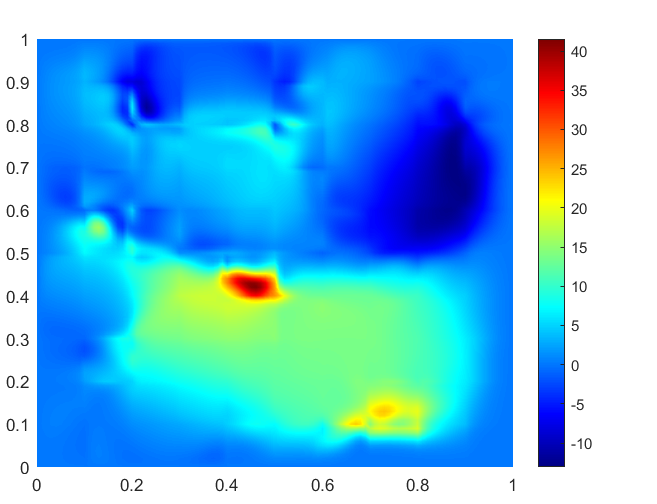}
  \caption{
  Contour plots of solutions for Test B.
  The reference solutions: (a) $u_1^{ref}$ (b) $u_2^{ref}$ and (c) $\theta^{ref}$;
  The CGMsFEM solutions: (d) $u_1^{cgm}$ (e) $u_2^{cgm}$ and (f) $\theta^{cgm}$;
  The GMsFEM  solutions: (g) $u_1^{gm}$ (h) $u_2^{gm}$ and (i) $\theta^{gm}$.
  }\label{fig:caseBsolution}
\end{figure}

\begin{table}[htbp]
	\footnotesize
	\caption{Test A: Relative energy errors of the CGMsFEM and GMsFEM with different contrast ratio of $\beta_{\text{max}}$ and $\beta_{\text{min}}$.}
	\label{sec4-tab-num2-1}
\begin{center}
	\begin{tabular}{ccccccc}
		\toprule
		$Ratio$ & $||E_\theta^{cgm}||_e$ & $||E_\theta^{gm}||_e$ & $||E_u^{cgm}||_e$ & $||E_u^{gm}||_e$ & $||E_w^{cgm}||_e$ & $||E_w^{gm}||_e$ \\
		\midrule
		$ 1\times 10^1 $       & 0.0474                 & 0.0946                & 0.0838            & 0.2682           & 0.0474            & 0.0947           \\
		$ 1\times 10^2 $        & 0.0638                 & 0.0946                & 0.1197            & 0.2538           & 0.0638            & 0.0947           \\
		$ 1\times 10^3 $         & 0.0805                 & 0.1647                & 0.1364            & 0.5539           & 0.0822            & 0.1828           \\
		$ 5\times 10^3 $         & 0.1306                 & 0.8367                & 0.1311            & 0.7455           & 0.1308            & 0.7926           \\
		$ 1\times 10^4 $         & 0.1512                 & 0.5599               & 0.1519            & 0.5171           & 0.1512            & 0.5360           \\
		\bottomrule
	\end{tabular}
\end{center}
\end{table}
\begin{table}[htbp]
  \footnotesize
  \caption{Test B: Relative energy errors of the CGMsFEM and GMsFEM with different variance $\sigma$ of $\beta$.}
  \label{sec4-tab-num2-2}
\begin{center}
	\begin{tabular}{ccccccc}
		\toprule
		$\sigma$ & $||E_\theta^{cgm}||_e$ & $||E_\theta^{gm}||_e$ & $||E_u^{cgm}||_e$ & $||E_u^{gm}||_e$ & $||E_w^{cgm}||_e$ & $||E_w^{gm}||_e$ \\
		\midrule
		2       & 0.2075                 & 0.3259                & 0.2124            & 0.2916           & 0.0474            & 0.0947           \\
		3       & 0.2336                 & 0.3410                & 0.2478            & 0.2871           & 0.2072            & 0.3116           \\
		4       & 0.2739                 & 0.4202                & 0.2119            & 0.2665           & 0.2371            & 0.3152           \\
		5       & 0.2656                 & 0.4194                & 0.2135            & 0.2677           & 0.2225            & 0.3048           \\
		6       & 0.4144                 & 0.5782                & 0.1492            & 0.2250           & 0.1781            & 0.2728           \\
		\bottomrule
	\end{tabular}
\end{center}
\end{table}

\begin{figure}[H]
  \centering
  {\tiny(a)}\includegraphics[width=0.3\textwidth]{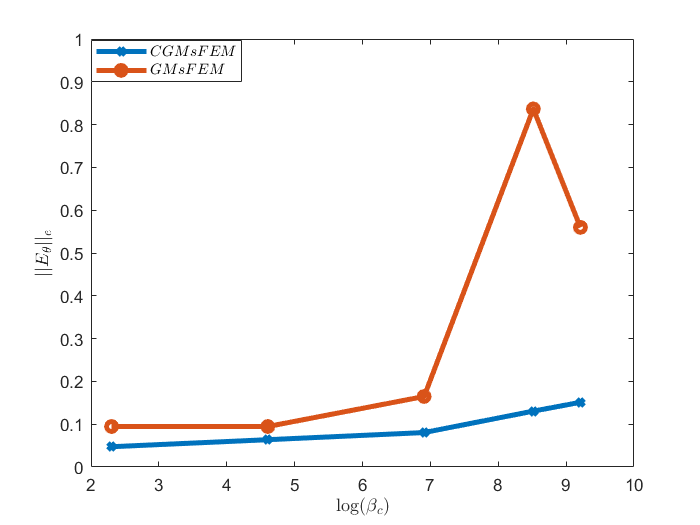}
  {\tiny(b)}\includegraphics[width=0.3\textwidth]{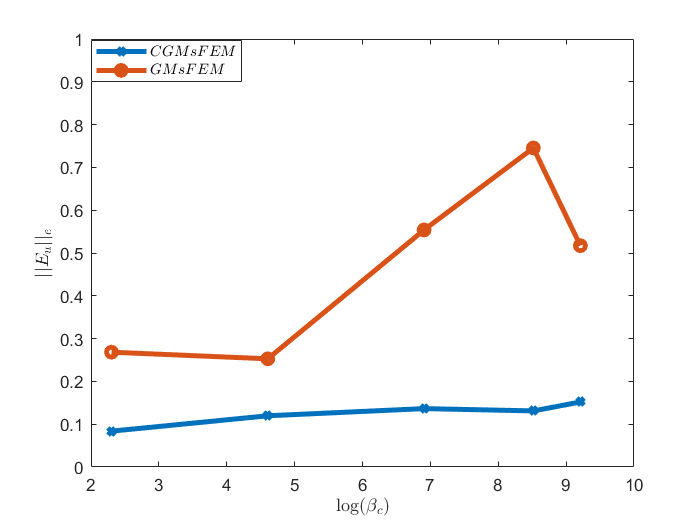}
  {\tiny(c)}\includegraphics[width=0.3\textwidth]{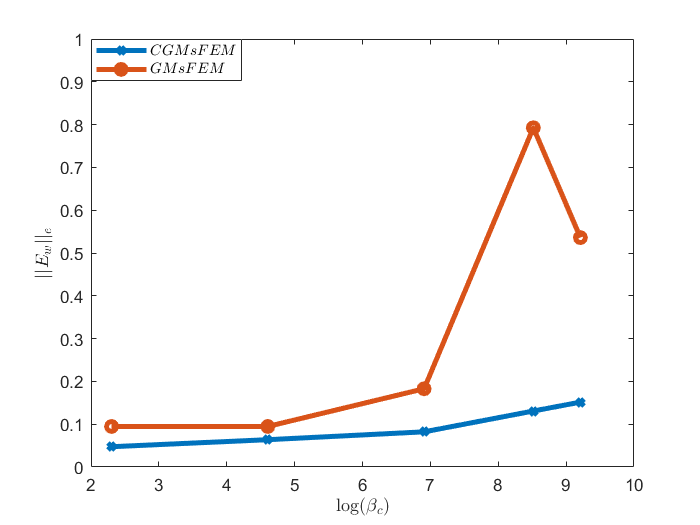}
  \caption{
    Test A: Comparison of relative energy errors for the CGMsFEM and GMsFEM with different contrast ratio of $\beta_{\text{max}}$ and $\beta_{\text{min}}$
    (a) $E_\theta$ (b) $E_u$, and (c) $E_w$, where $\beta_c = \frac{\beta_{\text{max}}}{\beta_{\text{min}}}$.
  }\label{fig:caseAerror}
\end{figure}

\begin{figure}[h]
  \centering
  {\tiny(a)}\includegraphics[width=0.3\textwidth]{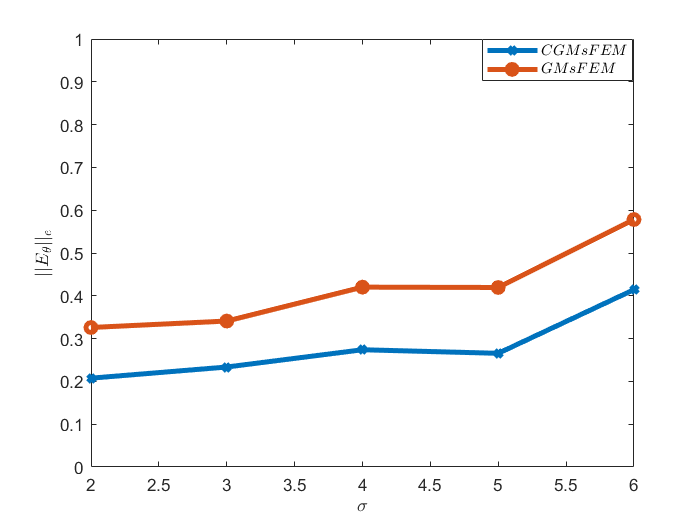}
  {\tiny(b)}\includegraphics[width=0.3\textwidth]{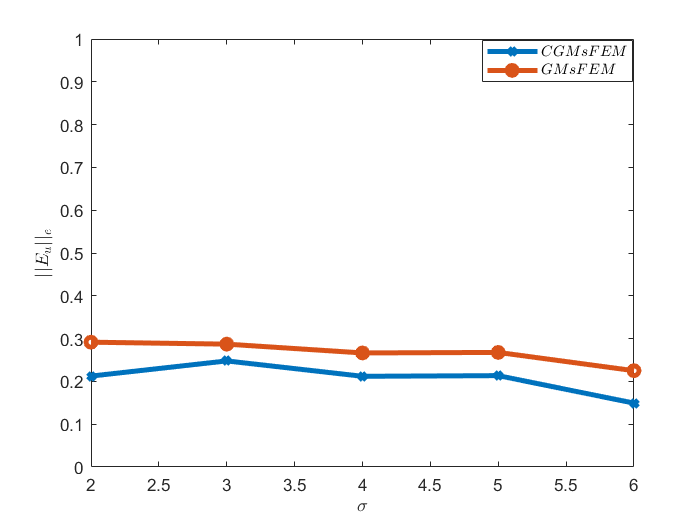}
  {\tiny(c)}\includegraphics[width=0.3\textwidth]{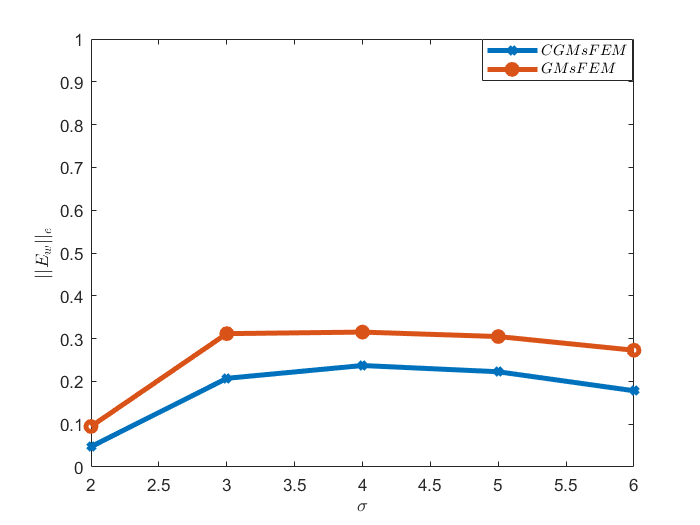}
  \caption{
    Test B: Comparison of relative energy errors of the CGMsFEM and GMsFEM with different variance $\sigma$ of $\beta$.
    (a) $E_\theta$ (b) $E_u$, and (c) $E_w$;
  }\label{fig:caseBerror}
\end{figure}

The reference solutions $(\uu^{ref},\theta^{ref})$, the CGMsFEM solutions $(\uu^{cgm},\theta^{cgm})$, and the GMsFEM solutions $(\uu^{gm},\theta^{gm})$ of Test A and Test B are presented in Figure \ref{fig:caseAsolution} and Figure \ref{fig:caseBsolution}. It can be observed that the CGMsFEM solutions have higher accuracy than the GMsFEM solutions by comparing the reference solutions, which is consistent with the periodic case. Moreover, in order to explore the influence of different coefficients on the results, we compare the energy errors of the CGMsFEM and GMsFEM with the change of $\beta$. In Test A, the contrast ratio of $\beta_{\text{max}}$ and $\beta_{\text{min}}$ constantly varies, and in Test B, the variance $\sigma$ is also constantly changed. Table \ref{sec4-tab-num2-1} and Table \ref{sec4-tab-num2-2} report the energy errors $E_\theta$, $E_u$, and (c) $E_w$ of the CGMsFEM and GMsFEM for Test A and Test B in detail. To make our point more clear, we also show the trend of the energy errors with $\beta_{\text{max}}$ and $\beta_{\text{min}}$ as contrast ratios in Figure \ref{fig:caseAerror} for Test A and with $\sigma$ as the variance of $\text{ln}(\beta)$ in Figure \ref{fig:caseBerror} for Test B. For both the CGMsFEM and GMsFEM, we observe that the energy errors of the CGMsFEM are almost at the same level, while those of the GMsFEM change significantly. Although the whole system becomes extremely complex with increasing the contrast ratio or the variance of two material coefficients, the energy error of the whole system stays in a stable state, which is hardly affected by the complexity of the system. This is also consistent with that in Section \ref{subsec:verification}. In addition, the accuracy of the CGMsFEM is always better than that of the GMsFEM for Tests A and B, even in the more complex cases. Those numerical results demonstrate that the CGMsFEM is computationally quite efficient and accurate, with wide applicability in many scenarios.

\section{Conclusions}
\label{sec:conclusions}
In this paper, a novel CGMsFEM is proposed to efficiently simulate the thermomechanical behaviors of heterogeneous media, which is suitable for both weak and strong coupling settings. To the best of our knowledge, this paper is the first example of designed and analyzed coupling multiscale basis functions for the GMsFEM of multiphysics problems. Two relaxation coefficients are innovatively designed for local regularized coupling spectral problems. The constructed multiscale basis functions can more accurately capture the multiphysical coupling information of the original thermomechanical problems, which obviously reduces the order of the global stiffness matrix. Moreover, the corresponding convergence analysis is derived in detail, where the error of CGMsFEM is closely related to the eigenvalue decay in each local coarse block and the upper error bound is independent of the two relaxation coefficients. Numerical experiments show that the proposed CGMsFEM is stable and effective and can provide enough numerical accuracy, which supports the theoretical results of this paper.

Despite a number of simplifying hypotheses adopted in this study, the CGMsFEM presented here is quite general and provides a generalized framework to design the coupling basis functions for fully coupled multiphysics multiscale problems under the ultimate load. It's worth noting that this regularization technique is not only effective within the framework of the GMsFEM but also applicable to other multiscale methods. Future work will focus on the analysis of the optimal relaxation coefficients for solving the local coupling spectral problems, although we find that they can be flexibly chosen without affecting the efficiency and accuracy of the CGMsFEM through numerical experiments.

\section*{Acknowledgement}
The authors were supported by the National Science Foundation of China (Nos. 12271409 and 12271408), the Natural Science Foundation of Shanghai (No. 21ZR1465800), the Interdisciplinary Project in Ocean Research of Tongji University and the Fundamental Research Funds for the Central Universities.

\FloatBarrier

\appendix

\section{The proof of Lemma \ref{Lemma:localinter}}
\label{sec:appendix:lemma}
Firstly, the local energy error of the interpolation operator $ \mathcal{I}_{\text{ms}} $ can be given as
\begin{equation}
  \begin{aligned}
    a^{K} & \left( \uu-\left(\mathcal{I}_{\text{ms}} \ww \right)_u , \uu-\left(\mathcal{I}_{\text{ms}} \ww \right)_u \right) +d^{K} \left(\theta-\left(\mathcal{I}_{\text{ms}} \ww \right)_{\theta} ,\theta-\left(\mathcal{I}_{\text{ms}} \ww \right)_{\theta} \right)           \\
          & = \int_{K} 2\mu \bepsilon \left(  \uu-\left(\mathcal{I}_{\text{ms}} \ww \right)_u \right):  \bepsilon \left(  \uu-\left(\mathcal{I}_{\text{ms}} \ww \right)_u \right) + \int_{K} \lambda \nabla \cdot \left( \uu-\left(\mathcal{I}_{\text{ms}} \ww \right)_u \right) \\
          & \quad \quad \nabla \cdot \left(  \uu-\left(\mathcal{I}_{\text{ms}} \ww \right)_u \right)
    + \int_{K} \kappa \nabla \left( \theta-\left(\mathcal{I}_{\text{ms}} \ww \right)_{\theta} \right) \cdot \nabla \left( \theta-\left(\mathcal{I}_{\text{ms}} \ww \right)_{\theta} \right)                                                                                      \\
          & \equiv I_1 + I_2 +I_3.
  \end{aligned}
  \label{eq:est_aKdK}
\end{equation}
Define $ E(u) : \mathbb{R}^d \rightarrow \mathbb{R}^{d \times d}$ with $ E(u) \left. \right|_{ij} = \frac{1}{2} \left(u_i+ u_j \right)$ for all $u \in V_u(\Omega) $. 
Let $\circ $ is Hadamard product. For any $K \in \mathcal{T}_H$, we have
\begin{equation} \nonumber
  \begin{aligned}
    \bepsilon \left( \uu-\left(\mathcal{I}_{\text{ms}} \ww \right)_u \right) & = \bepsilon \left( \sum_{y_i \in K} \chi_{i,u}^T I_{d} \left( \uu-\left(\mathcal{I}_{L_i}^{\omega_i} \ww \right)_u \right) \right) \\
                                                                             & = \sum_{y_i \in K} \left( E(\chi_{i,u}) \circ \bepsilon \left(\uu-\left(\mathcal{I}_{L_i}^{\omega_i} \ww \right)_u \right) + E\left(\uu-\left(\mathcal{I}_{L_i}^{\omega_i} \ww \right)_u \right) \circ \bepsilon \left( \chi_{i,u} \right) \right).
  \end{aligned}
\end{equation}
Then $I_1$ can be estimated as follows
\begin{equation}
  \begin{aligned}
    I_1 & \leq 2N_K \sum_{y_i \in K} \left( \int_{K} 2\mu E(\chi_{i,u}) \circ \bepsilon \left(\uu-\left(\mathcal{I}_{L_i}^{\omega_i} \ww \right)_u \right) :E(\chi_{i,u}) \circ \bepsilon \left(\uu-\left(\mathcal{I}_{L_i}^{\omega_i} \ww \right)_u \right) \right.                \\
        & \quad \quad \left. + \int_{K} 2\mu E\left(\uu-\left(\mathcal{I}_{L_i}^{\omega_i} \ww \right)_u \right) \circ \bepsilon \left( \chi_{i,u} \right): E\left(\uu-\left(\mathcal{I}_{L_i}^{\omega_i} \ww \right)_u \right) \circ \bepsilon \left( \chi_{i,u} \right)   \right) \\
        & \leq 2N_K \sum_{y_i \in K} \left(    \int_{K} 2\mu \bepsilon \left(\uu-\left(\mathcal{I}_{L_i}^{\omega_i} \ww \right)_u \right) : \bepsilon \left(\uu-\left(\mathcal{I}_{L_i}^{\omega_i} \ww \right)_u \right)          \right.                                           \\
        & \quad \quad \left. + \frac{C_1^2}{H^2}\int_{K} 2\mu   E\left(\uu-\left(\mathcal{I}_{L_i}^{\omega_i} \ww \right)_u \right):  E\left(\uu-\left(\mathcal{I}_{L_i}^{\omega_i} \ww \right)_u \right)       \right),                                                            \\
        & \leq 2N_K \sum_{y_i \in K} \left(    \int_{K} 2\mu \bepsilon \left(\uu-\left(\mathcal{I}_{L_i}^{\omega_i} \ww \right)_u \right) : \bepsilon \left(\uu-\left(\mathcal{I}_{L_i}^{\omega_i} \ww \right)_u \right)          \right.                                           \\
        & \quad \quad \left. + \frac{C_1^2}{H^2}2\int_{K} 2\mu   \left(\uu-\left(\mathcal{I}_{L_i}^{\omega_i} \ww \right)_u \right)\cdot  \left(\uu-\left(\mathcal{I}_{L_i}^{\omega_i} \ww \right)_u \right)      \right).
  \end{aligned}\label{eq:est_I1}
\end{equation}
Next, $I_2$ can be estimated for any $K \in \mathcal{T}_H$ as follows
\begin{equation} \nonumber
  \begin{aligned}
    \nabla \cdot \left( \uu-\left(\mathcal{I}_{\text{ms}} \ww \right)_u \right)  
                       & = \sum_{y_i \in K} \nabla \cdot \left(  \chi_{i,u}^T I_{d} \left( \uu-\left(\mathcal{I}_{L_i}^{\omega_i} \ww \right)_u \right) \right)                                                                                                  \\
    = \sum_{y_i \in K} & \left( \text{diag}(\chi_{i,u}) : \nabla \left(\uu-\left(\mathcal{I}_{L_i}^{\omega_i} \ww \right)_u \right) + \text{diag} \left(\uu-\left(\mathcal{I}_{L_i}^{\omega_i} \ww \right)_u \right) : \nabla \left( \chi_{i,u} \right) \right),
  \end{aligned}
\end{equation}
where $ \text{diag}\left( \cdot \right)$ represents the main diagonal elements of square matrix. Thus,
\begin{equation}
  \begin{aligned}
    I_2 & \leq 2N_K \sum_{y_i \in K} \left( \int_{K} \lambda \left( \text{diag}(\chi_{i,u}) : \nabla \left(\uu-\left(\mathcal{I}_{L_i}^{\omega_i} \ww \right)_u \right)  \right)^2 \right.
    \\& \quad \quad \left. +\left( \text{diag}(\uu-\left(\mathcal{I}_{L_i}^{\omega_i} \ww \right)_u) : \nabla \left(\chi_{i,u} \right)  \right)^2  \right)            \\
        & \leq 2N_K \sum_{y_i \in K} \left(    \int_{K} \lambda \left( \nabla \cdot \left(\uu-\left(\mathcal{I}_{L_i}^{\omega_i} \ww \right)_u \right) \right)^2 \right.
    \\& \quad \quad \left. + \frac{2C_1^2}{H^2}  \left(\uu-\left(\mathcal{I}_{L_i}^{\omega_i} \ww \right)_u \right) \cdot  \left(\uu-\left(\mathcal{I}_{L_i}^{\omega_i} \ww \right)_u \right) \right).
  \end{aligned}\label{eq:est_I2}
\end{equation}
Combining Eqs. (\ref{eq:est_I1}) with (\ref{eq:est_I2}), it follows that
\begin{equation}
  \begin{aligned}
    I_1+I_2 & \leq 2N_K \sum_{y_i \in K} \left(  \left\| \uu-\left(\mathcal{I}_{L_i}^{\omega_i} \ww \right)_u  \right\|_{a,K}^2  \right.                                                                                                         \\
            & \quad \quad \left. +\frac{2C_1^2}{H^2} \int_{K} \left( 2\mu+\lambda \right)   \left(\uu-\left(\mathcal{I}_{L_i}^{\omega_i} \ww \right)_u \right)\cdot  \left(\uu-\left(\mathcal{I}_{L_i}^{\omega_i} \ww \right)_u \right) \right).
  \end{aligned}\label{eq:est:I1I2}
\end{equation}
Estimating $I_3$ is similar, and we have
\begin{equation} \nonumber
  \begin{aligned}
    \nabla \left( \theta-\left(\mathcal{I}_{\text{ms}} \ww \right)_{\theta} \right) &                                                                                                                                                                                                                                                  
    = \sum_{y_i \in K} \nabla \left(  \chi_{i,\theta} \left( \theta-\left(\mathcal{I}_{L_i}^{\omega_i} \ww \right)_{\theta}  \right) \right)                                                                                                                                                                                           \\
                                                                                    & = \sum_{y_i \in K} \left( \chi_{i,\theta} \nabla \left(\theta-\left(\mathcal{I}_{L_i}^{\omega_i} \ww \right)_{\theta}  \right) + \left(\theta-\left(\mathcal{I}_{L_i}^{\omega_i} \ww \right)_{\theta} \right)  \nabla  \chi_{i,\theta}  \right).
  \end{aligned}
\end{equation}
Substituting it into $I_3$, we obtain
\begin{equation}
  \begin{aligned}
    I_3 & \leq 2N_K \sum_{y_i \in K} \left( \int_{K} \kappa \left( \chi_{i,\theta} \nabla \left(\theta-\left(\mathcal{I}_{L_i}^{\omega_i} \ww \right)_{\theta} \right)  \right) \cdot  \left( \chi_{i,\theta} \nabla \left(\theta-\left(\mathcal{I}_{L_i}^{\omega_i} \ww \right)_{\theta} \right)  \right)  \right. \\
        & \left.  \quad \quad +  \int_{K} \kappa \left(\left(\theta-\left(\mathcal{I}_{L_i}^{\omega_i} \ww \right)_{\theta} \right)  \nabla  \chi_{i,\theta} \right) \cdot  \left( \left(\theta-\left(\mathcal{I}_{L_i}^{\omega_i} \ww \right)_{\theta} \right)  \nabla  \chi_{i,\theta} \right)   \right)          \\
        & \leq 2N_K \sum_{y_i \in K} \left(    \int_{K} \kappa \nabla \left(\theta-\left(\mathcal{I}_{L_i}^{\omega_i} \ww \right)_{\theta} \right) \cdot   \nabla \left(\theta-\left(\mathcal{I}_{L_i}^{\omega_i} \ww \right)_{\theta} \right) \right.                                                              \\ 	& \left.  \quad \quad +\frac{2C_1^2}{H^2}  \int_{K} \kappa \left(\theta-\left(\mathcal{I}_{L_i}^{\omega_i} \ww \right)_{\theta} \right)^2 \right).
  \end{aligned}\label{eq:est:I3}
\end{equation}
Combining Eqs. (\ref{eq:est:I1I2}) and (\ref{eq:est:I3}), it follows that
\begin{equation}
  \begin{aligned}
    I_1+I_2+I_3 & \leq  2N_K \sum_{y_i \in K} \left(  \left\| \uu-\left(\mathcal{I}_{L_i}^{\omega_i} \ww \right)_u  \right\|_{a,K} ^2+  \left\| \theta-\left(\mathcal{I}_{L_i}^{\omega_i} \ww \right)_{\theta}  \right\|_{d,K} ^2+   \right.            \\
                & \quad \left. \frac{2C_1^2}{H^2}  \left(  \left\| \uu-\left(\mathcal{I}_{L_i}^{\omega_i} \ww \right)_u  \right\|_{L_a,K} +  \left\| \theta-\left(\mathcal{I}_{L_i}^{\omega_i} \ww \right)_{\theta}  \right\|_{L_d,K}  \right) \right).
  \end{aligned}\label{eq:est:I1I2I3}
\end{equation}
Moreover, define $ J=\left\| \uu-\left(\mathcal{I}_{L_i}^{\omega_i} \ww \right)_u  \right\|_{a,K} ^2+  \left\| \theta-\left(\mathcal{I}_{L_i}^{\omega_i} \ww \right)_{\theta}  \right\|_{d,K} ^2$, then
\begin{equation} \nonumber
  \begin{aligned}
    J & = \mathcal{A}^K\left( \ww-\mathcal{I}_{L_i}^{\omega_i} \ww,\ww-\mathcal{I}_{L_i}^{\omega_i} \ww\right) + \left( \gamma_1-\gamma_2 \right) b^K\left(\uu-\left(\mathcal{I}_{L_i}^{\omega_i} \ww \right)_u,\theta-\left(\mathcal{I}_{L_i}^{\omega_i} \ww \right)_{\theta} \right)                              \\
      & \leq \mathcal{A}^K\left( \ww-\mathcal{I}_{L_i}^{\omega_i} \ww,\ww-\mathcal{I}_{L_i}^{\omega_i} \ww\right) + | \gamma_1 - \gamma_2| C_0 \left\| \uu-\left(\mathcal{I}_{L_i}^{\omega_i} \ww \right)_u \right\|_{a,K}    \left\| \theta-\left(\mathcal{I}_{L_i}^{\omega_i} \ww \right)_{\theta} \right\|_{L_d,K} \\
      & \leq \mathcal{A}^K\left( \ww-\mathcal{I}_{L_i}^{\omega_i} \ww,\ww-\mathcal{I}_{L_i}^{\omega_i} \ww\right) + \frac{| \gamma_1 - \gamma_2| C_0}{2}                                                                                                                                                              \\
      & \quad \quad \left( \left\| \uu-\left(\mathcal{I}_{L_i}^{\omega_i} \ww \right)_u \right\|_{a,K}^2 +  \left\| \theta-\left(\mathcal{I}_{L_i}^{\omega_i} \ww \right)_{\theta} \right\|_{L_d,K}^2 \right)                                                                                             \\
      & \leq  \mathcal{A}^K\left( \ww-\mathcal{I}_{L_i}^{\omega_i} \ww,\ww-\mathcal{I}_{L_i}^{\omega_i} \ww\right) + \frac{| \gamma_1 - \gamma_2| C_0}{2} \left( J +  \left\| \theta-\left(\mathcal{I}_{L_i}^{\omega_i} \ww \right)_{\theta} \right\|_{L_d,K}^2 \right).
  \end{aligned}
\end{equation}
Then let $\ds C_2= \frac{2}{2-| \gamma_1 - \gamma_2| C_0},~C_3 = \frac{| \gamma_1 - \gamma_2| C_0}{2-| \gamma_1 - \gamma_2| C_0}$, and assume $\ds \frac{| \gamma_1 - \gamma_2| C_0}{2} < 1$, it follows that
\begin{equation}
  J \leq C_2 \mathcal{A}^K\left( \ww-\mathcal{I}_{L_i}^{\omega_i} \ww,\ww-\mathcal{I}_{L_i}^{\omega_i}  \ww\right) +C_3  \left\| \theta-\left(\mathcal{I}_{L_i}^{\omega_i} \ww \right)_{\theta} \right\|_{L_d,K}^2.
  \label{eq:est_J}
\end{equation}
Finally, combining Eqs. (\ref{eq:est_aKdK}), (\ref{eq:est:I1I2I3}) with (\ref{eq:est_J}), this proof is complete.

\section{The proof of Theorem \ref{theorem:energyerror}}
\label{sec:appendixB}
  By the definition of Riesz projection operator $\mathscr{R}_H$, and we define
  \begin{equation} \nonumber
    \eta_{Hu}^n =\mathscr{R}_{Hu} \left( \uu^n,\theta^n \right) - \uu^n_H, \quad  \eta_{H \theta}^n = \mathscr{R}_{H \theta} \left(\theta^n \right) - \theta^n_H.
  \end{equation}
  Combining Eq. (\ref{eq:discrevariation}), we have
  \begin{equation} \nonumber
    \begin{aligned}
      a\left(\eta_{Hu}^n, \vv_{uH} \right) & - b\left(\vv_{uH} , \eta_{H \theta}^n\right) + c\left( \eta_{H \theta}^n- \eta_{H \theta}^{n-1}, v_{\theta H} \right) + b \left( \eta_{Hu}^n- \eta_{Hu}^{n-1} ,v_{\theta H} \right) \\
                                           & + \tau_n d\left(  \eta_{H \theta}^n,v_{\theta H}\right)
      = \langle \tilde{\ff}^n-\tilde{\ff}^n_H, \vv_{uH}\rangle_a +  \tau_n \langle \tilde{g}^n-\tilde{g}^n_H, v_{\theta H}\rangle_d                                                                                              \\
                                           & + c\left( \delta_{H \theta}^n,v_{\theta H} \right) +b\left( \delta_{Hu}^n,v_{\theta H} \right), \quad  \forall \ \left(\vv_{uH},v_{\theta H}\right) \in V_{\text{cgm}},
    \end{aligned}
  \end{equation}
  where
  \begin{equation} \nonumber
    \begin{aligned}
      \delta_{H \theta}^n & =\mathscr{R}_{H \theta} \left(\theta^n \right)  - \mathscr{R}_{H \theta} \left(\theta^{n-1} \right) - \tau_n \partial_t \theta^n    \\
      \delta_{Hu}^n       & = \mathscr{R}_{Hu} \left(\uu^n,\theta^n \right)  - \mathscr{R}_{Hu} \left(\uu^{n-1},\theta^{n-1} \right) - \tau_n \partial_t \uu^n.
    \end{aligned}
  \end{equation}
  Define $ v_{uH}= \eta_{Hu}^n- \eta_{Hu}^{n-1} \in V_{uH}$ and $ v_{\theta H} = \eta_{H \theta}^n$, it follows
  \begin{equation}
    \begin{aligned}
      \frac{1}{2} \left\| \eta_{Hu}^n \right\|_a^2 + \frac{1}{2} \left\| \eta_{Hu}^n-\eta_{Hu}^{n-1} \right\|_a^2+\frac{1}{2} \left\| \eta_{H \theta}^n \right\|_c^2 +\frac{1}{2} \left\| \eta_{H \theta}^n-\eta_{H \theta}^{n-1} \right\|_c^2+\tau_n  \left\| \eta_{Hd \theta}^n \right\|_d^2 \\
      = \frac{1}{2} \left\| \eta_{Hu}^{n-1} \right\|_a^2 + \frac{1}{2} \left\| \eta_{H \theta}^{n-1} \right\|_c^2
      +  \langle \tilde{\ff}^n-\tilde{\ff}^n_H,\eta_{Hu}^n- \eta_{Hu}^{n-1} \rangle_a                                                                                                                                                                                                          \\
      + \tau_n \langle \tilde{g}^n-\tilde{g}^n_H,  \eta_{H \theta}^n \rangle_d
      + c\left( \delta_{H \theta}^n, \eta_{H \theta}^n \right) +b\left( \delta_{Hu}^n, \eta_{Hu}^n \right).
    \end{aligned}
  \end{equation}
  Similar to Eq. (\ref{eq:est:C4}), by use of $ \left\|  \eta_{H \theta}^n \right\|_c \leq C_4 \left\|  \eta_{H \theta}^n \right\|_d$, we have
  \begin{equation} \nonumber
    \begin{aligned}
      \frac{1}{2} \left\| \eta_{Hu}^n \right\|_a^2 + \frac{1}{2} \left\| \eta_{H \theta}^n \right\|_c^2 + \frac{1}{4} \tau_n  \left\| \eta_{H \theta}^n \right\|_d^2 \leq \frac{1}{2} \left\| \eta_{Hu}^{n-1} \right\|_a^2 + \frac{1}{2} \left\| \eta_{H \theta}^{n-1} \right\|_c^2 + C^n(\ff,g) \\
      + 4 C_4^2 \tau_n^{-1} \left\| \delta_{H \theta}^n  \right\|_c^2 + 4 C_4^2 \tau_n^{-1} C_0^2 \left\| \delta_{Hu}^n  \right\|_a^2.\label{eq:est_tn}
    \end{aligned}
  \end{equation}
  Based on the regularity assumption of problem and Eq. $ \left(\ref{eq:interpest_c}\right) $, we have
  \begin{equation}
    \begin{aligned}
      \left\| \delta_{H \theta}^n  \right\|_c & = \left\| - \int_{T_n} \left[ \partial_t \theta(s) -\mathscr{R}_{H\theta} \left( \partial_t \theta(s) \right) \right] ds - \int_{T_n} \left(s-t_{n-1} \partial_{tt}\theta(s)\right)dt\right\|_c,                                                                    \\
                                              & \leq \tau_n  C_4 C\left(H,\Lambda_{L+1}\right) \left\| \partial_t \ww  \right\|_{L^{\infty}\left(T_n, \interleave \cdot \interleave_{2,\Omega} \right)} +\frac{1}{2} \tau_n^2 \left\| \partial_{tt} \theta \right\|_{L^{\infty}\left(T_n, \| \cdot \|_{c} \right)}.
    \end{aligned}\label{eq:deltatheta}
  \end{equation}
  Similarly, using Eq. $ \left(\ref{eq:interpest_a}\right) $, we obtain
  \begin{equation}\label{eq:deltau}
    \left\| \delta_{Hu}^n  \right\|_a \leq \tau_n \text{max}\{ 1,C_4 C_0\} C\left(H,\Lambda_{L+1}\right) \left\| \left(\partial_t \ww \right)  \right\|_{L^{\infty}\left(T_n, \interleave \cdot \interleave_{2,\Omega} \right)} +\frac{1}{2} \tau_n^2 \left\| \partial_{tt} \uu \right\|_{L^{\infty}\left(T_n, \| \cdot \|_{a} \right)}.
  \end{equation}
  Combining Eq. (\ref{eq:deltatheta}) with Eq. (\ref{eq:deltau}), it follows
  \begin{equation} \nonumber
    4 C_4^2 \tau_n^{-1} \left\| \delta_{H \theta}^n  \right\|_c^2 + 4 C_4^2 \tau_n^{-1} C_0^2 \left\| \delta_{Hu}^n  \right\|_a^2 \leq \tau_n C^2\left(H,\Lambda_{L+1}\right) C_1^n(\ww) + \tau_n^3 C_2^n(\ww).
  \end{equation}
  Let $ u_H^0 =  \mathscr{R}_{Hu} \left( \uu^0,\theta^0 \right) $ and $ \theta_H^0 =  \mathscr{R}_{H\theta} \left(\theta^n \right) $, then $ \eta_{Hu}^0 =0 $ and $ \eta_{H\theta}^0 = 0 $, we have
  \begin{equation}
    \begin{aligned}
      \frac{1}{2} \left\| \eta_{Hu}^n \right\|_a^2 + \frac{1}{2} \left\| \eta_{H \theta}^n \right\|_c^2 +  \sum_{m=1}^{n} \frac{1}{4} \tau_m  \left\| \eta_{H \theta}^m \right\|_d^2 \leq \sum_{m=1}^{n} \left[ C^m\left(\ff,g\right) \right. \\
        \left. + \tau_m C^2\left(H,\Lambda_{L+1}\right) C_1^m(\ww) + \tau_m^3 C_2^m(\ww) \right].
    \end{aligned}
  \end{equation}
  Finally, the proof is complete by using the triangle inequality.

\bibliography{thermoref}

\end{document}